\documentclass[10pt]{amsart}

\usepackage{amssymb,amsmath,amscd,amsfonts,a4,psfrag,graphicx,
latexsym,epsfig,color}
\usepackage[active]{srcltx}
\usepackage[all]{xy}

\everymath{\displaystyle}

\newtheorem{teo}{Theorem}[section]
\newtheorem{lem}[teo]{Lemma}
\newtheorem{cor}[teo]{Corollary}

\newtheorem{prop}[teo]{Proposition}
\newtheorem{defi}[teo]{Definition}

\newtheorem{quess}[teo]{Questions}
\newtheorem{conj}[teo]{Conjecture}
\newtheorem{remark}[teo]{Remark}
\newtheorem{remarks}[teo]{Remarks}

\newcommand{\cvd}{\hfill$\Box$}

\newcommand{\mr}{\mathbb{R}}
\newcommand{\mc}{\mathbb{C}}
\newcommand{\mz}{\mathbb{Z}}

\newcommand{\R}{\mathbb{R}}
\newcommand{\C}{\mathbb{C}}
\newcommand{\Z}{\mathbb{Z}}

\newcommand{\Aa}{{\mathcal A}}
\newcommand{\Bb}{{\mathcal B}}

\newcommand{\vect}{\overrightarrow}

\newcommand{\bw}{{\bf w}}

\newcommand{\bv}{{\bf v}}

\newcommand{\Cc}{{\mathcal C}}
\newcommand{\Dd}{{\mathcal D}}

\newcommand{\Ff}{{\mathcal F}}

\newcommand{\Hh}{{\mathcal H}}

\newcommand{\Ll}{{\mathcal L}}

\newcommand{\Nn}{{\mathcal N}}

\newcommand{\Pp}{{\mathcal P}}
\newcommand{\Rr}{{\mathcal R}}
\newcommand{\Ss}{{\mathcal S}}
\newcommand{\Tt}{{\mathcal T}}
\newcommand{\Uu}{{\mathcal U}}
\newcommand{\Ww}{{\mathcal W}}

\newcommand{\Zz}{{\mathcal Z}}

\newcommand{\WW}{{\bf W}}

\newcommand{\hG}{{\mathfrak h}}

\newcommand{\sG}{{\mathfrak s}}

\newcommand{\ra}{\rightarrow}
\newcommand{\Dim}{{\it Proof.\ }}

\def\cvd{\hfill$\Box$}

\title[Non ambiguous structures on $3$-manifolds and quantum symmetry defects 
]{Non ambiguous structures on $3$-manifolds and quantum symmetry defects}

\author{St\'ephane Baseilhac, Riccardo Benedetti}

\begin{document}


\maketitle


\noindent $^1$ Institut Montpelli\'erain Alexander Grothendieck,
Universit\'e de Montpellier, Case Courrier 51, 34095 Montpellier Cedex
5, France (sbaseilh@univ-montp2.fr)

\smallskip

\noindent $^2$ Dipartimento di Matematica, Universit\`a di Pisa, Largo
Bruno Pontecorvo 5, 56127 Pisa, Italy (benedett@dm.unipi.it)


\begin{abstract} The state sums defining the quantum hyperbolic
  invariants (QHI) of hyperbolic oriented cusped $3$-manifolds can be
  split in a ``symmetrization" factor and a ``reduced'' state sum.  We
  show that these factors are invariants on their own, that we call
  ``symmetry defects'' and ``reduced QHI'', provided the manifolds are
  endowed with an additional ``non ambiguous structure'', a new type
  of combinatorial structure that we introduce in this paper. A
  suitably normalized version of the symmetry defects applies to
  compact $3$-manifolds endowed with $PSL_2(\mc)$-characters, beyond
  the case of cusped manifolds. Given a manifold $M$ with non empty
  boundary, we provide a partial ``holographic'' description of the
  non-ambiguous structures in terms of the intrinsic geometric
  topology of $\partial M$. Special instances of non ambiguous
  structures can be defined by means of taut triangulations, and the
  symmetry defects have a particularly nice behaviour on such ``taut
  structures''. Natural examples of taut structures are carried by any
  mapping torus with punctured fibre of negative Euler characteristic,
  or by sutured manifold hierarchies.  For a cusped hyperbolic
  $3$-manifold $M$ which fibres over $S^1$, we address the question of
  determining whether the fibrations over a same fibered face of the
  Thurston ball define the same taut structure. We describe a few
  examples in detail. In particular, they show that the symmetry
  defects or the reduced QHI can distinguish taut structures
  associated to different fibrations of $M$. To support the guess that
  all this is an instance of a general behaviour of state sum
  invariants of $3$-manifolds based on some theory of $6j$-symbols,
  finally we describe similar results about reduced Turaev-Viro
  invariants.
\end{abstract}
\medskip



\tableofcontents

\section{Introduction} 
A first aim of this paper is to produce refinements of the quantum
hyperbolic invariants (QHI) of hyperbolic cusped $3$-manifolds
(defined in \cite{AGT}, extending \cite{GT, AGT0}). But our thesis is
that such refinements are instances of a sort of ``universal''
phenomenon concerning the quantum invariants of $3$-manifolds based on
some theory of $6j$-symbols and defined by means of state sums over
triangulations. At least this holds for every example to our knowledge
(for instance the family of generalized Turaev-Viro invariants based
on the $6j$-symbols of any unimodular category, constructed in
\cite{Tu}). The arguments of these new invariants, that we call {\it
  reduced invariants}, are the same as those of the ``unreduced" ones,
but they apply to $3$-manifolds equipped with an additional {\it non
  ambiguous structure}, a new type of combinatorial structure on
compact oriented $3$-manifolds that we introduce and investigate,
pointing out in particular a strong relationship with the theory of
taut triangulations.  

\smallskip

Let us describe qualitatively this phenomenon.

\smallskip

$\bullet$ There is a ``background theory'', usually given by a
category of finite dimensional representations of a Hopf algebra; the
``$6j$-symbols" of the theory can be organized at first to produce a
kind of ``basic tensor'', that is, linear isomorphisms
$$\Bb: V_3\otimes V_1 \to V_2\otimes V_0$$
carried by an oriented $3$-simplex equipped with a further
``decoration'', say $d$, dictated by the background theory. Usually
the $V_j$ are complex vector spaces of the same dimension. A
$3$-simplex is a tetrahedron $\Delta$ with ordered vertices $v_0,
v_1,v_2,v_3$; this can be encoded by a system $b$ of edge orientations
(called a {\it (local) branching}) such that there are $j$-incoming edges at
the vertex $v_j$. The $2$-faces $F_j$ of $(\Delta,b)$ are ordered
accordingly with the opposite vertices $v_j$, and we associate to the
face $F_j$ the space $V_j$. This is the basic use of the branching to
associate the tensor $\Bb=\Bb(\Delta,b,d)$ to $(\Delta,b)$. But
depending on the background theory, there are further subtler uses
concerning the decoration $d$. For example, in quantum hyperbolic
theory an ingredient of the decoration is a triple of {\it shape parameters}, 
which are scalars in $\mc\setminus \{0,1\}$, associated to the triple of couples of
opposite edges. These naturally occur with a {\it cyclic ordering}
depending on the orientation of $\Delta$. The branching $b$ is used to
select one {\it linear} ordering compatible with the cyclic one.
\smallskip

$\bullet$ Given a triangulation $T$ of some compact oriented
$3$-manifold $M$, we give each oriented tetrahedron $\Delta$ of $T$ a
branching and a decoration, so that the system formed by such local
data verifies certain {\it global constraints} (also dictated by the
background theory). For example, for what concerns the local branchings we could
require that the edge orientations of the tetrahedra match under the $2$-face pairings in $T$, to produce a {\it global
  branching} $b$ of $T$. This is equivalent to promote $T$ to be a
$\Delta$-complex (a generalized simplicial complex), a kind of object
familiar in algebraic topology (see \cite{HATCHER}). But it turns out that this is too
demanding.  On another hand, every oriented $3$-simplex $(\Delta,b)$
carries a weaker structure, that is, a system of transverse
co-orientations of its $2$-faces such that two co-orientations are
incoming and two are outgoing.  We call it the local {\it
  pre-branching} $(\Delta,\omega_b)$ induced by the branching $b$.  We
say that a system $\tilde b$ of local branchings of the tetrahedra of $T$ is a {\it weak
  branching} if their $2$-face co-orientations match under the $2$-face pairings, and we call such a 
  system $\omega_{\tilde b}$ of co-orientations a global {\it
  pre-branching} of $T$. Pre-branchings occur, for instance, in the
definition of taut triangulations \cite{L}. But one eventually
realizes that the notion of pre-branching is also the most fundamental
global enhancement of $3$-manifold triangulations in order to deal
with ``quantum state sums'' (for any background theory). Given a
weakly branched triangulation $(T,\tilde b)$ endowed with a global
decoration, say $\Dd$, we get a {\it tensor network} by associating
the above basic tensor $\Bb(\Delta,b,d)$ to every decorated
$3$-simplex $(\Delta,b,d)$ of $(T,\tilde b,\Dd)$. The state sum
$\Bb(T,\tilde b, \Dd)$ is by definition the {\it total contraction} of
this network. We call it a {\it reduced} (or basic) state sum of the
theory.

\smallskip

$\bullet$ In order to obtain invariants of the $3$-manifold $M$
(possibly equipped with additional structures, as it happens in quantum
hyperbolic theory) by means of such state
sums, we have to mod out the arbitrary choices we have made, in
particular the choice of the weak branching $\tilde b$. This gives
rise to a more or less delicate procedure of {\it symmetrization} of
the basic tensors, producing tensors of the same type, say
$\Ss(\Delta,b,d)$, such that the corresponding state
sums (ie. network total contraction) $\Ss(T,\tilde b, \Dd)$ have the required
invariance properties.

\smallskip

$\bullet$ One might wonder anyway whether the reduced state sums
$\Bb(T,\tilde b, \Dd)$ define some kind of $3$-dimensional invariant. 
For example it holds that, keeping the decoration $\Dd$ fixed,
the value of $\Bb(T,\tilde b, \Dd)$ depends only on the underlying
pre-branching $\omega_{\tilde b}$, not on $\tilde b$ itself. Moreover,
the state sums verify a highly non trivial system of functional
identities which apparently is formally the same for every background
theory; this system corresponds to a restricted system
of ``moves'' on pre-branched triangulations, called {\it non ambiguous
  transits}.  Then the notion of {\it non ambiguous structure} on $M$
arises as an equivalence class of pre-branched triangulations of
$M$ up to non ambiguous transits.

\medskip

In order to substantiate our thesis, in the present paper we spell 
out the case of the quantum hyperbolic state sums, and in the
Appendix, the most fundamental prototype of this business, the
Turaev-Viro state sums \cite{TV}.

\medskip

Having in mind this strong motivation (at least in our opinion), the theory
of non ambiguous structures can be introduced and developed by itself,
without any reference to any specific $6j$-symbols theory. It is
remarkable nevertheless that in doing it, some issues of the quantum
hyperbolic machinery emerge. We will develop two instances of
non ambiguous structures; the most important is based on {\it
 ideal} triangulations of $3$-manifolds that are the interior of
compact connected oriented $3$-manifolds $M$ with non empty boundary
(Sections 2 to 6); the other is based on relative ``distinguished
triangulations'' $(T,H)$ of $(M,L)$, where $M$ is a compact closed
oriented $3$-manifold and $L$ is a non empty link in $M$ (Section \ref{ML}).

\medskip

In the rest of this introduction we describe more features of the non ambiguous 
structures (Sections \ref{on-id-na} and \ref{on-rel-na}) and our results 
about reduced quantum hyperbolic invariants (Section \ref{on-red-qhi}).

\subsection{On ideal non ambiguous structures}\label{on-id-na}
After some generalities on $3$-manifold triangulations and the
different notions of ``branchings'' (Section \ref{TRIANG}), the
combinatorial definition of ideal non ambiguous structures is given in
Section \ref{NAIS}. We prefer to consider an ideal triangulation $T$
of Int$(M)$ as a triangulation of the compact space $\hat M$ obtained
by adding a point at infinity at each end of Int$(M)$, requiring that
the set of vertices of $T$ coincides with the set of these added
points. These ``naked'' ideal triangulations are considered up to the
equivalence relation generated by the well known $2\leftrightarrow 3$
(MP) and $0\leftrightarrow 2$ (lune) moves. Dealing with ideal
pre-branched triangulations $(T,\omega)$, there are natural notions of
$pb$-{\it transits} that enhance the above naked moves.  Roughly
speaking, a $pb$-transit $(T,\omega)\rightarrow (T',\omega')$ is ``non
ambiguous'' if both the transit and its inverse $(T,\omega)\leftarrow
(T',\omega')$ are the unique $pb$-enhancements of the underlying naked
move with the given initial configuration. Once the combinatorial definition has 
been established, our effort is to point out some {\it intrinsic geometric topological 
content} and natural families of (ideal)
 non ambiguous structures.

\smallskip

In Section \ref{seccharges} we remark that the non ambiguous
structures on $M$ have an intrinsic cohomological content, strictly
related to the theory of charges and ``cohomological weights'' 
 which play an important role in the definition of the QHI for
hyperbolic cusped manifolds.

\smallskip

In Section \ref{2D} we develop a partial, though rather illuminating,
``holographic" approach to the non ambiguous structures on $M$, based
on a suitably defined restriction of the $3D$ structures on $\partial
M$.  In particular, we discover that these $2D$ structures carry
certain {\it singular combings} on $\partial M$ (defined in intrinsic
geometric topological terms) that are eventually invariants of the
$3D$ structures. A conjectural holographic classification of ideal non ambiguous structures is given in Conjecture \ref{Ricconj}. 

\smallskip
 
When $\partial M$ is a collection of tori, as for hyperbolic cusped
$3$-manifolds, we easily realize that the (possibly empty) set of {\it
  taut} pre-branched ideal triangulations $(T,\omega)$ of $\hat M$ in
the sense of \cite{L} is closed with respect to the non ambiguous
transits. We call {\it taut structures} the corresponding non
ambiguous structures.  If $M$ is realized as a mapping torus $M_\psi$
($\psi$ being an automorphism of a punctured surface $\Sigma$ with
$\chi(\Sigma)<0$, considered up to isotopy), or more generally if $M$
carries a {\it sutured manifold hierarchy} $\mathcal{Z}$ (which exists
for example when $M$ is a hyperbolic cusped manifold), then one finds
in \cite{L} a procedure to construct taut triangulations $(T,\omega)$
of $\hat M$ which depend on $M_\psi$ or $\mathcal{Z}$ and also on
other arbitrary choices. For example
in the case of $M_\psi$ these choices are an ideal triangulation of a fiber,  say
$R$, and a sequence of elementary diagonal exchanges
(``flips'') that connects $R$ with $\psi(R)$; with these data one
constructs a so called taut ``layered'' ideal triangulation of $\hat
M_\psi$.  One eventually realizes that these further choices are immaterial up to
non ambiguous transits, and that we have the following
result (see Section \ref{TAUT}). Let us call ``Thurston ball of $M$''
the unit ball of the Thurston norm on $H_2(M,\partial M;\R)$.

\begin{teo}\label{existNA} (1) Every mapping torus $M_\psi$ with
  punctured fibre of negative Euler characteristic carries a natural
  taut structure $\sG_{\psi}$, represented by any {\rm layered}
  triangulation, constructed by means of any ideal triangulation of a
  fiber.
  
  (2) Assume that $M$ fibers over $S^1$. Then, any two fibrations of
  $M$ such that the corresponding mapping tori $M_\psi$ and $M_\phi$
  satisfy $\sG_\phi = \sG_\psi$ lie in the cone over a same face of
  the Thurston ball of $M$. Moreover, mapping tori corresponding to
  fibrations lying on a same ray from the origin of $H_2(M,\partial
  M;\R)$ satisfy $\sG_\phi = \sG_\psi$.
     
  (3) More generally, every compact oriented $3$-manifold $M$ with non
  empty boundary equipped with a sutured manifold hierarchy
  $\mathcal{Z}$ carries a natural taut structure $\sG_{\mathcal{Z}}$.
\end{teo}

In the simpler case when $M$ fibres over $S^1$, it can happen that
different (non multiple) fibrations carry the same non ambiguous structure.  An
interesting case is provided by the following Theorem which easily
follows from a result of Agol \cite{A,A2}. For completeness we will
discuss the proof in Section \ref{TAUT}.

\begin{teo}\label{Agol} 
  To every couple $(W,\Ff)$, where $W$ is a complete hyperbolic
  $3$-manifold of finite volume that fibres over $S^1$ and $\Ff$ is a
  fibred face of the Thurston ball of $W$, one can associate in a
  canonical way a couple $(M,\Ff')$, where $M$ is cusped manifold
  (obtained by removing a suitable link from $W$) and $\Ff'$ is a
  fibred face of the Thurston ball of $M$, such that all fibrations
  $M_\psi$ of $M$ in the cone over $\Ff'$ define the same taut
  structure. Hence the taut structure $\sG_{\Ff'} := \sG_{\psi}$ is
  well defined.
\end{teo}

In general, for any oriented $3$-manifold $M$ bounded by tori and
which fibers over $S^1$, Theorem \ref{existNA} (2) implies that for every 
fibred face $\Ff$ of the Thurston ball of $M$, there is a
well defined map $\sG_\Ff$ which associate to any rational point $p\in
{\rm Int}(\Ff)$ the taut structure $$\sG_\Ff(p):=\sG_\psi$$ where
$\psi$ is the monodromy of any fibration of $M$ in the ray spanned by
$p$.  A very attractive and probably demanding problem is to study
this map in general.  For example, we can state

\begin{quess}{\rm (1) Is $\sG_\Ff$ always constant? (2) Otherwise, is the image
    of $\sG_\Ff$ always finite?}
\end{quess}

\subsection{On relative non ambiguous structures}\label{on-rel-na}
Dealing with not necessarily ideal triangulations of $\hat M$, for
example when $M$ is closed, we must complete the naked triangulation
moves with the so called $0\leftrightarrow 2$ {\it bubble} move, which
modifies the number of vertices (see Section \ref{ML}). Looking at the
associated $pb$-transits we realize that none is ``non ambiguous'' in
a strict sense.  We need some further input to select one. We can do
it in the framework of relative {\it distinguished triangulations}
$(T,H)$ of $(M,L)$, qualified by the fact that $H$ is a Hamiltonian
subcomplex of the $1$-skeleton of $T$ isotopic to the link $L$. These
are considered up to relative ``distinguished'' versions of the naked
moves (bubble move included). This kind of triangulation has been already
used to define QHI for pairs $(M,L)$.  In Section \ref{ML} we
develop the relative non ambiguous structures somewhat in parallel to
what we have done in the case of ideal triangulations. In particular
we will introduce the notion of {\it relative taut structure} and
indicate some procedures to construct examples.

\medskip

Let us ouline now some specific features of the {\it reduced} quantum
hyperbolic invariants.

\subsection{On reduced quantum hyperbolic invariants}
\label{on-red-qhi} 
Although the QHI can be defined in more general situations (see
\cite{AGT0}), in this paper we focus on two main instances of compact
oriented $3$-manifolds: {\it cusped manifolds} $M$ such that the non
empty boundary is made by tori and  the
interior has a finite volume complete hyperbolic structure (see \cite{GT,AGT0,AGT}); pairs
$(M,L)$, where $M$ is closed, $L$ is a non empty link in $M$, and $M$
is equipped with a $PSL_2(\mc)$-character (see \cite{Top, GT,
  LINK}). In such situations, for every odd integer $N\geq 3$ the QHI
of $M$ or $(M,L)$ at {\it level $N$} is a complex number defined up to
multiplication by $2N$th-roots of unity. Let us assume for a while
that $M$ is a cusped manifold. Then its QHI depend on a choice of
conjugacy class $\rho$ of representations of $\pi_1(M)$ in
$PSL(2,\mc)$, and two pairs of so called {\it bulk} and {\it boundary
  weights} $h:=(h_f,h_c)$ and $k:=(k_f,k_c)$ given by cohomology
classes
\begin{equation}\label{notationweights}
(h_f,h_c)\in H^1(M;\Z/2\Z)^2\ ,\ (k_f,k_c)\in H^1(\partial M;\C)
\times H^1(\partial M; \Z)
\end{equation}
satisfying certain natural compatibility conditions. In particular,
$i^*(h_c)=k_c$ mod$(2)$, where the map $i:\partial M \to M$ is the inclusion;
$k_f$ encodes a sort of ``logarithm'' of the class in $H^1(\partial
M;\mc/2\pi i \mz)$ (with multiplicative coefficients) 
defined by the restriction of $\rho$ on $\partial
M$. Also, we assume that when Int($M$) has one cusp, $\rho$ varies in the 
irreducible component $X$ of the variety of $PSL(2,\mc)$-characters of $M$ 
containing the character of the discrete faithful holonomy $\rho_{hyp}$; if there are several cusps, 
$X$ has to be replaced by its so called {\it eigenvalue} subvariety (see \cite{KT} for this notion).

\begin{remark}\label{manycusp}{\rm In \cite{AGT} we treated only the case of 
one-cusped manifolds because in this case we could develop a rigidity argument for systems 
of shape parameters $w$ with holonomies in $X$, based on a result of N. Dunfield in \cite{Dun} (
see below for the notion of shape parameters). This result has been extended in \cite{KT} to the case of an
arbitrary number of cusps, and we can adapt the rigidity argument of \cite{AGT} to this setup as 
well by using the eigenvalue variety mentioned above instead of $X$. If a reader prefers to dispose 
of a detailed reference like \cite{AGT}, she/he can restrict to
one-cusped manifolds, without substantially effecting the discussion of the
present paper.}
\end{remark}

For every odd integer $N\geq 3$, the QHI $\Hh_N(M,\rho,h,k)$ is
computed by state sums $\Hh_N(\Tt)$ over so called {\it QH
  triangulations} $\Tt=(T,\tilde b,w,f,c)$, which are ideal weakly
branched triangulations $T$ of $\hat M$ ``decorated'' with a heavy
apparatus $\Dd=(w,f,c)$ of combinatorial structures encoding $\rho$
and $(h,k)$. In fact $w$ encodes a system of shape parameters on the
abstract tetrahedra $\Delta$ of $T$ verifying the Thurston compatibility
condition around every edge of $T$ (that is $w$ determines a point in the
``gluing variety'' carried by the triangulation $T$); $f$ and $c$ are integer valued
labellings of the couples of opposite edges of every $\Delta$, called
{\it flattening} and {\it charge} respectively, which contribute to
determine a system $\bw$ of $N$-th roots of the shape parameters $w$
(verifying suitable local and global constraints, at every $\Delta$ and
around every edge of $T$). Referring to the qualitative picture
depicted at the beginning of this Introduction, in the present
situation all spaces $V_j=\C^N$, the basic tensors have an
explicit matrix form, are called {\it basic matrix dilogarithms}, and
are denoted by $\Ll_N(\Delta,b,d)$. They are derived from the
$6j$-symbols of the cyclic representations of a Borel subalgebra of
$U_q(sl_2)$ ($q$ being a $N$-th root of $1$), firstly derived in the
seminal Kashaev's paper \cite{K}. The ``symmetrized''
tensors have the same type, are called {\it matrix dilogarithms}, and
are denoted by $\Rr_N(\Delta,b,d)$. It is a specific feature of the
quantum hyperbolic setting that every symetrized tensor is equal to 
the corresponding basic one up to a scalar factor, that is
\begin{equation}\label{basicvsnon}
\Rr_N(\Delta,b,d):= \alpha_N(\Delta,b,d)\Ll_N(\Delta,b,d)
\end{equation}
where $\alpha_N(\Delta)$ is a scalar called the {\it local
  symmetrization factor} of $(\Delta,b,d)$ (see Section \ref{SYMD}); hence the state
sums can be factorized as
\begin{equation}\label{facto}
\Hh_N(\Tt) := \alpha_N(\Tt)\Hh_N^{red}(\Tt)
\end{equation}  
where
\begin{equation}\label{localfactor}
\alpha_N(\Tt):= \prod_{\Delta\in \Tt^{(3)}}\alpha_N(\Delta,b,d)
\end{equation}
is called the {\it global symmetrization factor}, while $\Hh_N^{red}(\Tt)$
is the {\it reduced} state sum, which involves only the basic 
matrix dilogarithms. Note that in general (for instance in the case
of Turaev-Viro state sums considered in the Appendix) there is not
such a simple factorization.

\begin{remarks}\label{varie}
{\rm (1) The definition of the QHI by means of {\it weakly}
branched triangulations is an achievement of \cite{AGT}. In our
previous papers we used more demanding {\it branched} triangulations. 
Given a weakly branched triangulation $(T,\tilde b)$ the $2$-face
pairings which produce $T$ from its set of ``abstract tetrahedra'' are encoded 
by colorings of the $2$-faces $F$ of $T$ by colors $s(F)\in \mz/3$ (we will show it in practice in the examples of Section \ref{eightsec} and \ref{sistersec}). The tensor network whose contraction
is the state sum $\Hh_N(\Tt)$ includes ``face'' tensors $Q^{s(F)}$, $Q$ being an automorphism 
of $\C^N$ with an explicit matrix form. If $\tilde b$ is a genuine branching, then such face tensors are immaterial 
(all colors $s(F)=0$). So, keeping the same notation, we stipulate that the basic matrix dilogarithms incorporate the face 
tensors, in the sense that we add to the decoration $d$ of $(\Delta,b)$ the $\Z/3$-colors of the two $2$-faces  
with outgoing transverse co-orientation, with respect to the associated pre-branching $\omega_b$, 
and contract $\Ll_N(\Delta,b,d)$ with the face tensors associated to these $2$-faces.
\smallskip

(2) We stress that $QH$ triangulations over ideal triangulations $T$ of $\hat M$ make sense and may exist 
beyond the case of cusped manifolds, that is, assuming just that $M$ has a 
non empty boundary made by tori. Also in this general case, a $QH$ triangulation 
$\Tt=(T,\tilde b,w,f,c)$ encodes a $PSL(2,\C)$-valued character $\rho$ of $\pi_1(M)$, and 
a system of weights $h:=(h_f,h_c)$ and $k:=(k_f,k_c)$. There are no a priori restrictions on $\rho$.   
Hence (reduced) state sums and symmetrization factors are defined as well.}
\end{remarks}

In Section \ref{SYMD} we prove the following results. All terms are defined
precisely in Section \ref{secondproofs}. 
\begin{teo}\label{Intro-inv} Let $\Tt=(T,\tilde
b,w,f,c)$ be a QH triangulation encoding a tuple
$(M,\rho,h,k)$, where $M$ is any compact connected oriented $3$-manifold with
non empty boundary made by tori (according to Remark \ref{varie} (2) above). 
Denote by $\omega$ the pre-branching underlying $(T,\tilde b)$. Then we have:
\smallskip

(1) The value of $\alpha_N(\Tt)$ does not depend on the choice of
  $\tilde b$ among the weak branchings compatible with $\omega$, and
  it does not depend on the choice of $c$ among the charges encoding
  $(h_c,k_c)$ up to multiplication by $4$-th roots of $1$. On another
  hand, in general it varies with the flattening $f$ by a $4N$-th root
  of $1$.
  \smallskip
    
    (2) Let $\Tt$ and $\Tt'$ be two QH triangulations such that the
  underlying pre-branchings $\omega$ and $\omega'$ represent a same
  non ambiguous structure on $M$. If $\Tt$ and $\Tt'$ are connected by
  a sequence of QH transits lifting a sequence of non-ambiguous
  transits between $(T,\omega)$ and $(T',\omega')$, then
  $\alpha_N(\Tt)=\alpha_N(\Tt')$.

  (3) The conclusions of (1) and (2) hold true up to multiplication
  by $4N$-th roots of $1$ by replacing $\alpha_N$ with the reduced
  state sums $\Hh_N^{red}$. Moreover, $\Hh_N^{red}(\Tt)$ does not
  depend on the choice of bulk weight $h=(h_f,h_c)$, and as a function
  of $k_f$ and $k_c$ it depends only on $k_f-\pi i k_c$ mod$(\pi i
  N)$.
\end{teo}

\begin{defi}\label{fused}{\rm We call the above class $\kappa:= k_f-\pi i k_c\in H^1(\partial M;\mc/\pi i N\mz)$ 
a {\it fused weight}. }
\end{defi} 
Let us restrict now to cusped manifolds $M$. Recall that if $M$ has a single cusp we denote by $X$ the irreducible component of the variety of $PSL(2,\mc)$-characters of $M$ containing the character of the discrete faithful holonomy $\rho_{hyp}$. If there are several cusps, we consider instead the eigenvalue subvariety of $X$ (see Remark \ref{manycusp}). 
Fix a further notation: For any integer $n$, denote by $\mu_{n}$ the group of $n$-th roots of $1$ acting on $\mc$ by multiplication.
Then we will deduce from Theorem \ref{Intro-inv} (again all terms are defined in Section \ref{secondproofs}):
\begin{cor}\label{corinv2}  (1) For every non ambiguous structure $\sG$ on $M$, the value of $\alpha_N(\Tt)$ on any rich QH triangulation $\Tt$ encoding
  $(M,\rho,h,k)$ and $\sG$ does not depend on the choice of $\Tt$ and
  $(h_c,h_f,k_f)$ up to multiplication by $4N$-th roots of $1$. Also,
  the reduced state sums $\Hh_N^{red}(\Tt)$ do not depend on the
  choice of $\Tt$ and $(h_c,h_f)$, and both define invariants
  $\alpha_N(M,\rho,k_c;\sG)$ and $\Hh_N^{red}(M,\rho,\kappa;\sG)$,
  where $\kappa:= k_f-\pi i k_c$ is the fused weight as above.
\smallskip

  (2) Assume that $M$ has only one cusp. Then there
  exists a determined $(\mz/N\mz)^2$-covering space $\tilde{X}_N$ of
  $X$ such that, by fixing $\sG$ and $k_c$, and varying $\rho$ in $X$
  and $\kappa$ among the fused weights compatible with $\rho$,
  $\alpha_N(M,\rho,k_c;\sG)$ defines a function on $X$ that lifts to a
  rational function $\alpha^{k_c,\sG}: \tilde{X}_N\rightarrow
  \mc/\mu_{4N}$, and $\Hh_N^{red}(M,\rho,\kappa;\sG)$ defines a
  rational function $\Hh_N^{red,\sG}:\tilde{X}_N \rightarrow
  \mc/\mu_{4N}$. Similar results hold true when $M$ has several cusps by replacing $X$ with the eigenvalue variety.
  \end{cor}

\smallskip

We call $\alpha_N(M,\rho,k_c;\sG)$ and $\Hh_N^{red}(M,\rho,\kappa;\sG)$
the {\it symmetry defects} and {\it reduced QHI} respectively. They
have the same ability to distinguish different non ambiguous
structures $\sG$, by the formula \eqref{facto} and the fact that the (unreduced) QHI do
not depend on $\sG$. The symmetry defects involve only products of
simple scalars, and so they are much simpler to compute than the
reduced QHI. This is useful in studying non ambiguous structures.

\begin{remark}\label{intert}{\rm There should be strong connections, 
    that deserve to be fully understood in future investigations,
    between the reduced QHI of fibred cusped manifolds and the
    intertwiners of local representations of the quantum
    Teichm\"uller spaces, introduced in \cite{BBL}.}
\end{remark}


Beside the symmetrization factors $\alpha_N(\Tt)$, it is also
meaningful to consider {\it normalized symmetrization factors}
associated to a pair of ``base" $c$-weights $(h_c^0,k_c^0)$. They are
defined by
\begin{equation}\label{normsymdef}
\alpha_{N,c_0}(\Tt):=\alpha_N(\Tt)/\alpha_N(\Tt_{c_0})
\end{equation}
where $\Tt_{c_0}$ is obtained from $\Tt$ by replacing the charge $c$
(encoding the weights $(h_c,k_c)$) with any charge $c_0$, encoding a weight
$(h_c^0,k_c^0)$. By Theorem \ref{Intro-inv} (2) we have clearly
$\alpha_{N,c_0}(\Tt) = \alpha_{N,c_0}(\Tt')$. On another hand,
$\alpha_{N,c_0}(\Tt)$ has better invariance properties with respect to
$(w,f,c)$, so that Theorem \ref{Intro-inv} (1) becomes:

\begin{teo}\label{Intro-invbis} The value of $\alpha_{N,c_0}(\Tt)$
  does not depend on the choice of $\tilde b$ among the weak
  branchings compatible with $\omega$, and it does not depend on the
  choice of tuple $(w,c,f)$ and charge $c_0$ encoding $(\rho,h,k)$ and
  $(h_c^0,k_c^0)$ up to multiplication by $4$-th roots of
  $1$. Moreover, as a function of $(\rho,h,k)$ and $(h_c^0,k_c^0)$ it
  depends only on $k_f$ and $k_c-k_c^0$.
\end{teo}
Then we will get the following generalization of Corollary \ref{corinv2}. 
Note that its range goes beyond the case of cusped manifolds,
according to Remark \ref{varie} (2).
\begin{cor}\label{corinv2bis} Let $M$ be an arbitrary compact oriented
  $3$-manifold $M$ such that $\partial M$ is a collection of tori and
  $\rho$ can be represented on the gluing variety of an ideal
  triangulation of $\hat M$. Let $(h_c^0,k_c^0)$ and $\sG$ be any
  $c$-weights and non ambiguous structure on $M$. Then, for any
  weights $(h,k)$ of $(M,\rho)$, the value of $\alpha_{N,c_0}(\Tt)$ is
  independent of the choice of $c_0$ among the charges encoding
  $(h_c^0,k_c^0)$, and independent of the choice of $\Tt$ among the QH
  triangulations encoding $(M,\rho,h,k)$ and $\sG$, up to
  multiplication by $4$-th roots of $1$. As a function of $\rho$,
  $(h,k)$ and $(h_c^0,k_c^0)$ it depends only on $k_f$ and
  $k_c-k_{c}^0$, and hence defines an invariant
  $\alpha_{N,k_c^0}(M,k_f,k_c;\sG)$. If $M$ is a cusped manifold, it extends to a rational
  function $\alpha^{k_c-k_c^0,\sG}: H_1(\partial M;\mc^*)\rightarrow
  \mc/\mu_{4}$.
\end{cor}

We call $\alpha_{N,k_c^0}(M,k_f,k_c;\sG)$ a {\it normalized symmetry
  defect}. Perhaps its residual ambiguity by $4$-th roots of $1$ is not sharp, 
  but this is not the point of the
paper; a similar issue was solved for the QHI sign ambiguity in
\cite{AGT}, Section 8.
\smallskip

Clearly $\alpha_{N,k_c^0}(M,k_c,k_f;\sG)=1$ whenever
$k_c=k_c^0$. A ``universal'' natural choice can be $k_c^0=0$.
Another natural choice is possible for taut structures. Every taut
triangulation $(T,\omega)$ carries a ``tautological'' charge $c_0$.
The very definition of $\alpha_N(\Tt)$ implies that
$\alpha_{N,c_0}(\Tt) = \alpha_{N}(\Tt)$ for {\it any} QH triangulation
$\Tt=(T,\tilde b,w,f,c)$ such that $(T,\omega_{\tilde b})$ is a taut
triangulation and $c_0$ is the charge tautologically carried by
$(T,\omega_b)$.  Then Corollary \ref{corinv2bis} implies:
\begin{cor}\label{corinv3} For any taut structure $\sG$ the symmetry
  defect $\alpha_N(M,\rho,k_c;\sG)$ depends only on the restriction of
  $\rho$ to $\partial M$ and lifts to an invariant
  $\alpha_N(M,k_f,k_c;\sG)$ depending on $k_f$ and well-defined up to
  multiplication by $4$-th roots of $1$. It is defined for any $M$ and
  $\rho$ as in Corollary \ref{corinv2bis}, and satisfies
  $\alpha_{N}(M,k_f,k_c;\sG) = \alpha_{N,k_c^0}(M,k_f,k_c;\sG)$ where
  $k_c^0$ is the boundary $c$-weight tautologically carried by $\sG$.
\end{cor}

In a sense this nice behaviour of the symmetry defect indicates that
taut structures are the most natural non ambiguous structures.

\smallskip

A few words about the proofs of these results. We adopt again a kind of 
``holographic'' approach (see Section \ref{on-id-na}).
Roughly, every QH triangulation $\Tt$ of $M$ induces a ``$2D$ QH
triangulation'' $\partial \Tt$ of $\partial M$ where we can compute a
scalar $\alpha^0_N(\partial \Tt)$ such that $\alpha_N(\Tt)^4 =
\alpha^0_N(\partial \Tt)$. The proof of Theorem \ref{Intro-inv} (1)
deals with $\alpha^0_N(\partial \Tt)$, for which it turns out that, at
a fixed $w$, only the boundary weight $k_c$ is relevant. Moreover,
under the normalization \eqref{normsymdef} and up to the weaker
ambiguity by $4$-th roots of $1$, the same argument shows that at a
fixed $k_f$ also the choice of $w$ and $f$ is immaterial. This proves
Theorem \ref{Intro-invbis} and Corollary \ref{corinv3}. The
conclusions of Corollary \ref{corinv3} are not true in general for the
(non normalized) symmetry defect of arbitrary non ambiguous
structures, as can be seen by explicit computations (eg. when $M$ is
the sister of the figure eight knot complement, see Section \ref{EXAMPLES}). 
The first claim of Theorem \ref{Intro-inv} (3) follows
immediately from (1), (2) and the factorization formula \eqref{facto};
the second claim is easy. In order to deduce Corollary \ref{corinv2}
we combine Theorem \ref{Intro-inv} with a {\it rigidity} argument
about the shape parameters $w$, that we had already used in the
invariance proof of the QHI in \cite{GT,AGT}. It is based on a result of N. Dunfield (\cite{Dun}), 
extended in \cite{KT} as mentioned in Remark \ref{manycusp}.
On another hand, we will see that Corollary \ref{corinv2bis} follows almost 
immediately from the statement analogous to Theorem \ref{Intro-inv} (2) for the
normalized symmetrization factors, with no need of any rigidity
argument.
\smallskip

{\bf On relative reduced QHI.}  The symmetry defects and reduced QHI
of pair $(M,L)$ equipped with relative non ambiguous structures are
treated in Section \ref{MLdefect}. In the case of pairs $(M,L)$, the
(unreduced) QHI depend on a arbitrary conjugacy class $\rho$ of
representations of $\pi_1(M)$ in $PSL(2,\mc)$, and the weights $h$ and
$k$ reduce to the ``bulk'' weights $h_f$ and $h_c$.  Differently from
the case of cusped manifolds, we will see that in general the
non normalized symmetry defects (hence the reduced QHI) are ill-defined,
while the normalized ones are well defined but trivial. For relative
taut structures the reduced QHI are well defined but are not able to
distinguish them.  Precisely we have:
\begin{prop} \label{invcaseML} (1) For every relative non ambiguous structure
  $\sG$, every normalized symmetry defect
  $\alpha_{N,k_c^0}(M,L,\rho,h;\sG)=1$, up to multiplication by a
  $4$th root of $1$.

  (2) For every relative taut structure $\sG$, the reduced invariants
  $\Hh_N^{red}(M,L,\rho,h;\sG)$ are well defined and do not depend on
  the choice of $h$ and $\sG$.
\end{prop} 

Again this holds true because the normalized defects are functions of
boundary data on the spherical links of the vertices of the
triangulation.  

\begin{remark}\label{sQH} {\rm Point (2) of Theorem \ref{Intro-inv}
    suggests another possible notion of non ambiguous structure. While
    the ``ordinary" one that we use is defined via non ambiguous transits
    of triangulations just endowed with a pre-branching, we can
    consider QH transits of QH triangulations which enhance non
    ambiguous pre-branching transits. Let us denote by $\sG^{QH}$ such
    a kind of ``QH" non ambiguous structure. Note that it dominates an
    ordinary non ambiguous structure $\sG$, and incorporates some triple
    $(\rho,h,k)$, but different $\sG^{QH}$'s can incorporate the same
    $\sG$. Then, via point (2) of Theorem \ref{Intro-inv}, it is almost
    immediate that we can defined invariants $\alpha_N(M,\sG^{QH})$ and
    $\Hh_N^{red}(M,\sG^{QH})$. However, this definition is not so
    interesting for the following reasons: 
\begin{enumerate} 
\item In the case of cusped manifolds, the invariants
      $\alpha_N(M,\sG^{QH})$ factorize through the invariants
      $\alpha_N(M,\rho,k_c;\sG)$ which are much stronger.  

\item Also in the
      case of pairs $(M,L)$ the invariants  $\alpha_N(M,L,\sG^{QH})$
are well defined. However, it happens that infinitely many QH-non
      ambiguous structures $\sG^{QH}$, distinguished by the respective
      invariants $\alpha_N(M,L,\sG^{QH})$, dominate the same basic $\sG$ and 
incorporate the same tuple $(M,L,\rho,h,k)$ (see Section \ref{MLdefect}).

\item The ordinary non ambiguous structures $\sG$ should support the 
``universal phenomenon'' depicted at the beginning of this Introduction.
\end{enumerate}}
\end{remark}

In Section \ref{EXAMPLES} we analyse several examples in details.  We
consider at first the trivial bundle over $S^1$ with fiber a torus
with one puncture; there is one taut structure associated to the
infinite family of multiples of the natural fibration, and we show
that the symmetry defects are constant on QH triangulations which are
layered for these fibrations (as it must be). Then we describe several
examples of non-ambiguous structures (in particular some taut ones) on
some simple cusped manifolds: the figure-eight knot complement, its
sister, and the Whitehead link complement. In each case we show that
the non-ambiguous structures are distinguished by the symmetry
defects, and for the Whitehead link complement, the symmetry defects
distinguish taut structures associated to fibrations lying over non
opposite faces of the Thurston ball. So, they would separate these
faces if, for instance, the map $s_\Ff$ discussed above were constant
over them.  

\subsection{On reduced Turaev-Viro invariants}\label{RTV}
In the Appendix we quickly verify our thesis  
for the most fundamental prototype of $3$-dimensional state sums,
the Turaev-Viro ones \cite{TV}. As an application we indicate 
a procedure to construct reduced TV invariants of fibred knots in $S^3$.

\medskip

{\bf Acknowledgments.} We had very useful discussions with I. Agol,
N. Dunfield, S. Schleimer, and H. Segerman on the matter discussed in
Section \ref{FIBER}. We also thank the referees, whose suggestions allowed us to improve the exposition of our results.

\section{Generalities on triangulations}\label{TRIANG}
We will work on a given compact connected oriented smooth $3$-manifold
$M$. We denote by $\hat M$ the space obtained by collapsing to one
point each boundary component. Equivalently, $\hat M$ is obtained by
compactifying the interior ${\rm Int}(M)=M\setminus \partial M$ of $M$
by adding one point ``at infinity'' at each end. Hence $M \neq \hat M$
if and only if $\partial M \neq \emptyset$ and in such a case the {\it
  non manifold} points of $\hat M$ are the points of $\hat M \setminus
{\rm Int}(M)$ corresponding to the non spherical components of
$\partial M$. We use triangulations $T$ of $\hat M$ which are not
necessarily regular, that is, self and multiple adjacencies of
tetrahedra are allowed, and such that the set $V$ of vertices contains
$\hat M \setminus {\rm Int}(M)$.  A triangulation is {\it ideal} if $M
\neq \hat M$ and
$$ V= \hat M \setminus {\rm Int}(M). $$
The cell decomposition obtained by removing the vertices from an ideal
triangulation of $\hat M$ is also called an ``ideal triangulation'' of
${\rm Int}(M)$. Every triangulation of $\hat M$ is realized by smooth
cells in Int$(M)$ and is considered up to isotopy. It is often
convenient to consider a triangulation of $\hat M$ as a collection of
oriented tetrahedra $\Delta_1,\dots, \Delta_s$ equipped with a
complete system $\sim$ of pairings of their $2$-faces via orientation
reversing affine isomorphisms, and a piecewise smooth homeomorphism
between the oriented quotient space $$T := \coprod_{i=1}^s
\Delta_i/\sim $$ and $\hat M$, preserving the orientations. Then we
will distinguish between the ``abstract'' $j$-faces, $j=0,1,2,3$, of
the disjoint union $\textstyle \coprod_{i=1}^s \Delta_i$, and the $j$-faces
of $T$ after the $2$-face pairings. In particular we denote by
$E(\{\Delta_i\})$ and $E(T)$ the set of edges of $\textstyle
\coprod_{i=1}^s \Delta_i$ and $T$ respectively, and we write $E\to e$
to mean that an edge $E\in E(\{\Delta_i\})$ is identified to $e \in
E(T)$ under the $2$-face pairings. The $2$-faces of each tetrahedron
$\Delta_i$ have the boundary orientation defined by the rule: {\it
  ``first the outgoing normal''}. We will also consider triangulations
of a closed surface $S$ with the analogous properties.
\medskip

\begin{figure}[ht]
\begin{center}
 \includegraphics[width=3.5cm]{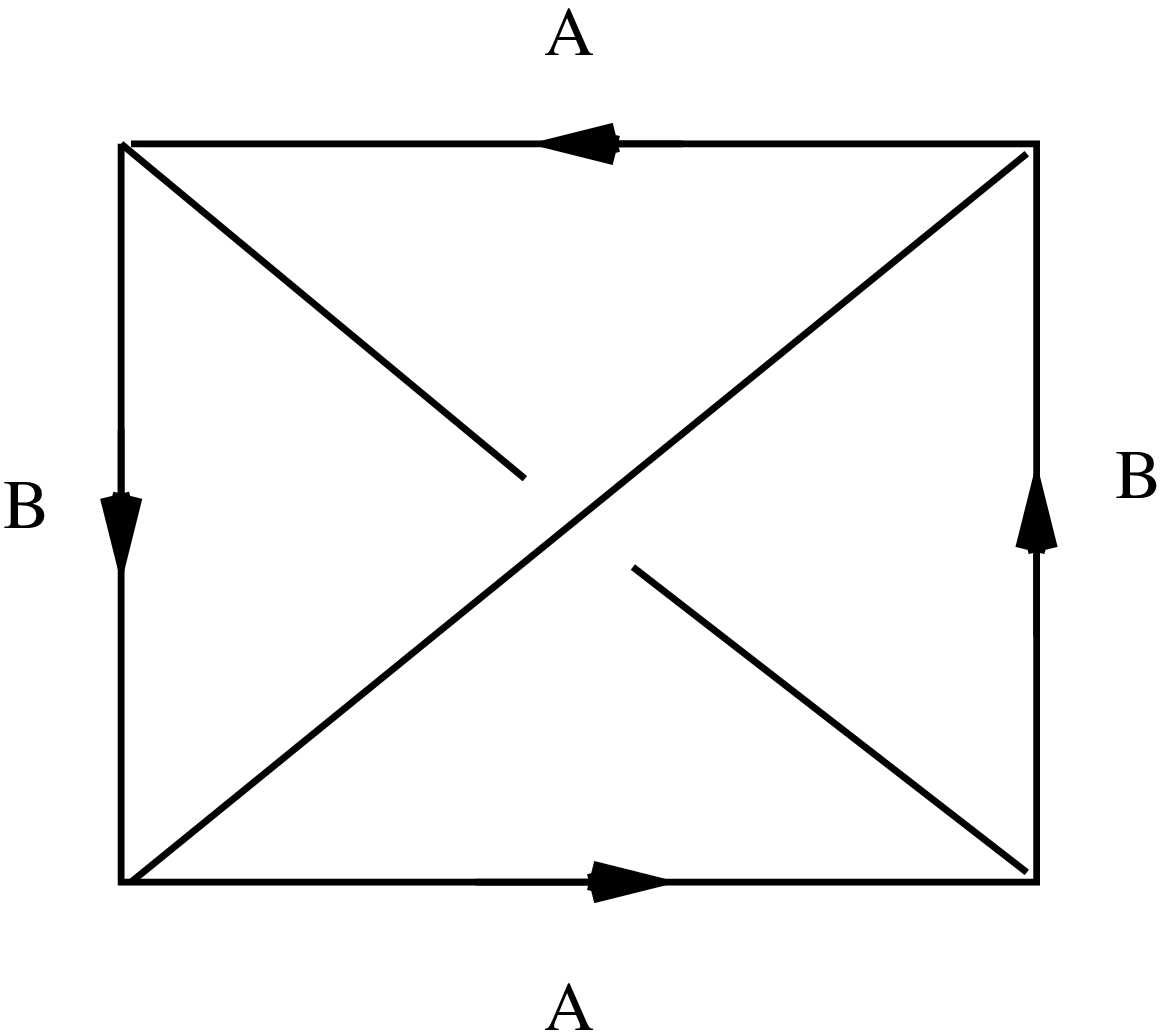}
\caption{\label{pre_b} }
\end{center}
\end{figure}

{\bf On branchings.}  We consider here with more details the notions
already mentioned in the Introduction.  A pre-branched triangulation
$(T,\omega)$ of $\hat M$ is a triangulation $T$ equipped with a {\it
  pre-branching} $\omega$; this assigns a transverse orientation to
each $2$-face of $T$ (also called a {\it co-orientation}), in such a
way that for every abstract tetrahedron $\Delta$ of $T$ two
co-orientations are ingoing and two are outgoing.  As $M$ is oriented,
a pre-branching can be equivalently expressed as a system of ``dual''
orientations of the $2$-faces of $T$. A (local) pre-branching on
$\Delta$ is illustrated in Figure \ref{pre_b}; it shows the
tetrahedron embedded in $\mr^3$ and endowed with the orientation
induced from the standard orientation of $\mr^3$. The pre-branching is
determined by stipulating that the two $2$-faces above (resp. below)
the plane of the picture are those with outgoing (resp. ingoing)
co-orientations. This specifies two {\it diagonal edges} and four {\it
  square edges}. Every square edge is oriented as the common boundary
edge of two $2$-faces with opposite co-orientations. So the square
edges form an oriented quadrilateral. Using the orientation of
$\Delta$, one can also distinguish among the square edges two pairs of
opposite edges, called {\it $A$-edges} and {\it $B$-edges}
respectively. The orientation of the diagonal edges is not
determined. Note that the total inversion of the co-orientations
preserves the pair of diagonal edges as well as the colors $A$, $B$ of
the square edges.

\begin{figure}[ht]
\begin{center}
 \includegraphics[width=7cm]{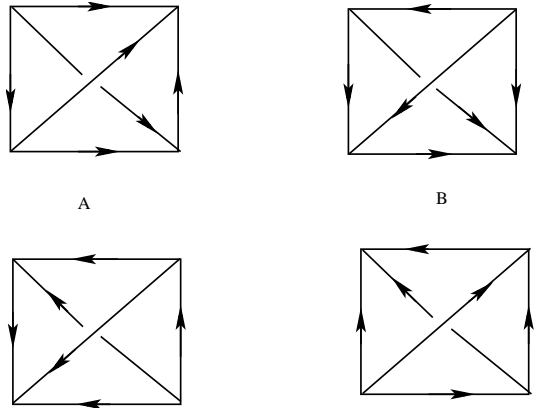}
\caption{\label{Branched_Delta} Branched tetrahedra inducing the same
pre-branched tetrahedron.}
\end{center}
\end{figure}

A tetrahedron $\Delta$ becomes a $3$-{\it simplex} by ordering its
vertices. This is equivalent to a system $b$ of orientations of the
edges, called a (local) {\it branching}, such that the vertex $v_j$
has $j$ incoming edges ($j\in \{0,\ldots,3\}$). The $2$-faces of
$(\Delta, b)$ are ordered as the opposite vertices, and $b$ induces a
branching $b_F$ on each $2$-face $F$. The branchings $b$ and $b_F$
define orientations on $\Delta$ and $F$ respectively, the {\it $b$-}
and {\it $b_F$-orientations}, defined by the vertex orderings up to even permutations. If $\Delta$ is already oriented, then
the $b$-orientation may coincide or not with the given orientation. We
encode this by a sign, $*_b\in \{-1,+1\}$. The boundary orientation
and the $b_F$-orientation agree on exactly two $2$-faces. Hence $b$
induces a (local) pre-branching $\omega_b$. On another hand, given a
pre-branching $\omega$ on $\Delta$ there are exactly four branchings
$b$ such that $\omega_b =\omega$. They can be obtained by choosing an
$A$- (resp. $B$-) edge, reversing its orientation, and completing the
resulting orientations on the square edges to a branching $b$ (this
can be done in a single way; see Figure \ref{Branched_Delta}). Note
that $*_b=1$ (resp. $*_b=-1$) if and only if one chooses an $A$
(resp. $B$) square edge, and this square edge is eventually $[v_0,v_3]$. The
diagonal edges are $[v_0,v_2]$ and $[v_1,v_3]$.

A {\it weakly-branched triangulation} $(T,\tilde b)$, with abstract
tetrahedra $\{\Delta_j\}$, consists of a system of branched
tetrahedra $\{(\Delta_j, b_j)\}$ such that the induced pre-branched
tetrahedra $\{(\Delta_j,\omega_{b_j})\}$ match under the 2-face pairings to 
form a (global) pre-branched
triangulation $(T,\omega)$. We write $\omega=\omega_{\tilde b}$.

A {\it branched triangulation} $(T,b)$, with abstract
tetrahedra $\{\Delta_j\}$, consists of a system of branched
tetrahedra $\{(\Delta_j, b_j)\}$ such that the branchings (ie. the edge
$b$-orientations) match under the 2-face pairings. 

Let us recall a few facts (see \cite{AGT} or \cite{BP} for more details):

\begin{itemize}
\item Every triangulation $T$ of $\hat M$ carries pre-branchings $\omega$;

\item For every pre-branching $\omega$ there is weak branching $\tilde
  b$ on $T$ such that $\omega=\omega_{\tilde b}$.

\item A branching is a weak-branching of a special kind. Endowing $T$
  with a branching is equivalent to promote $T$ to a $\Delta$-{\it
    complex} in the sense of \cite{HATCHER}. In general there are
  naked triangulations which do {\it not} carry any branching. 
But for every $M$ 
  there are branched (possibly ideal) triangulations $(T,b)$
  of $\hat M$.
\end{itemize}

\section{Non ambiguous ideal structures }\label{NAIS}
In this Section we restrict to {\it ideal} triangulations of a given $\hat
M$ (hence $\partial M \neq \emptyset$). 
These {\it naked} ideal triangulations are considered up to the
equivalence relation generated by isotopy relative to the set of vertices
$V$, the $2 \leftrightarrow 3$ (also called {\it MP}) {\it move}, and
the $0 \leftrightarrow 2$ (also called {\it lune}) {\it move}. These moves are {\it embedded}, and keep $V$ fixed pointwise. We call this equivalence relation the {\it (naked) ideal transit equivalence}. It is a fundamental, well known fact (due to Matveev,
Pachner, and Piergallini) that the quotient set of naked ideal 
triangulations up to ideal transit
  equivalence consists of {\it one}
point. In presence of additional structures on $T$, we consider enhanced
versions of the transit equivalence. In what follows we will often
confuse two possible meanings of a triangulation move: as a local
modification on a portion of a given triangulation, or as  an 
``abstract'' modification
pattern that can be implemented to modify a global triangulation. Then, when we will say
that a (possibly enhanced) move preserves a certain property,
we will mean that this holds true whenever we implement the move
on any triangulation verifying that property. 

\medskip

\begin{figure}[ht]
\begin{center}
\includegraphics[width=7cm]{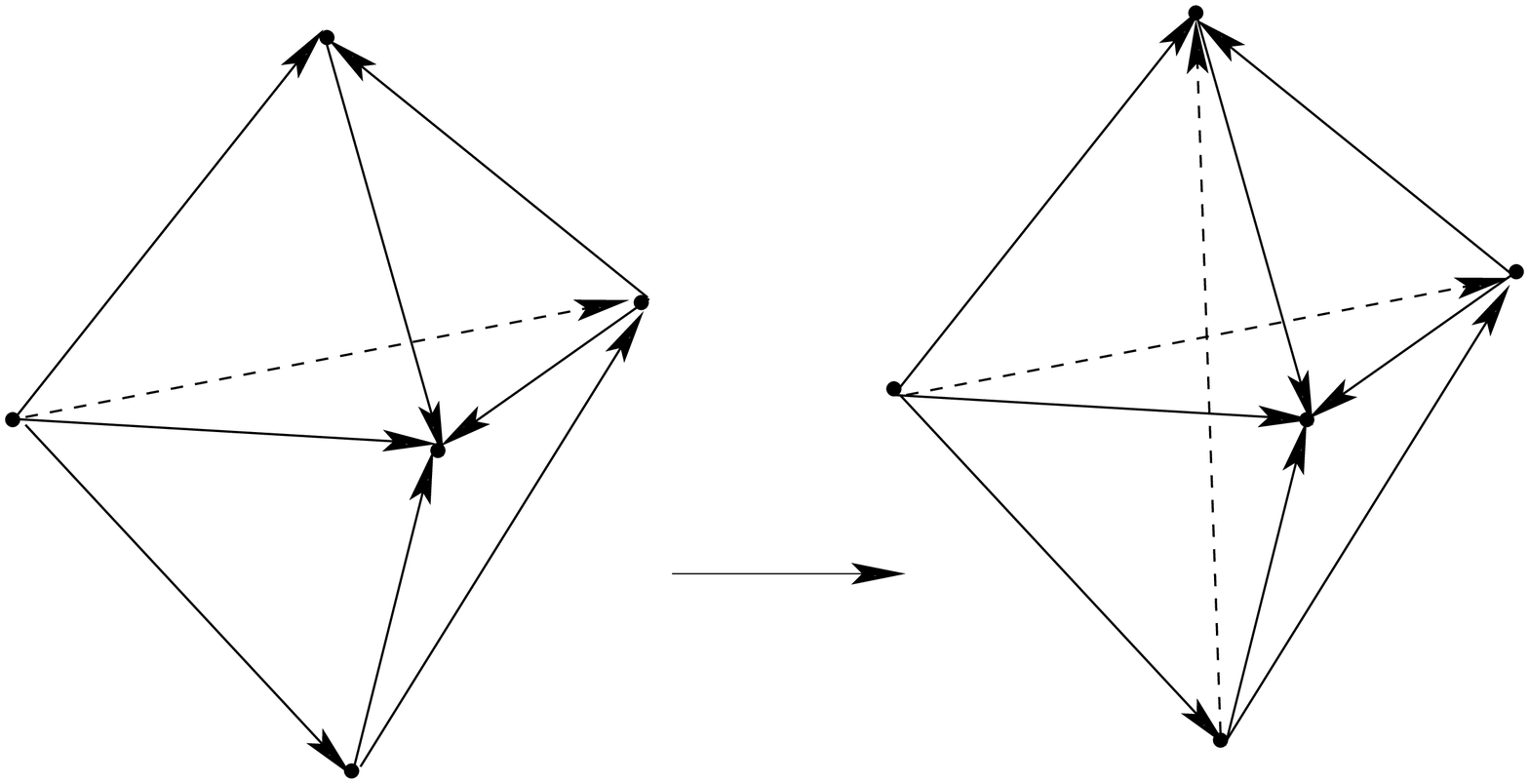}
\caption{\label{NAb} }
\end{center}
\end{figure}

{\bf On $2\leftrightarrow 3$ (MP) transits.}
Let $T \rightarrow T'$ be a $2\leftrightarrow 3$ triangulation move
between naked ideal triangulations of $\hat M$. The ``positive''
$2\rightarrow 3$ move is shown in Figure \ref{NAb} 
(the branching shown in the picture will be used later).  
Given pre-branchings $\omega$ on $T$ and $\omega'$ on
$T'$, we say that $(T,\omega)\rightarrow (T',\omega')$ is a {\it
  pre-branching transit} if for every $2$-face $F$ which is common to
$T$ and $T'$ the $\omega$ and $\omega'$ co-orientations of $F$
coincide.

Assume that we are given a pre-branched triangulation $(T,\omega)$ of
$\hat M$. Consider a naked $2\rightarrow 3$ move $T\rightarrow T'$. Denote
by $e_0$ the edge of $T'$ produced by the move, and by
$(\Delta_1,\omega_1)$, $(\Delta_2,\omega_2)$ the two (abstract) pre-branched
tetrahedra of $(T,\omega)$ involved in the move, having a common $2$-face
in $T$; recall that their
edges are either diagonal edges, or square edges colored by $A$ or
$B$. Then, a quick inspection shows that $T\rightarrow T'$ supports
always some pre-branching transit $(T,\omega)\rightarrow
(T',\omega')$, and that one of the following exclusive possibilities
is eventually realized:
\begin{itemize}
\item {\bf (NA-transit)} The pre-branched tetrahedra
  $(\Delta_j,\omega_j)$, $j=1,2$, have exactly one square edge $e$ in
  common on the shared $2$-face. Necessarily, $e$ is {\it
    monochromatic}, in the sense that the two (abstract) square edges
  identified along $e$ have the same color. In such a situation,
  $T\rightarrow T'$ supports a {\it unique} pre-branching transit
  $(T,\omega)\rightarrow (T',\omega')$; we call it a {\it non
   ambiguous MP transit}. Among the three abstract edges identified
  along $e_0$, two are diagonal edges, and the color of the square
  edge depends on the color of the monochromatic edge $e$.

\item {\bf (A-transit)} The pre-branched tetrahedra
  $(\Delta_j,\omega_j)$ have two square edges $e,e'$ in common on the
  shared $2$-face. Necessarily, both are {\it not} monochromatic. In
  such a situation $T\rightarrow T'$ supports exactly two
  pre-branching transits $(T,\omega)\rightarrow (T',\omega'_1),
  (T',\omega'_2)$. In both cases, all the abstract edges of
  $(T',\omega'_j)$ identified along $e_0$ are square edges, and $e_0$
  is not monochromatic. The two transits are distinguished by the
  prevailing color at $e_0$. We call them {\it ambiguous MP transits}.
\end{itemize}

Concerning the negative $2 \leftarrow 3$ transits we have:
\begin{itemize}  

\item A negative $2 \leftarrow 3$ pre-branching transit  
$(T,\omega) \leftarrow (T',\omega')$ is by definition {\it non ambiguous}
if it is the inverse of a positive non ambiguous transit.

\item Given a pre-branching $\omega'$ on $T'$, a ``negative'' naked
$2\leftarrow 3$ move $T \leftarrow T'$ does not support any
pre-branching transit if and only if all abstract edges around $e_0$
are square edges and $e_0$ is monochromatic. In this case we say that
$(T',\omega')$ gives rise to a stop.
\end{itemize}

\medskip

\begin{figure}[ht]
\begin{center}
\includegraphics[width=6cm]{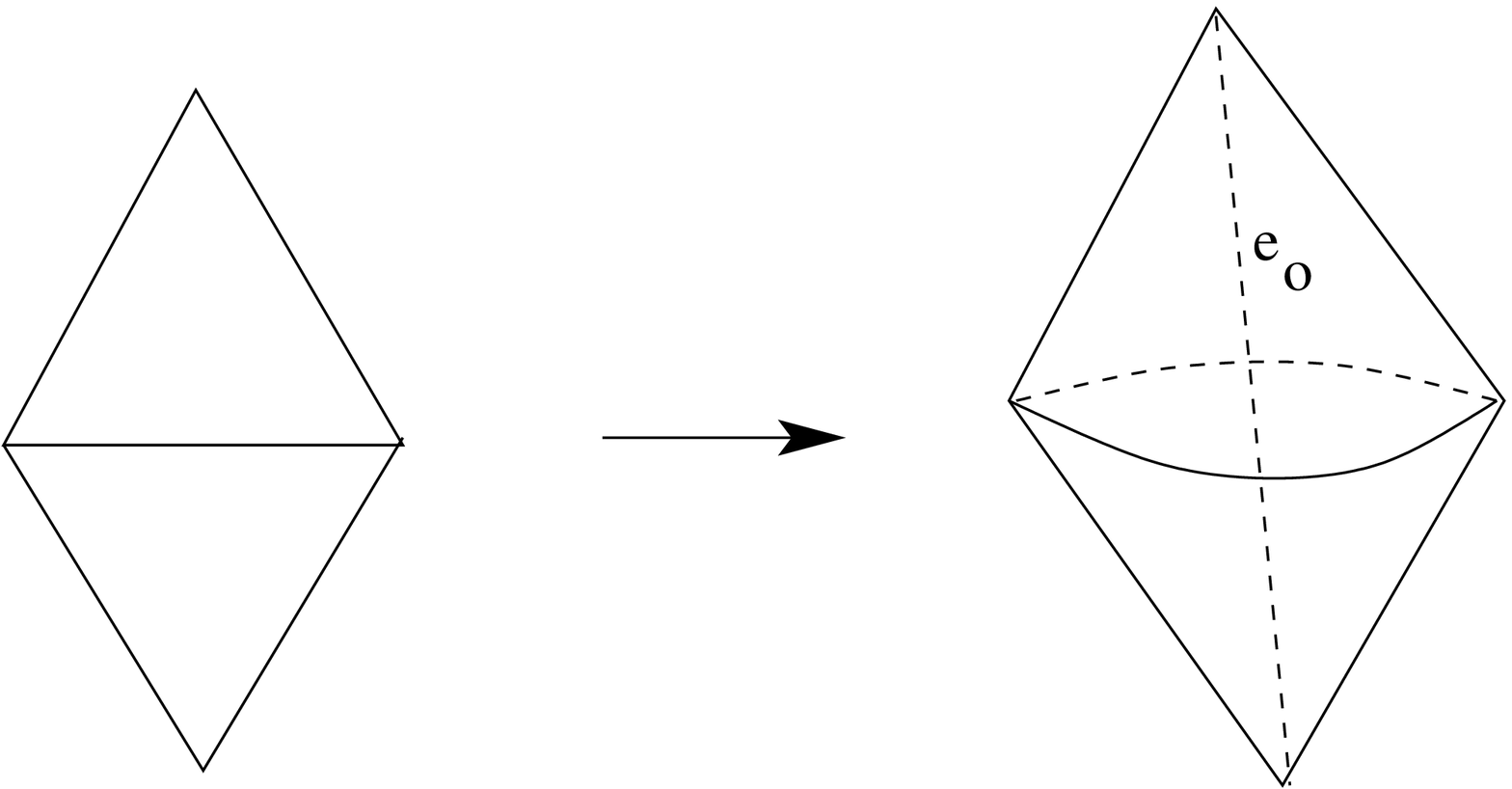}
\caption{\label{lune} }
\end{center}
\end{figure}

{\bf On $0\leftrightarrow 2$ (lune) transits.} The positive naked lune
move is shown in Figure \ref{lune}.

Let $(T,w)$ and $(T',\omega')$ be two pre-branched triangulations of
$\hat M$ such that the naked triangulations $T$ and $T'$ are related
by a positive lune move $T \rightarrow T'$. The move applies at the
union of two (abstract) $2$-faces $F_1$, $F_2$ of $T$ with a common
edge, and produces a $3$-ball $B$ triangulated by two tetrahedra glued
along two $2$-faces in $T'$ with a common edge $e_0$. The boundary of
$B$ is triangulated by two copies of $F_1 \cup F_2$ glued along their
quadrilateral common boundary. We say that $(T,\omega)\rightarrow
(T',\omega')$ is a {\it pre-branching transit} if for every $2$-face
$F$ which is common to $T$ and $T'$, the $\omega$ and $\omega'$
co-orientations of $F$ coincide, and if the restriction of $\omega'$
on the boundary of $B$ consists of two copies of the restriction of
$\omega$ to $F_1\cup F_2$. For a negative lune move, the latter
condition is replaced by: the restriction of $\omega'$ on the boundary
of $B$ consists of two copies of a same pair of co-orientations on
$F_1\cup F_2$.

It is easy to check that for every pre-branched triangulation
$(T,\omega)$, every positive lune move $T\rightarrow T'$ supports a
pre-branching transit, and that one of the following exclusive
possibilities is eventually realized:
\begin{itemize}
\item {\bf (NA-lune transit)} The $\omega$-co-orientations of $F_1$
  and $F_2$ are compatible, that is, they define a global co-orientation of
  $F_1\cup F_2$. Necessarily, the two abstract edges of $(T',\omega')$
  identified along $e_0$ are diagonal edges. In such a situation,
  $T\rightarrow T'$ supports a {\it unique} pre-branching transit
  $(T,\omega)\rightarrow (T',\omega')$; we call it a {\it non
    ambiguous} lune transit.
\item {\bf (A-lune transit)} The $\omega$-co-orientations of $F_1$ and
  $F_2$ are ``opposite''. Then $T\rightarrow T'$ supports exactly two
  pre-branching transits $(T,\omega)\rightarrow (T',\omega'_1),
  (T',\omega'_2)$, and in both cases the two edges identified along
  $e_0$ are square edges and $e_0$ is not monochromatic. We call them 
  {\it ambiguous} lune transits.
\end{itemize}

Concerning the negative lune transits we have:
\begin{itemize}  
\item   A negative  pre-branching lune transit  
$(T,\omega) \leftarrow (T',\omega')$ is by definition {\it non ambiguous}
if it is the inverse of a positive non ambiguous transit.

\item Given a pre-branching $\omega'$ on $T'$, a negative naked lune
move $T \leftarrow T'$ does not support any pre-branching transit if
and only if the two edges identified along $e_0$ are square edges,
$e_0$ is monochromatic, and the two tetrahedra have no common
diagonal edge. Again, in this case we say that $(T',\omega')$ gives
rise to a stop.
\end{itemize}

\begin{defi} \label{NAS} {\rm The {\it non ambiguous ideal $pb$-transit
    equivalence} on the set of ideal pre-branched triangulations of $\hat M$ is
    generated by isotopy relative to the set of vertices $V$, the non
    ambiguous $2 \leftrightarrow 3$ (MP) and the non ambiguous $0
    \leftrightarrow 2$ (lune) pre-branching transits. We denote by
    $\mathcal{NA}^{id}(M)$ the quotient set. We call a coset in
    $\mathcal{NA}^{id}(M)$ a {\it non ambiguous structure} on $M$.}
\end{defi}

Clearly, by allowing arbitrary ideal
pre-branching transits, we get the general {\it ideal $pb$-transit equivalence},
with quotient set $\Pp\Bb^{id}(M)$, which is in fact a quotient  of
$\Nn\Aa^{id}(M)$. 

\begin{remark}\label{total-inv} {\rm The {\it total inversion} of a 
global pre-branching,  
that is, the simultaneous reversal of all the
$2$-face co-orientations,  induces an involution on  $\Nn\Aa^{id}(M)$.
We set the question whether this involution is different from the
identity. Anyway it is immaterial with respect to most future
developments, so that one could even incorporate the total
inversion among the generators of the  non ambiguous $pb$-transit
equivalence.}
\end{remark}

As we have defined the notion of non ambiguous structure
purely in terms of pre-branching, we see that the latter is
the fundamental triangulation enhancement underlying our discussion.
However, it is useful (and necessary when we deal with the quantum hyperbolic state sums)
to treat the pre-branchings also in terms of other enhanced transits.

The (non ambiguous) ideal $pb$-transit equivalence can be somewhat
tautologically rephrased in terms of weak branchings. 
By definition a (non ambiguous)
ideal $wb$-{\it transit} $(T,\tilde b)
\leftrightarrow (T',\tilde b')$ is such that it dominates an associated 
(non ambiguous) 
$pb$-transit $(T,\omega_{\tilde b})
\leftrightarrow (T',\omega_{\tilde b'})$, and moreover $\tilde b$ and $\tilde
b'$ coincide on the common tetrahedra of $T$ and $T'$. Consider on the set of
weakly branched ideal triangulations of $\hat M$ the equivalence
relation generated by the $wb$-transits, imposing furthermore that
$(T,\tilde b)\sim (T,\tilde b')$ if $\omega_{\tilde b}=\omega_{\tilde
  b'}$. It is then obvious that the correspondence $\tilde b \to
\omega_{\tilde b}$ induces a bijection between the quotient set, say
$\Ww\Bb^{id}(M)$, and $\Pp\Bb^{id}(M)$. By restricting to non ambiguous transits,
we get a way to treat $\Nn\Aa^{id}$ in terms of weak branchings.

\smallskip

A naked ideal move $T\rightarrow T'$ supports a $b$-{\it transit}
$(T,b)\rightarrow (T',b')$, for some branchings $b$ and $b'$, if at
every common edge $e$ of $T$ and $T'$, the $b$- and $b'$-orientations
coincide. It is immediate that  every $b$-transit
dominates a $pb$-transit $(T,\omega_b)\rightarrow (T',\omega_{b'})$ of
the underlying pre-branchings.

\begin{defi}\label{NAB}{\rm An ideal $b$-transit  $(T,b)\rightarrow (T',b')$
is {\it non ambiguous} if the associated $pb$-transit is non ambiguous.}
\end{defi}

This definition is coherent with the above discussion if we consider
a branching as a special kind of weak branching. On the other hand, 
we stress that it is {\it not} the immediate definition of non
ambiguous $b$-transit one would wonder. It is actually stronger.  Let
us say that a positive ideal $b$-transit $(T,b)\rightarrow (T',b')$ is
{\it forced} if it is the unique $b$-transit that enhances the naked move 
$T\to T'$ starting with $(T,b)$. We have:

\begin{lem}\label{forcedB}
(1) If $(T,b)\to (T,b')$ is non ambiguous, then it is forced. 
  
(2) If $(T,b)\to (T,b')$ is not forced, then there are exactly two
$b$-enhancements $(T,b)\to (T',b_1)$ and $(T,b)\to (T',b_2)$ which
dominate the respective ambiguous $pb$-transits $(T,\omega_b)\to
(T',\omega_{b_1})$ and $(T,\omega_b)\to (T',\omega_{b_2})$.

(3) Given a negative naked ideal move $T\leftarrow  T'$, a branching
$(T',b')$ gives rise to a stop (ie. there is no $b$-enhancement
starting with $(T',b')$) if and only if $(T',\omega_{b'})$ gives rise to a 
stop of $pb$-transits.

(4) There are instances of forced positive ideal $b$-transits 
$(T,b)\to (T',b')$ which nevertheless are ambiguous.
\end{lem}
    
The proof of (1)--(3) is easy by direct inspection. As for (4), a
$2\rightarrow 3$ example is given by the following
configuration: let $(\Delta_1,b_1)$ and $(\Delta_2,b_2)$ be the two
(abstract) branched tetrahedra of $T$ involved in the move. Denote by
$v_1$ and $v_2$ the vertices of $\Delta_1$ and $\Delta_2$ respectively
which are opposite to the common $2$-face. Then the branching is such
that $v_1$ is a source while $v_2$ is a pit. 
There are similar lune move examples. 
\medskip

{\bf Two remarkable non ambiguous $b$-transits.}
In Figure \ref{NAb} we show a non ambiguous $b$-transit
which dominates a non ambiguous $2\rightarrow 3$ pre-branching transit
such that the common square edge of the two tetrahedra of $(T,\omega)$
is $A$-monochromatic. This $b$-transit has the nice property that all
branched tetrahedra are positively $b$-oriented ($*_b=1$); when the common
square edge of the two tetrahedra of $(T,\omega)$ is
$B$-monochromatic, there is a similar non ambiguous $b$-transit where
the five tetrahedra are negatively $b$-oriented ($*_b=-1$).

\medskip

A weak branching on $T$ can be a genuine branching on some portion of
$T$. In particular an ideal $wb$-transit $(T,\tilde b)\to (T',\tilde
b')$ is said {\it locally branched} if the restriction of $\tilde b$
or $\tilde b'$ to the portions of $T$ and $T'$ involved in the move is
a genuine branching. In such a case the discussion about the non
ambiguous $wb$-transits is (locally) confused with the one in terms of
$b$-transits. With this notion, we can elaborate a little bit on the
above bijection $ \Ww\Bb^{id}(M)\to \Pp\Bb^{id}(M)$ induced by the
correspondence $\tilde b \to \omega_{\tilde b}$, which is
useful in the following form (for example it has been used in \cite{BP}, Proposition
3.4, and in \cite{AGT}):

\begin{lem}\label{wb} The equivalence relation that realizes
  $\Ww\Bb^{id}(M)$ (hence $\Pp\Bb^{id}(M)$) is generated by the local
  {\rm tetrahedral moves} which change the local branching on any
  tetrahedron of $(T,\tilde b)$ by preserving the induced
  pre-branching, and the {\rm locally branched} ideal
  $wb$-transits. Moreover, via tetrahedral moves we can realize every
  locally branched configuration compatible with a given $pb$-transit.
  By restricting to locally branched non ambiguous transits, we have a
  similar realization of $\Nn\Aa^{id}(M)$.
\end{lem}

\section{Charges and taut structures} \label{seccharges} In this
section we point out the cohomological content and remarkable
specializations of the non ambiguous structures. An issue is to stress
the compatibility of the notions of transit which arise in
different contexts.

Let us recall the notion of charge (called $\mz$-charge hereafter)
which plays an important role in QH geometry.  Given a naked (not
necessarily ideal) triangulation $T$ of $\hat M$, a $\mz$-{\it charge}
$c$ assigns to every edge $E$ of every abstract tetrahedron $\Delta$
of $T$ a color $c(E)\in \Z$, in such a way that opposite edges have
the same color and the following local and global conditions are
satisfied:
\begin{itemize}
\item[(i)] For every $\Delta$ the sum of the three colors is equal to
  $1$.
\item[(ii)] Every edge $e$ of $T$ has {\it total charge} $C(e)=2$,
  where $C(e)$ is the sum of the colors $c(E)$ of the abstract edges
  which are identified along $e$ (ie. $E\to e$).
\end{itemize}

If $c$ is a charge, then $\pi c$ is often called a {\it $\mz$-angle
  system} on $T$.  

\smallskip

One defines similarly the notion of $\Z/2$-{\it charge} $\gamma$ by
taking colors in $\mz/2$ and considering the above conditions with
coefficients in $\mz/2$. 

\begin{defi}\label{taut-charge} {\rm A $\mz/2$-charge
    (resp. $\mz$-charge) is {\it locally taut} if for every $\Delta$,
    one color is $1$ and the others are $0$. It is {\it $\Z/2$-taut}
    (resp. $\mz$-{\it taut}) if moreover for every edge $e$ of $T$, there
    are at least two (resp. exactly two) $1$-colors ``around" $e$, that is, 
    $1$-colored abstract edges $E$ such that $E\rightarrow e$ (NB: a
    {\it locally taut} $\mz$-charge is automatically $\mz$-taut).
    \smallskip

If $c$ is a taut $\Z$-charge, then $\pi c$ is also called a {\it taut
  $\mz$-angle system}. }
 \end{defi}
 
 As above, we will often write that some configuration of abstract edges takes place 
 ``around'' an edge $e$ of $T$ when it is realized by the abstract edges $E$ such 
that $E\rightarrow e$.

\begin{defi}\label{taut-triang}{\rm A pre-branched triangulation
    $(T,\omega)$ is {\it taut} (\cite {L}) if for every edge $e$ of
    $T$ there are exactly two diagonal edges around $e$. 
Let us call it {\it $\Z/2$-taut} 
    if around every edge there are abstract diagonal edges 
(then they are at least two).}
\end{defi} 

\smallskip

Let $(T,\omega)$ be a pre-branched triangulation of $\hat M$. For
every abstract edge $E$, set $\gamma_\omega(E)=1$ if $E$ is a diagonal
edge, and $\gamma_\omega(E)=0$ if it is a square edge; here $0,1$
belong to either $\Z/2$ or $\Z$ in accordance with the context.

The following results (together with those about taut triangulations
in Section \ref{TAUT}) express conditions under which the various
kinds of charges defined above exist, and their relations with
pre-branched and taut triangulations. 

\begin{prop}\label{N2_via_pb}
  (1) Let $(T,\omega)$ be a pre-branched triangulation of $\hat
  M$. Then $\gamma_\omega$ is a locally taut $\mz/2$-charge.
  Conversely, every locally taut $\Z/2$-charge $\gamma$ on a
  triangulation $T$ of $\hat M$ lifts to a locally taut $\Z/2$-charge
  $\gamma^*$ on a pre-branched triangulation $(T^*,\omega^*)$ of a
  double covering $M^* \to M$, such that $\gamma^* =
  \gamma_{\omega^*}$.
  
  \smallskip
  
  (2) If $(T,\omega)$ is a $\Z/2$-taut (resp. taut) triangulation,
  then $\gamma_\omega$ is a taut $\Z/2$-charge ($\Z$-charge).

\smallskip
  
(3) If a triangulation $T$ of $\hat M$ admits a taut $\Z/2$-charge,
then $M\neq \hat M$, $T$ is an ideal triangulation, and $\partial M$
has no spherical component. If moreover all the components of
$\partial M$ are tori, then the charge is the reduction mod$(2)$ of
a taut $\Z$-charge.

\smallskip
    
  (4) A triangulation $T$ of $\hat M$ admits a $\Z$-charge if and only
  if $\hat M \neq M$, every boundary component of $M$ is a torus, and
  $T$ is an ideal triangulation.

\end{prop}

\Dim The first claim in (1) follows easily from an analysis of the
boundary configurations in the star of each vertex; we postpone this
to Section \ref{PRETOB}.  As for the second claim, every locally taut
$\Z/2$-charge $\gamma$ on $T$ determines on every abstract tetrahedron
$\Delta$ of $T$ a local pre-branching which is unique up to total
inversion. Then $\gamma=\gamma_\omega$ for some $\omega$ if and only
if we can choose a compatible family of such local pre-branchings. The
obstruction is a $1$-cohomology class mod$(2)$ which vanishes up to
passing to a double covering.  
Point (2) follows immediately from the definitions.  
The direct implications in (3) and (4) use a Gauss-Bonnet argument
(see also Section \ref{train-trak} below). Namely, given a taut
$\Z/2$-charge on $T$, the boundary of a small open star of a vertex
$v$ of $T$ is a triangulated surface $S_v$ whose number $t$ of
triangles is equal to the number of germs of $1$-colored abstract
edges ending at $v$, which is bigger than twice the number $s$ of
vertices of $S_v$. Hence $t-2s=-2\chi(S_v)\geq 0$, which shows that
$v$ is a singular point of $\hat M$ and $T$ is an ideal
triangulation. For a $\mz$-charge we have $t=2s$. The converse
implication of (4) is much harder, and proved by applying the
arguments of W. Neumann's ``flattening'' theory (\cite{N}, see also
\cite{AGT}, Section 4.3).\cvd
\medskip
 
There is a natural notion of transit $(T,c) \rightarrow (T',c')$
between two $\Z$-charged triangulations of $\hat{M}$, which is widely
used in the theory of QHI: one requires that $c'$ coincides with $c$
on every common abstract tetrahedron of $T$ and $T'$, and moreover,
considering the restrictions of $c$ and $c'$ to the polyhedron supporting the move,
one requires that the local charge condition is verified on every
abstract tetrahedron, the global one is verified around every internal
edge, and at every boundary edge the value of the total charge is
preserved. One defines similarly the notion of transit between
$\Z/2$-charged triangulations. By dealing with $\Z/2$-charges we can
use arbitrary triangulations, hence we should include also the bubble
transits (see Section \ref{ML}); for $\Z$-charges, by Proposition
\ref{N2_via_pb}, it is necessary to restrict to ideal triangulations
and to consider only ideal transits. For simplicity, let us restrict
anyway to the ideal setting. 
For both $\mz/2$- and $\mz$-charges, as in Definition
\ref{NAS}, one can consider the quotient sets under the equivalence
relations generated by the relevant charge transits. Let us denote them by
\begin{equation}\label{sets}
  c^{id}(M,\Z/2)\ , \ c^{id}(M,\Z)
\end{equation}
respectively. These quotient sets have
a clear intrinsic topological meaning:
\begin{prop}\label{Z2cmeaning}  
  (1) The set $c^{id}(M,\Z/2)$, encodes $H^1(M;\Z/2)$: for every
  ideal triangulation $T$ of $\hat M$ and every class $\alpha \in
  H^1(M;\Z/2)$, there is a $\Z/2$-charge on $T$ which realizes
  $\alpha$, and any two $\Z/2$-charged triangulations of $\hat M$
  realizing a same class $\alpha$ are equivalent under $\Z/2$-charge
  transits.
 \smallskip
  
  (2) The set $c^{id}(M,\Z)$ encodes
$$W_c(M):=\{(h,k)\in H^1(M;\Z/2\Z) \times H^1(\partial M; \Z) \ | 
\ k= i^*(h) \ {\rm mod}(2)\}$$ where $i: \partial M \to M$ is the
inclusion. That is, for every ideal triangulation $T$ of $\hat M$ and
every $\alpha \in W_c(M)$ there is a $\Z$-charge on $T$ which
represents $\alpha$, and any two $\Z$-charged triangulations of $\hat
M$ realizing a same $\alpha$ are equivalent under $\Z$-charge
transits.
\end{prop}

\Dim These results are extracted from the
theory of {\it cohomological weights} developed in \cite{GT,AGT},
largely elaborating on W. Neumann's theory of ``flattening'' \cite{N}.
How a relevant charge encodes a class in $H^1(M;\Z/2)$ or $W_c(M)$ can
be found in \cite{AGT}, Proposition 4.8, as well as the fact that
every cohomology class can be encoded in this way. Let us remind such
an encoding. Represent any non zero class in $H_1(\partial M;\mz)$ by
{\it normal loops}, that is, a disjoint union of oriented essential
simple closed curves in $\partial M$, transverse to the edges of
the triangulation $\partial T$ induced by $T$ on $\partial M$, and such that no curve enters and exits a triangle by a
same edge (the triangulation $\partial T$ is considered in detail in Section \ref{PRETOB}). The intersection of a normal loop, say $C$, with a
triangle $F$ of $\partial T$ consists of a disjoint union of arcs,
each of which turns around a vertex of $F$; if $F$ is a cusp section
of the tetrahedron $\Delta$ of $T$, for every vertex $v$ of $F$ we
denote by $E_v $ the edge of $\Delta$ containing $v$. We write $C
\rightarrow E_v$ to mean that some subarcs of $C$ turn around $v$. We
count them algebraically, by using the orientation of $C$: if there
are $s_+$ (resp. $s_-$) such subarcs whose orientation is compatible
with (resp. opposite to) the orientation of $\partial M$ as viewed
from $v$, then we set $ind(C,v) :=s_+-s_-$.  For every $\Z$-charge $c$
on $T$, one defines the cohomology class $\gamma(c)\in
H^1(\partial M;\mz)$ by ($C$ is a normal loop on $\partial M$):
\begin{align}
 \gamma(c)([C]) & := \sum_{C \ra E_v} ind(C,v) c(E_v).\label{weightc}
\end{align} 
Another class $\gamma_2'(c) \in H^1(M;\mz/2\mz)$
is defined similarly, by using normal loops in $T$ and taking the sum
mod$(2)$ of the charges of the edges we
face along the loops. We have $\gamma(c)=i^*(\gamma_2'(c))$ mod$(2)$. The pair $(\gamma(c),\gamma_2'(c))\in W_c(M)$ is the class associated to the $\mz$-charge $c$.

Let us consider now the $\mz$-charges on an ideal triangulation
$T$ as integral vectors with entries indexed by the abstract edges of
$T$.  Then the claims about the transits are consequences of two
results: the difference between two $\mz$-charges $c_1$ and $c_2$ on
$T$ having equal class $\alpha$ in $W_c(M)$ is an integral linear
combination of determined integral vectors $d(e)$ associated to the
edges $e$ of $T$ (this follows from the exact sequence (43) in
\cite{AGT}, where $c_1-c_2$ represents a class in ${\rm
  Ker}(\gamma',\gamma'_2) = 0\in H$, the zero class of $H$ being
represented in ${\rm Im}(\beta)$ by such linear combinations); if $c$
is a $\mz$-charge on $T$ and $T\rightarrow T'$ a positive
$2\rightarrow 3$ move, the affine space of $\mz$-charges $c'$ produced
by all possible transits of $\mz$-charges $(T,c)\rightarrow (T',c')$
starting with $c$ coincides with the space generated by the above edge
vectors $d(e)$.  A concrete description of the vectors $d(e)$ is given
in the proof of Proposition \ref{3D-weight-inv}.  Then, by using
finite sequences of $2\rightarrow 3$ moves $T\rightarrow \ldots
\rightarrow T$ starting and ending at $T$, and that blow down and then
up any given edge of $T$ at a certain intermediate step, one can
change a given $\mz$-charge on $T$ to any other one with same class
$\alpha$.  (This argument is fully detailed in \cite{Top}, Section
4.1, in the context of the distinguished triangulations of pairs
$(M,L)$ that we discuss in Section \ref{ML}.)\cvd
\medskip

The following lemma states that the sets of taut and $\Z/2$-taut
triangulations are closed under ideal non ambiguous transits. The next
one indicates the relation between pre-branched transits and charge
transits. The proofs are easy, basically a direct consequence of the
definitions. Here we consider a transit as an ``abstract" pattern that
can be implemented to locally modify an ideal triangulation.

\begin{lem}\label{taut-preserved}  For an ideal pre-branching transit
  the following facts are equivalent:
  
(a) The transit is non ambiguous.
  
(b) The transit sends $\Z/2$-taut triangulations to $\Z/2$-taut
  triangulations.
  
(c) The transit sends taut triangulations to taut triangulations.
  \end{lem}
  
\begin{lem}\label{on_transits}
  (1) Any pre-branching transit $(T,\omega)\leftrightarrow
  (T',\omega')$ induces a transit of locally taut $\Z/2$-charges 
  $(T,\gamma_{\omega})\leftrightarrow
  (T',\gamma_{\omega'})$.
    
  (2) Any (necessarily non ambiguous) ideal transit of taut
  (resp. $\Z/2$-taut) triangulations $(T,\omega)\leftrightarrow
  (T',\omega')$ induces a transit of $\Z$-taut (resp.  $\Z/2$-taut)
  charges $(T,\gamma_{\omega})\leftrightarrow (T',\gamma_{\omega'})$.
  \end{lem}

  \begin{remark}\label{vering}{\rm Recall (\cite{A}, \cite{GF,HRST})
      that a taut triangulation $(T,\omega)$ is {\it veering} if every
      edge $e$ is either $A$- or $B$-monochromatic (the diagonal
      abstract edges around $e$ being considered as
      achromatic). Though the non ambiguous transits preserve
      tautness, they do not preserve the property of being veering
      (see also Remark \ref{partialA}).}
    \end{remark}

    As an immediate consequence of Proposition \ref{Z2cmeaning} and
    Lemma \ref{taut-preserved} and \ref{on_transits}, the following
    Proposition summarizes the results of this Section.
  
  \begin{prop}\label{taut-and-homology}  
    (1) The taut and $\Z/2$-taut triangulations of $\hat M$
    respectively determine well defined (possibly empty) subsets
    $\tau(M)$ and $\tau(M,\Z/2)$ of $\Nn\Aa^{id}(M)$. They are called
    respectively the set of {\rm taut structures} and the set of {\rm
      $\Z/2$-taut structures} on $M$.  

\smallskip

(2) There are well defined maps 
$$ \hG: \Nn\Aa^{id}(M)\to H^1(M;\Z/2)$$
(which restricts to $ \tau(M,\Z/2)$) and    

$$(\hG,\partial \hG): \tau(M)\to H^1(M;\Z/2)\times H^1(\partial M;\Z),$$
defined via $[(T,\omega)] \to  [\gamma_\omega]$, where $[\gamma_\omega]$ is 
the class defined by the charge tautologically carried by the pre-branching, and similarly for $(\hG,\partial \hG)$.
\end{prop}  

Conditions under which $\tau(M)$ or $\tau(M,\Z/2)\ne \emptyset$ are
discussed in Section \ref{TAUT}.
\begin{remarks}\label{transit-equiv}{\rm In the present paper we are mainly
interested in the role of non ambiguous structures in the definition
of invariant ``reduced'' quantum state sums. However, the study of
(ideal) variously branched triangulations up to different 
``transit equivalences'' will be not exhausted and we are currently working on this topic. 
Here we limit ourselves to a few further information (without proofs). For every ideal pre-branched triangulation $(T,\omega)$ of $\hat M$,
the $1$-skeleton, say $X$, of the dual cell decomposition of 
$\hat M$ is oriented
by $\omega$  and becomes a (cellular) integral $1$-cycle $(X,\omega)$ in $M$.
It is easy to see that the correspondence $(T,\omega)\to [(X,\omega)]$
well defines a map
$$ {\bf h}: \Pp\Bb^{id}(M)\to H_1(M;\Z)$$
which naturally lifts to $\Nn\Aa^{id}(M)$ via the natural projection 
$\pi: \Nn\Aa^{id}(M)\to \Pp\Bb^{id}(M)$. We can add to the
generators of the $pb$-transit equivalence a further
(non local) so called ``circuit move'' acting on the pre-branchings
of any given $T$ so that the further quotient set $\Pp\Bb^{id}(M)/\sim$
reduces to one point (see \cite{BP}, \cite{AGT}). By studying the class
$[(X,\omega')-(X,\omega)]$ when $(T,\omega)$ and
$(T,\omega')$ differ from each other by a circuit move, we can prove
for example that {\it if  $H_1(M;\Z)$ is infinite, then $\Pp\Bb^{id}(M)$ is infinite}
(hence also $\Nn\Aa^{id}(M)$). 

It is also interesting to study the quotient sets $\Bb^{id}(M)$  and $\Bb\Nn\Aa^{id}(M)$ of branched ideal triangulations of $\hat M$ consider up to $b$-transits and non-ambiguous $b$-transits, respectively. For example, there is a natural
``forgetting'' map $\phi: \Bb^{id}(M)\to \Pp\Bb^{id}(M)$ and one would like to
understand its image. It is easy to see that if $\alpha\in {\rm Im}(\phi)$, then {\bf h}$(\alpha)=0 \in H_1(M;\Z)$. If $H_2(M;\Z/2)=0$, one realizes that this necessary
condition is also sufficient. In general this is not true.}
\end{remarks}

\section{Holographic approach to non ambiguous structures} 
\label{2D}

Assume that $\partial M$ is non empty, and let $T$ be an ideal
triangulation of $\hat M$.  We are going to show that the ``restrictions''
to $\partial M$ of the structures we have considered on
$3$-dimensional triangulations $T$ have a clear intrinsic meaning. In
the same time we will easily realize that the resulting
$2$-dimensional structures make sense by themselves, even when they
are not induced by $3$-dimensional ones. This ``free'' $2$-dimensional
theory deserves to be studied by itself. On the other hand, the
restrictions of $3$-dimensional structures present specific coherent,
or ``entangled", behaviours, so that the $2D$ information eventually
leads to meaningful intrinsic features of the $3$-dimensional
theory. We will see another instance of this approach in Section
\ref{SYMD}.

\subsection{From pre-branchings to boundary branchings}\label{PRETOB}
\begin{figure}[ht]
\begin{center}
 \includegraphics[width=7cm]{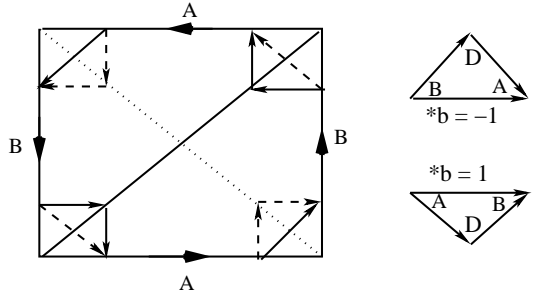}
\caption{\label{pb-boundary} The boundary triangulation.} 
\end{center}
\end{figure}

Every ideal triangulation $T$ of $\hat M$ induces a cellulation of $M$
made of truncated tetrahedra, and thus it defines a triangulation
$\partial T$ of $\partial M$, whose triangles are the triangular
faces of the truncated tetrahedra. Every pre-branching $\omega$ on $T$
induces a branching $\partial \omega$ on $\partial T$ defined locally
as in Figure \ref{pb-boundary}; we realize easily that this definition
is globally compatible. The total inversion of $\omega$ changes
$\partial \omega$ by the {\it branching total inversion},
which reverses all edge orientations.

Every (abstract) branched triangle of $(\partial T, \partial \omega)$
has, as usual, a sign $*_{\partial \omega}$ with respect to the
boundary orientation. Since $M$ is oriented the pre-branching induces
an orientation and hence a sign $*_\omega$ on every $2$-face of every
abstract tetrahedron of $(T,\omega)$. Then, we see that for every triangle
$\sigma$ of $\partial T$, $*_{\partial \omega}(\sigma)$ coincides with
the sign $*_\omega(F)$ of the $2$-face $F$ of $T$ opposite to
$\sigma$.

\begin{figure}[ht]
\begin{center}
 \includegraphics[width=3.5cm]{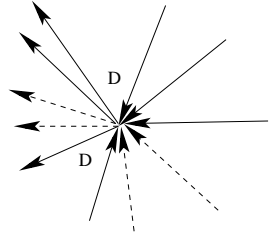}
\caption{\label{pb-boundary2} Vertex star germ.} 
\end{center}
\end{figure}

Every corner of every abstract triangle of $\partial T$ corresponds to
an abstract edge $E$ of $T$. Hence it inherits a color, either $D$ if
$E$ is a diagonal edge, or $A$ or $B$ if $E$ is a $A$- or $B$-square
edge. Forgetting the colors $A$ and $B$, Figure \ref{pb-boundary2} shows
a typical star-germ at a vertex of $(\partial T,\partial
\omega)$. Clearly there is an even number of $D$-corners around every
vertex of $\partial T$.  This proves the first claim of (1) of Lemma
\ref{N2_via_pb}. Note that the corner coloring depends only
on the sign $*_{\partial \omega}$ of each branched triangle. This is
also illustrated in Figure \ref{pb-boundary}.

\subsection{$2D$ charges, train tracks and singular
  combings}\label{train-trak}
Let $(K,b)$ be any branched triangulation of a closed oriented surface
$S$.

The notions of triangle signs and corner colorings introduced in
Section \ref{PRETOB} make sense as well for the two-dimensional 
triangulation $(K,b)$. Also, the notions
of locally taut or taut $\mz$- or $\mz/2$-charges defined in Section
\ref{seccharges} can be considered verbatim for maps $c$ assigning a
color $c(V)$ to every corner $V$ of every abstract triangle of $K$;
such a map $c$, called a {\it 2D charge}, satisfies the local charge condition (i)
on every triangle and the global charge condition (ii) about every vertex. 

Clearly, every branched triangulation $(K,b)$ carries a locally taut
$\mz/2$-charges $\gamma_b$, defined by labelling every $D$-corner with
$1$ and the $A$- and $B$-corners with $0$.  A triangulation $(K,b)$ is
$\Z/2$-{\it taut} if $\gamma_b$ is so, and $(K,b)$ is {\it taut} if
$\gamma_b$ is $\Z$-taut.  If $S=\partial  M$ and $K=\partial
T$ for some ideal triangulation $T$ of $\hat M$, it is evident that every
instance of $3D$ charge on $T$ restricts to an equally named $2D$
charge on $\partial T$, and that taut or $\Z/2$-taut triangulations $T$
restrict to triangulations $K$ qualified in the same way.

\medskip

Next we will point out how any $(K,b)$ carries further interesting
derived structures.  

\smallskip
{\bf Euler cochains.} Given a locally taut $\mz/2$-charge $\gamma$ on
$K$, we can define two integral cellular $2$-cochains ${\rm
  Eu_\gamma}$ and ${\rm Eu}^-_\gamma$ on $S$, with respect to the cell
decomposition dual to $K$. Every $2$-cell $C$ is dual to one vertex
$v$ of $K$, so set
 $${\rm Eu}_\gamma(C)= 1- \frac{r(v)}{2}\ ,
 \ {\rm Eu}^-_\gamma(C) = \min(0,1- \frac{r(v)}{2})$$ where $r(v)$ is
 the number of corners around $v$ with color $1$ given by $\gamma$.
 Note that ${\rm Eu}_\gamma(C)\in \mz$, and either ${\rm
   Eu}_\gamma(C)= 1$ or ${\rm Eu}_\gamma(C)\leq 0$; moreover, ${\rm
   Eu_\gamma} = {\rm Eu}^-_\gamma$ if and only if $\gamma$ is taut. If
 $\gamma=\gamma_b$, we will denote ${\rm
   Eu}_\gamma$, ${\rm Eu}^-_\gamma$ by ${\rm
   Eu}_b$, ${\rm Eu}^-_b$ respectively. 
 
 \begin{lem}\label{EU} (1) The Euler-Poincar\'e
 characteristic of every component $Z$ of $S$ is given by 
 $$\chi(Z)= {\rm Eu}(Z) = \sum_{C\subset Z} {\rm Eu}(C). $$
 
 (2) If $K$ supports a taut $\Z/2$-charge, then $\chi(Z)\leq 0 $ for
 every component $Z$ of $S$. If $K$ supports a taut $\Z$-charge, 
 then $\chi(Z)=0 $ for every component $Z$ of $S$; hence they are all tori.  
 
 (3) If $\gamma$ is a taut $\Z/2$-charge on $K$ and every component of
 $S$ is a torus, then $\gamma$ is the reduction mod$(2)$ of a taut
 $\Z$-charge.
 \end{lem}
 
Lemma \ref{EU} (1) can be proved like the ``discrete Gauss-Bonnet
 formula'' for surfaces triangulated by means of euclidean
 triangles. Via Hopf's index theorem, it will be also a consequence of
 the features of the tangent combings of $S$ defined below. Note that
 (2) implies (3) of Proposition \ref{N2_via_pb} (by considering $\partial
 T$), and that (3) applies in particular when $K=\partial T$ and $T$
 is an ideal triangulation of $\hat M$. 

\medskip

{\bf Singular combings.} Let $S$ be an oriented closed
surface as usual. We consider tangent vector fields $\vect \bv$ on $S$,
possibly having isolated zeros, where they locally
look like one of the following models, distinguished from each
other by the zero indices:
\begin{enumerate}
\item  The gradient of  $\pm (x^2+y^2)$, the index being
equal to $1$.
\item For every integer $n\geq 1$, consider the  $2n$-roots
  of unity $\alpha_1,-\alpha_1,\dots , \alpha_n,-\alpha_n$ in $\C\cong
  \R^2$. Let $y-a_j x=0$ be the equation of the straight line through
  $\alpha_j$ and $-\alpha_j$. Then the local model is the gradient of
  $\textstyle \pm(\prod_{j=1}^n (y-a_jx))$.  The zero index is equal
  to $1-n$.
\end{enumerate}

Every such a field $\vect \bv$ has a {\it type} given by the list of
its zero indices. The type satisfies the constraint given by the index
theorem, so that the sum of the indices equals $\chi(Z)$ for the restriction 
of the field on every component $Z$ of $S$. It is easy
to see that every tuple of integers $\leq 1$ satisfying this constraint is 
actually realized as the type of a field on $S$.

The fields $\vect \bv$ are considered up to homotopy through fields of
a given type (keeping the same notation). We denote by ${\rm Comb}(S)$
the set of equivalence classes, called {\it (singular)
  combings}. Clearly the combings are distributed by types. We say
that $\vect \bv$ is {\it taut} if all its zero indices are non
positive; we denote by ${\rm Comb}_\tau(S)\subset {\rm Comb}(S)$ the
subset of taut combings.

\begin{remark}\label{comb-inv}{\rm The {\it inversion} of every such a
    singular field $\vect \bv \leftrightarrow - \vect \bv$ preserves
    the type and induces the identity on ${\rm Comb}(S)$ as it is
    realized by the rotation by $\pi$ in every tangent plane, with
    respect to any auxiliary Riemann metric on the surface.}
 \end{remark}

 {\bf Realization via Abelian differentials.} For simplicity, assume again 
 that $S$ is connected. Give $S$ a structure of
 Riemann surface. Let $\phi$ be a {\it holomorphic Abelian
   differential} on $S$, with quadratic differential $\phi^2$. The
 horizontal measured foliation on $S$ defined by $\phi^2$ has
 orientable leaves. By fixing an orientation and using the field of
 oriented directions of the foliation we get a taut singular combing
 on $S$.  Every taut combing type can be obtained in this way (see
 \cite{MS}).

 \medskip
 
 {\bf On the structure of ${\rm Comb}(S)$.}
 Again for simplicity, assume that $S$ is connected. We have:
 
 \begin{prop}\label{affine-comb} The set ${\rm Comb}(S)$ has a partition 
 by subsets indexed by the singularity types, each
  subset being an affine space on $H_1(\Uu;\Z)$, where $\Uu$ is the
  complement of a system of disjoint $2$-disks in $S$ centred at the
  singular points. If $S$ is a torus, then ${\rm Comb}_\tau(S)$ is
  canonically identified to $H_1(S;\Z)$.
\end{prop}
\Dim  If $S$ is a torus, then every taut vector field on $S$ is non
singular. Given two such fields $\vect \bv_1$, $\vect \bv_2$, the
primary obstruction to determine the same combing is a
class $$\sigma(\vect \bv_1-\vect \bv_2)\in H_1(S;\Z)$$ such that
$\sigma(\vect \bv_1-\vect \bv_2)=-\sigma(\vect \bv_2-\vect \bv_1)$. As
$S0(2)\cong S^1$ and $\pi_2(S^1)=0$, this is in fact the complete
obstruction, so that ${\rm Comb}_\tau(S)$ is an affine space on
$H_1(S;\Z)$. Moreover, every oriented simple closed curve $s$ on $S$
determines an oriented foliation by parallel curves, and hence a field
$\vect \bv_s$. These special fields are all equivalent, and thus they
fix a base point in ${\rm Comb}_\tau(S)$, which is eventually
identified with $H_1(S;\Z)$.

In general, let $\vect \bv_1$, $\vect \bv_2$ be two (singular) vector
fields on $S$ with the same type of singularities. Up to isotopy, we
can assume that the two fields have the same zeros and coincide near
them. Then the complete obstruction to determine the same combing is a
class
$$\sigma(\vect \bv_1-\vect \bv_2)\in H_1(\Uu;\Z)$$
where $\Uu$ is as in the statement.\cvd 
\medskip

{\bf From triangulations to combings.} Given a branched triangulation
$(K, b)$ of $S$, the $1$-skeleton of the cell decomposition dual to
$K$ naturally carries a structure of {\it co-oriented} (hence oriented)
  {\it train-track} $\theta_b$ on $S$, by the following rule: {\it at every
  intersection point, an oriented edge of $(K,b)$ followed by the dual
  oriented branch of $(\theta,b)$ realize the orientation of $S$.}
 
Consider a regular neighbourhood $U$ (with smooth boundary) of
$\theta_b$. The closure of each component of $S \setminus U$ is a
$2$-disk $D_C$ contained in a $2$-cell $C$ of the cell decomposition
dual to $K$; this establishes a bijection between components of
$S\setminus U$ and $2$-cells. Indeed $U$ has a natural cellulation
made by the ``truncated triangles'' of $K$, and every component 
of $\partial U$ is the ``link'' of one vertex of $K$.
The neighborhood $U$ carries a 
tangent vector field $\vect \bv_b$ that is {\it traversing}, in the sense that:
\begin{enumerate}
\item Every integral line of $\vect \bv_b$ is a non degenerate closed
  interval which intersects transversely $\partial U$ at its
  endpoints. Generic integral lines are properly embedded into $(U,\partial U)$;
  \item There is a finite number of {\it exceptional} integral lines which are 
  simply tangent to $\partial U$ at a finite number of internal points.
 \end{enumerate}
 
 Moreover, $\vect \bv_b$ is {\it generic}, that is, every exceptional
 integral line has just {\it one} tangency point.
\smallskip
 
\begin{figure}[ht]
\begin{center}
 \includegraphics[width=7cm]{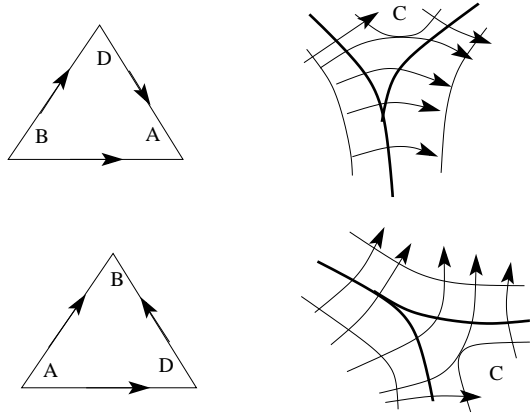}
\caption{\label{train-track} Train-tracks and combings.} 
\end{center}
\end{figure}

\begin{figure}[ht]
\begin{center}
 \includegraphics[width=7cm]{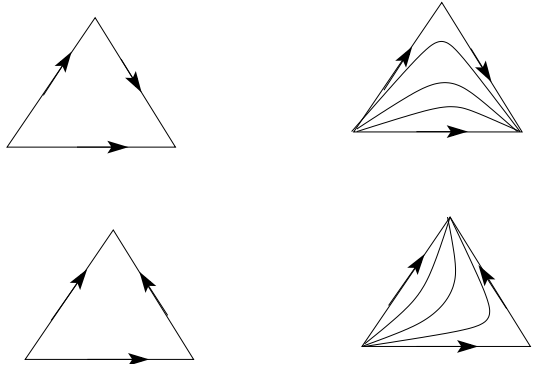}
\caption{\label{c-puzzle} Whitney fields on branched triangles.} 
\end{center}
\end{figure}

We easily realize that there is a tangency point in correspondence with
each $D$-colored corner of $(K,b)$, and that along every component of
$\partial U$ the exceptional integral lines occur with alternating
orientations. Hence the field $\vect \bv _b$ extends to a tangent
vector field (we keep the same name) on $S$ of the kind fixed above,
which has one zero of index equal to ${\rm Eu}_b(C)$ in the interior
of each component of $S \setminus U$, whenever ${\rm Eu}_b(C)\neq
0$. The traversing field is uniquely determined up to homotopy through
generic traversing fields, hence its singular completion is uniquely
defined up to homotopy through fields of the given type. So $(K,b)$
{\it carries a well defined combing}. A representative of $\vect
\bv_b$ can be obtained also as a {\it puzzle} where every tile is a
branched triangle equipped a classical Whitney field which can be
defined explicitely in terms of barycentric coordinates (see
\cite{HT}).  All this is illustrated in Figures \ref{train-track} and
\ref{c-puzzle}; in the first $C$ indicates a germ of $2$-cell $C$
where the field contributes to ${\rm Eu}_b(C)$.

If the locally taut $\mz/2$-charge $\gamma_b$ on $(K,b)$ is taut, then
the combing $\vect \bv_b$ is taut; if moreover all the components of
$S$ are tori, then $\vect \bv_b$ is a non-singular combing.

\subsection{$2D$ transits and intrinsic structures} We consider local
moves on naked triangulations of the surface $S$.  These are the
$2\leftrightarrow 2$ ``diagonal exchange'' move, also called {\it flip}, 
the {\it bubble} $0\leftrightarrow 2$ move and the
$1\leftrightarrow 3$ move. The last one can be obtained as a
concatenation of a bubble move and a flip, but it is convenient to
consider it by itself. Flips preserve the number of vertices while the
other two positive moves increase it by $1$.  As in Section
\ref{seccharges}, all these moves can be naturally enhanced to
transits of 2D charges $(K,c) \leftrightarrow (K',c')$; the rule is
that $c$ and $c'$ coincide on every common triangle of $K$ and $K'$,
and that, by restricting them to the subcomplex supporting the move, $c$ and $c'$
satisfy the local charge condition (i) on each abstract triangle, the global one (ii)
around each internal vertex, and the total charge at every boundary vertex
is preserved. For instance, a {\it taut} transit of 2D charges is such
that the taut condition is satisfied at the locus of the move. So,
similarly to the $3$-dimensional case, we can define the quotient sets
$$c(S,\Z)\ , \ c(S,\Z/2)\ $$ of  $\mz$-charged and 
$\mz/2$-charged triangulations of $S$ up to charge transits. In the two dimensional
setting, the result analogous to Proposition \ref{Z2cmeaning} (and
easier to prove) is: 
\begin{prop} The set $c(S,\Z)$ encodes
$H^1(S;\Z)$, while $c(S,\Z/2)$ encodes $H^1(S;\Z/2)$.
\end{prop}
 
The $2D$ branching transits (``$b$-transits" for short) $(K,b) \to (K',b')$ are defined by
imposing that the orientation is preserved on every common edge of $K$ and $K'$. We classify now the $b$-transits with respect to the associated
combings on $S$. This is reminiscent of the study of $3D$ branched
spines \cite{BP2}.

\medskip
 
{\bf Branched flips.}  Given a branched triangulation $(K,b)$ and a
naked flip $K\rightarrow K'$, there always exists a $b$-transit
$(K,b)\rightarrow (K',b')$, called {\it branched flip} or {\it
  $b$-flip}.  A branched flip $(K,b) \rightarrow (K',b')$ is {\it
  forced} if it is the unique one supported by the naked flip
$K\rightarrow K'$ and starting with $(K,b)$. 

\begin{figure}[ht]
\begin{center}
 \includegraphics[width=9.5cm]{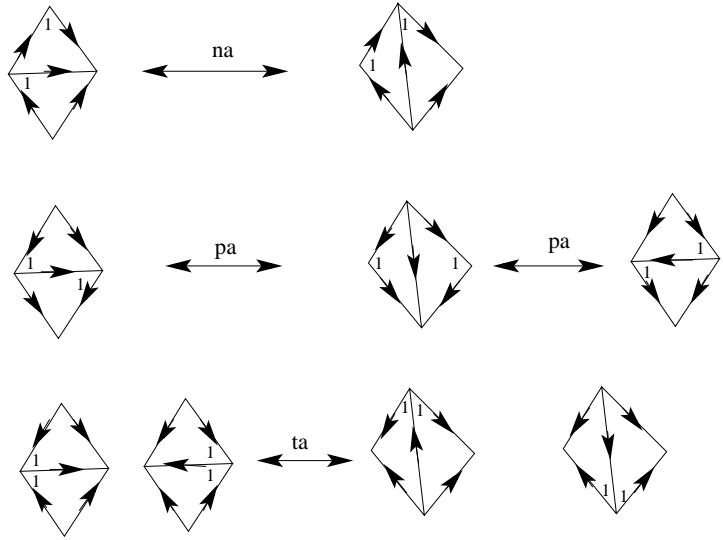}
\caption{\label{b-flip} Branched flips.} 
\end{center}
\end{figure}

\begin{figure}[ht]
\begin{center}
 \includegraphics[width=9cm]{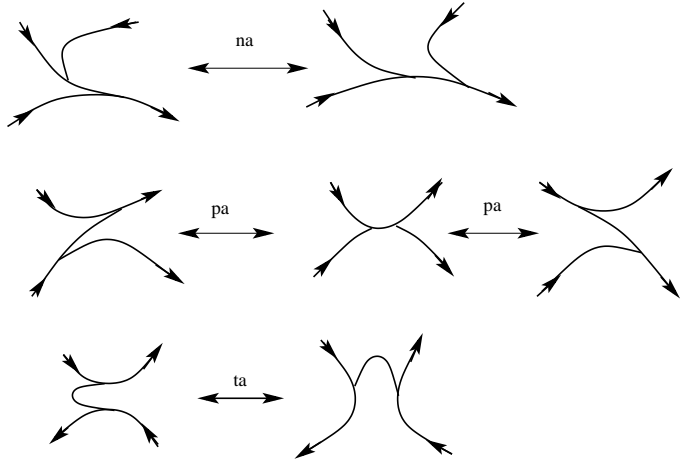}
\caption{\label{sliding_flip} Sliding and bump branched 
track flips.} 
\end{center}
\end{figure}

The branched flips
$(K,b)\rightarrow (K',b')$ are distributed in the following classes,
illustrated in Figures \ref{b-flip} and \ref{sliding_flip}. In Figure
\ref{b-flip} we have labelled by ``$1$'' the $D$-colored corners; we
will do the same in the next figures. In Figure \ref{sliding_flip} we
show dual train tracks; according to our orientation convention, they are obtained by total inversion of the branchings shown in Figure \ref{b-flip}.  Note that {\it the flip classification below is invariant under total inversion}:
 
  \begin{enumerate}
  \item {\it Non ambiguous}, such that $(K,b)\rightarrow (K',b')$ and 
the inverse 
  $b$-flip $(K,b)\leftarrow (K',b')$ are forced.
 \item {\it Forced ambiguous}, such that $(K,b)\rightarrow (K',b')$
 is forced but the inverse $b$-flip is not.
 \item {\it Sliding}, such that at least one among  
 $(K,b)\rightarrow (K',b')$ and its inverse is forced. 
\item {\it Bump}, such that both $(K,b)\rightarrow (K',b')$
and its inverse are not forced. 
\end{enumerate}

\begin{figure}[ht]
\begin{center}
 \includegraphics[width=7cm]{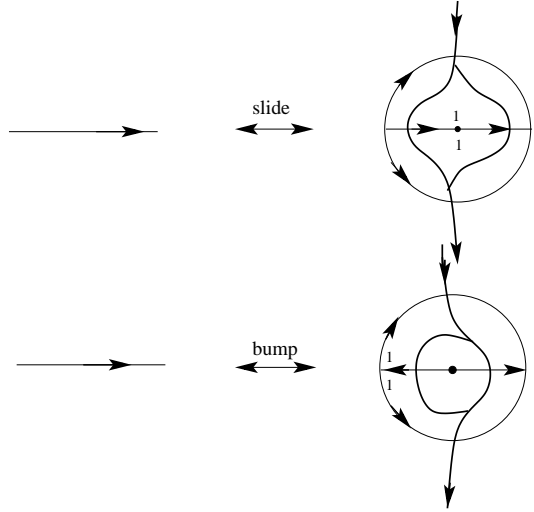}
\caption{\label{2Dbubble} Sliding and bump bubble transit.} 
\end{center}
\end{figure}

{\bf Bubble $b$-transits.} The (positive) 2D bubble $b$-transits are
distributed in two classes (see Figure \ref{2Dbubble}):
 \begin{enumerate}
 \item {\it Sliding}, such that the two corners at the central new
   vertex after the positive transit are $D$-colored.
 \item {\it Bump}, such that one central corner is $A$-colored,
   and the other is $B$-colored.
\end{enumerate}

\begin{figure}[ht]
\begin{center}
 \includegraphics[width=7cm]{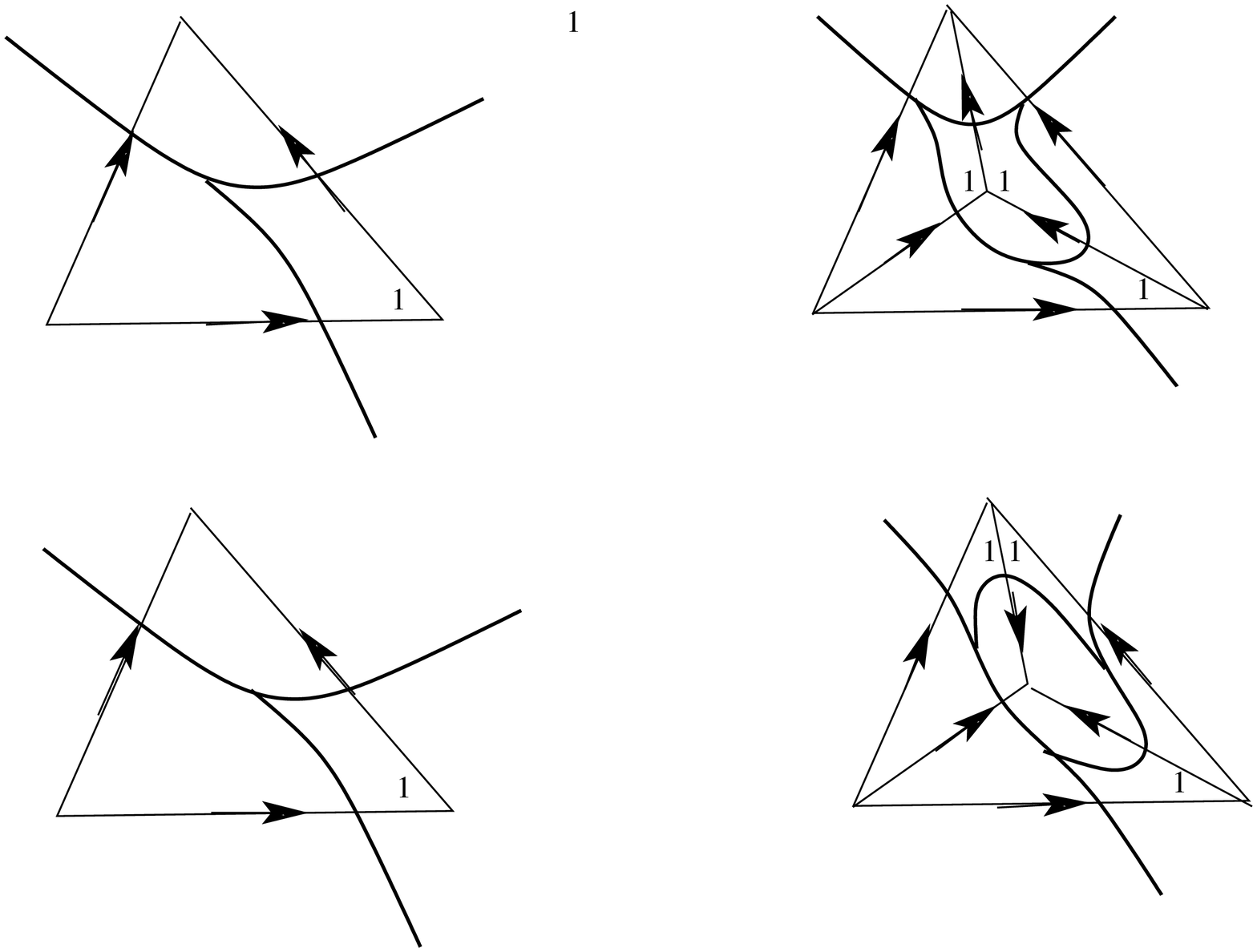}
\caption{\label{1to3} Sliding and bump $1\leftrightarrow 3$ transit.} 
\end{center}
\end{figure}

{\bf $1\leftrightarrow 3$ $b$-transits.} These are distributed in two
classes (see Figure \ref{1to3}):
\begin{enumerate}
\item {\it Sliding}, such that two corners at the central new vertex
  after the positive transit are $D$-colored, while the other corner
  can be either $A$- or $B$-colored.
\item {\it Bump}, such that no central corner is $D$-colored, the
  central vertex is not monochromatic, and either $A$ or $B$ can be
  the prevailing color.
\end{enumerate}

Recall the Euler $2$-cochain ${\rm Eu}_b^-$ defined in Section \ref{train-trak}, 
and denote by ${\rm Eu}_b^-(Z)$ its evaluation
on the fundamental class of a component $Z$ of $S$.

\begin{prop}\label{on2Dflip} 
  (a) Every $b$-transit $(K,b) \rightarrow (K',b')$ induces a locally
  taut $\Z/2$-charge transit $(K,\gamma_b)\rightarrow
  (K',\gamma_{b'})$.

\smallskip
  
  (b) For every $b$-transit
  $(K,b)\rightarrow (K',b')$, the following properties are
  equivalent:
\begin{enumerate}
\item It preserves the combing on $S$.
\item ${\rm Eu}_b^-(Z)={\rm Eu}_{b'} ^-(Z)$
on every  component $Z$ of $S$.
\item $(K,\gamma_b)\rightarrow (K',\gamma_{b'})$ is actually a
taut $\Z/2$-charge transit.
\item It is a sliding $b$-transit.
\end{enumerate}
\end{prop}
\Dim As above it is convenient to represent the combings carried by a
branched triangulation by the singular completions of suitable generic
traversing fields. By analyzing the local modifications of the
traversing fields supported by the transits of train tracks in Figure
\ref{sliding_flip}, we easily realize that before and after a sliding
$b$-flip we deal with fields defined on a same neighborhood $U$ and
homotopic through (not necessarily generic) traversing
fields. Precisely, we can construct a generic homotopy such that just
one traversing field is not generic, as it has just one exceptional
integral interval, which is simply tangent at two points of $\partial
U$. Then the singular completions define a same combing.  In the case
of a positive sliding bubble $b$-transit, up to homotopy the
traversing field after the move is the restriction of the traversing
field before the move, and we readily realize that again the singular
completions define the same combing. Similarly for a positive sliding
$1\rightarrow 3$ $b$-transit. On the other hand, any positive bump
$b$-transit introduces a new singular point of the combing of index
$1$. The Proposition straighforwardly follows from these
considerations.\cvd 

\medskip

Let us denote by $B_s(S)$ the quotient set of branched triangulations
of $S$ up to the relation generated by isotopy and the sliding
transits.  Similarly let $B_\tau(S,\Z/2)\subset B_s(S)$ be formed by
the classes represented by $\Z/2$-taut triangulations. By Proposition
\ref{on2Dflip}, the correspondence
$$(K,b)\mapsto \vect \bv_b$$
induces a well defined map
$$ \kappa: B_s(S) \rightarrow {\rm Comb}(S)$$
which factorizes through the set $B_s(S)/\pm$ obtained by adding to the
generators of the relation the branching inversion $b\to -b$. By
restricting $\kappa$ we have also a map
$$\kappa_\tau: B_\tau(S,\Z/2) \rightarrow {\rm Comb}_\tau(S)$$ 
which factorizes through   $B_\tau(S,\Z/2)/\pm$. 
\begin{teo}\label{sliding-combing} 
  The map $\kappa: B_s(S)/\pm \to {\rm Comb}(S)$ is bijective.
  Similarly for the restriction to $B_\tau(S,\Z/2)/\pm$.
\end{teo}

\Dim Refering to the proof of Proposition \ref{on2Dflip}, the main points are:
\begin{enumerate}
\item Every (taut) combing on $S$ can be realized by the vector field
$\vect \bv_b$ carried by some ($\Z/2$-taut) triangulation $(K,b)$
(that is, the maps $\kappa$ and $\kappa_\tau$ are onto).

\item Two triangulations $(K,b)$ and $(K',b')$ with the same set of vertices
carry generic traversing fields (defined on the same $U$) which are homotopic
through (non necessarily generic) traversing fields if and only if
they are related by a finite sequence of sliding $b$-flips.  
\end{enumerate}

Both facts follow from simplified versions of the arguments used
in \cite{BP2} and \cite{BP3} for the treatment of combings on
$3$-manifolds via branched spines. Let us indicate the main
ideas. The proof of (1) is based on Ishii's notion of {\it flow
  spines} (\cite{Is}); a detailed proof in $3D$ is given in Chapter 5
of \cite{BP2} in the case of closed manifolds, and in \cite{BP3} it is 
extended to manifolds with boundary.  As for (2), one
implication has been already remarked in Proposition \ref{on2Dflip};
for the other implication, by transversality we can assume that the
homotopy is generic, that is, it contains only a finite number of non
generic traversing fields, each one containing one exceptional
integral interval which is tangent at two points of $\partial U$. Then
we have to analyze how two generic traversing fields close to a non
generic one are related to each other. Finally one realizes that the
sliding $b$-flips cover all possible configurations.

Given two triangulations that carry the same combing, we can modify the 
one with fewer vertices by isotopy and a finite number of positive sliding 
bubble $b$-transits in order that the two resulting triangulations verify the
hypothesis of (2). By using these facts, taking into account the branching inversion, 
the injectivity of the maps readily follows.\cvd

\subsection{Non negative $1$-cycles and their transits}\label{1-cycles}
Let as usual $(K,b)$ be a branched triangulation of a closed oriented
surface $S$. The $1$-skeleton $K^{(1)}$ has oriented edges. A {\it non
  negative $1$-cycle on $(K,b)$} is a simplicial $\Z$-cycle
$\textstyle \gamma = \sum_{e\in K^{(1)}} a_e e$ such that $a_e \geq 0$
for every edge $e$.  Denote by $H^+(K,b)$ the image of the set of non
negative cycles in $H_1(S;\Z)$.  We want to point out the behaviour of
$H^+(K,b)$ under $2D$ sliding transits. We have:
\begin{prop}\label{transit-1cycle}
(1) If $(K,b)\to (K',b')$ is a non ambiguous flip or a sliding bubble 
$b$-transit
then $H^+(K,b)=H^+(K',b')$.

(2)  If $(K',b')\to (K,b)$ is a forced ambiguous $b$-flip, 
then $H^+(K,b)\subset H^+(K',b')$ and in general the inclusion
is strict.
\end{prop}

\begin{figure}[ht]
\begin{center}
 \includegraphics[width=5cm]{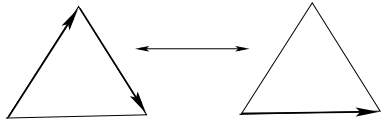}
\caption{\label{cycle-move} $1$-cycle move.} 
\end{center}
\end{figure} 

The proof consists in analyzing the possible local configurations of a
non negative cycle $\gamma$ of $(K,b)$, and to verify whether $\gamma$
transits in a unique way to a non negative cycle $\gamma'$ on
$(K',b')$ in the same homology class, up to the local move suggested
by Figure \ref{cycle-move} (which preserves the homology class). The
only situation which gives rise to a stop is when $(K',b')\to (K,b)$
is a forced ambiguous $b$-flip, and when the relevant local portion of $\gamma'$
is supported by the edge which is flipped
to produce $(K,b)$. This can be expressed also in terms of dual
measures on the train track $\theta_b$; we will spell it for flips in
the next Remark. So we have the {\it negative} conclusion that $H^+(K,b)$
is not in general a $2D$ sliding move invariant. On the other hand, we
will see that things go better in the $3D$-fillable situation (see
Section \ref{2-cycles}).

\begin{remarks}\label{partialA}{\rm 
Proposition \ref{transit-1cycle} indicates an interesting
difference between non ambiguous and forced ambiguous
flips.  We can better understand this difference dually in terms
of the measures carried by the train tracks $\theta_b$. Recall
that a {\it measure} $\mu$ on $\theta_b$ assigns to every edge a
real non negative weight in such a way that at every $3$-valent
vertex of $\theta_b$, the natural ``switching condition'' is
satisfied. As $\theta_b$ is oriented and is a spine of its regular
neighbourhood $U$, the measures on $\theta_b$ actually form a {\it
positive cone} of $H_1(U;\R)$ (recall that every $1$-homology
class of $U$ is {\it uniquely} represented by a real $1$-cycle on
$\theta_b$ which assigns to every edge a real weight and verifies
the same switching condition at vertices).

\begin{figure}[ht]
\begin{center}
 \includegraphics[width=8cm]{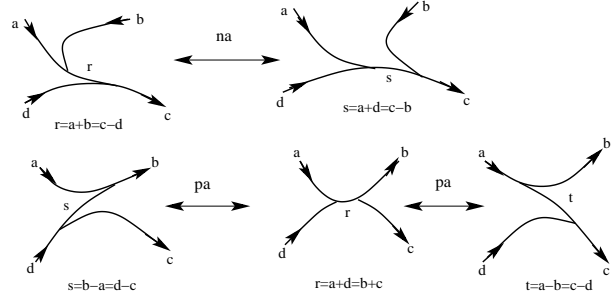}
\caption{\label{measured_flip} Measured flips.} 
\end{center}
\end{figure}    

Assume that $(K,b)$ is $\Z/2$-taut and that $\theta_b$ carries a
nowhere vanishing measure, say a {\it plain measure}.  In Figure
\ref{measured_flip}, top line, we see the transit of such a measure
supported by a typical non ambiguous flip. We realize that the {\it
inequalities} $c>d$ and $c>b$ are necessarily satisfied, and that
there is a {\it bijection} of the weights. In the middle of the bottom
line we see the result of two forced ambiguous flips related to each other.
 We realize that in the middle the {\it
  equality} $a+d=b+c$ is necessarily satisfied, and no further
inequalities must be imposed. On the other hand, on the left
(resp. right) side the inequalities $b>a$, $d>c$ (resp. $a>b$, $c>d$)
are necessarily satisfied and the transit {\it injects} the left
(resp. right) side set of weights onto a subset, say $W_L$
(resp. $W_R$), of the middle one. We easily see that $W_L$ and $W_R$
form a partition of this last set. In order to reverse the transits we
have to restrict to $W_L$ and $W_R$ respectively. Hence in a sense the
measure inequalities solve the partial ambiguities.
\smallskip

There is a natural $2D$ notion of {\it veering triangulation} $(K,b)$
of $S$: every vertex of $K$ must be either $A$- or
$B$-monochromatic. Clearly, if a $3D$ triangulation $(T,\omega)$ of
$\hat M$ is veering, then $(\partial T,\partial \omega)$ is
veering. We note that {\it the non ambiguous flips preserve this
  property}, while the forced ambiguous flips and the sliding bubble
$b$-transits do not. As mentioned in Remark \ref{vering}, the $3D$ non
ambiguous transits do not preserve the veering property; this reflects
in the fact that the associated systems of entangled $2D$ transits
(described in Proposition \ref{first_inv} below) involve all kinds of
sliding transits, not only non ambiguous flips.}
\end{remarks}
 
\subsection{The boundary maps}\label{hydden2}
\begin{figure}[ht]
\begin{center}
 \includegraphics[width=6cm]{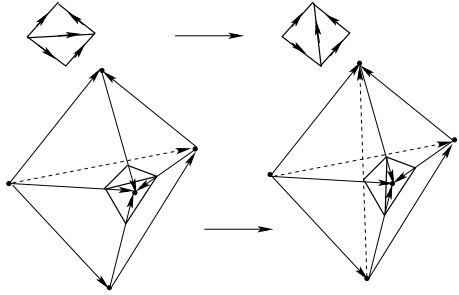}
\caption{\label{betaflip} From $3D$ non ambiguous towards  $2D$ 
sliding transits.} 
\end{center}
\end{figure}

Let us go back to dimension three. Let $(T,\omega)$ be a pre-branched
ideal triangulation of $\hat M$, and consider the boundary branched
triangulation $(\partial T, \partial \omega)$ of $\partial M$. Every
abstract tetrahedron of $T$ carries four triangles of $\partial
T$. Consider any naked ideal $3D$ transit $T\rightarrow T'$; it gives
rise to a system of $2D$ transits on $\partial T$. More precisely:
\begin{enumerate} 
\item Every naked positive $2\rightarrow 3$ move gives rise to three naked
  flips and two positive $1\rightarrow 3$ moves.
\item Every naked positive lune move gives rise to two positive $2D$
  bubble moves.
\item Every $3D$ pre-branching transit induces a system of $2D$
  branched transits, supported by the associated system of naked moves.
\end{enumerate}
\begin{prop}\label{first_inv} 
  Every $3D$ non ambiguous ideal $pb$-transit gives rise to a system of $
  2D$ sliding transits. Hence we have well defined maps $\partial :
  \mathcal{NA}^{id}(M) \rightarrow B_s(\partial M)$ and $\kappa
  \circ \partial : \mathcal{NA}^{id}(M) \rightarrow {\rm
    Comb}(\partial M)$, which restricts to $ \kappa \circ \partial :
  \tau(M,\Z/2) \rightarrow {\rm Comb}_\tau(\partial M)$. 
\end{prop}    
\Dim The proof goes through a direct analysis of all transits. In fact,
by using Lemma \ref{wb}, concerning the $2\leftrightarrow 3$ non ambiguous
$pb$-transits it is enough to study the two ones dominated by the remarkable 
$b$-transits indicated before the statement of Lemma \ref{wb}. 
For example, refering to the positive non ambiguous transit of Figure
\ref{betaflip}, we see that two of the associated flips are non
ambiguous, one is forced ambiguous, and the $1\to 3$ $2D$ transits are
sliding.  The other remarkable $2\rightarrow 3$ $b$-transit as well as
the non ambiguous lune $b$-transits have similar behaviour.\cvd
 \medskip
 
In this way we have obtained a geometric topological invariant for
$\mathcal{NA}^{id}(M)$ which lives on the boundary of $M$.
\medskip

The following Conjecture (perhaps better qualified as a ``question'')
sounds attractive and non trivial. Recall the natural projection 
$\pi: \Nn\Aa^{id}(M)\to \Pp\Bb^{id}(M)$.

\begin{conj} \label{Ricconj}{\rm For every ideal pre-branched triangulations $(T,\omega)$ 
and $(T',\omega')$ of $\hat M$, we have
$$[(T,\omega)]=[(T',\omega')]\in \Nn\Aa^{id}(M)$$ 
if and only if
$$\pi([T,\omega])=\pi([T',\omega']) \  {\rm and} \  
\kappa \circ \partial([T,\omega]) =  
\kappa \circ \partial([T',\omega']) \ . $$}
\end{conj} 

We are going to point out a few further invariants of non ambiguous
strucures.  
  
\subsection{Non negative $2$-cycles and their
  boundary}\label{2-cycles}
Let $(T,\omega)$ be a pre-branched ideal triangulation of $\hat M$.
Consider as usual the cell decomposition of $M$ formed by the
truncated tetrahedra of $T$. Every truncated $2$-face of $T$ is a
hexagon with edges alternatively on $\partial M$ or in the interior of
$M$. Recall that every hexagon $H$ is oriented by $\omega$.

A {\it non negative $2$-cycle} is a cellular relative $2$-cycle on
$(M,\partial M)$ of the form $\textstyle \Gamma= \sum_{H} c_H H$ such
that every coefficient $c_H \geq 0$.  The boundary $\partial \Gamma$
is a non negative $1$-cycle on $(\partial T, \partial \omega)$ (see
Section \ref{1-cycles}), which determines completely $\Gamma$. Denote
by $H^+(T,\omega)\times \partial H^+(\partial T, \partial \omega)$ the
subset of $H_2(M,\partial M;\Z) \times H_1(\partial M;\Z)$ determined
by the set of these pairs of cycles $(\Gamma,\partial \Gamma)$. The
key remark is that, in the system of boundary 2D sliding transits
associated to an ideal non ambiguous $3D$-transit, a boundary cycle
$\partial \Gamma$ never falls in the configuration that gives rise to
a stop, as described after Proposition \ref{transit-1cycle}. Then a
$2D$ $1$-cycle transit can be completed to a $3D$ $2$-cycle transit.
So we have:
\begin{prop}\label{invariant H+}  If $(T,\omega)$ and $(T',\omega')$ 
  represent the same non ambiguous structure, then the semigroups 
  $H^+(T,\omega)\times \partial H^+(\partial T, \partial \omega)$ and
  $H^+(T',\omega')\times \partial H^+(\partial T', \partial \omega')$
  are isomorphic.
\end{prop}

\begin{figure}[ht]
\begin{center}
 \includegraphics[width=3cm]{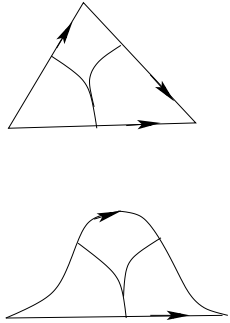}
\caption{\label{transv-track} Transverse tracks.} 
\end{center}
\end{figure}  

\begin{remark}\label{on taut}{\rm When the triangulation $(T,\omega)$
    is taut the union of the hexagons has a natural structure of
    branched surface called the {\it branched surface of the taut
      triangulation}. Its boundary is an oriented train track on
    $\partial M$, which is positively transverse to the one on
    $(\partial T, \partial \omega)$ defined in Section
    \ref{train-trak}. They determine equivalent combings.
The basic tiles of such tracks are shown in Figure \ref{transv-track}}.
\end{remark}   

\subsection {Invariant coloured arc-links }\label{inv_arclink}
Let $\sG=[(T,\omega)]\in \Nn\Aa^{id}(M)$. Every truncated edge of
$T$ is a properly embedded simple arc $\lambda$ in $(M,\partial M)$
joining two vertices of $\partial T$, say $v$, $v'$; they belong to two
regions $C$, $C'$ of $\partial M \setminus \theta_{\partial \omega}$
such that ${\rm Eu}_{\partial \omega}(C)={\rm Eu}_{\partial
  \omega}(C')$. We give the arc the same label $l(\lambda):= {\rm
  Eu}_{\partial \omega}(C)$. The union of such arcs with {\it non zero
  labels} forms a $\Z$-{\it colored link $L_\omega$} of properly
embedded arcs in $(M,\partial M)$.  If $\lambda \subset L_\omega$, then
we can assume that $v$, $v'$ are the singular points of a boundary
combing representative $\vect \bv_{\partial \omega}$. Clearly $L_\omega =
\emptyset$ if $(T,\omega)$ is taut. Let us denote by $\Ll(M,\partial M,\Z)$
the set of isotopy classes of $\Z$-colored links of properly embedded arcs
on $(M,\partial M)$. It is immediate that the non ambiguous transits
preserve the isotopy class of  $L_\omega$. Then we have:
\begin{prop}\label{arc-link}
The correspondence $\sG \mapsto [L_\omega]$, for every
  $\sG=[(T,\omega)]\in \Nn\Aa^{id}(M)$, well defines a map  
$ \Lambda: \Nn\Aa^{id}(M)\to \Ll(M,\partial M,\Z)$.
\end{prop}

A natural {\it realization problem} is to determine the image of the
map $\Lambda$. Immediate obstructions can be derived from Lemma
\ref{EU}. 

\section {Taut structures}\label{TAUT} 
In this section we consider in detail the set of taut structures
$\tau(M)$. We assume that $\partial M$ is a non empty collection of tori. 
The aim of this section is to prove Theorem \ref{existNA}. Let us recall 
the main qualitative
results of \cite{L}:
\begin{prop} {\rm \cite[Proposition 10 and Theorem 1]{L}} If $\hat M$
  admits a taut triangulation, then $M$ is irreducible and the
  boundary tori are incompressible. Moreover, every such a manifold $M$
  which is an-annular admits taut triangulations.
\end{prop}

Note that for hyperbolic cusped manifolds the sufficient existence
conditions of taut triangulations are satisfied. Probably the most
remarkable property of taut triangulations is:

\begin{prop} {\rm \cite[Theorem 3]{L}} Every compact surface $Z$
  properly embedded in $M$ and carried by a branched surface of a taut
  triangulation realizes the Thurston's norm of its homology class in
  $H_2(M,\partial M;\Z)$.
\end{prop}

As for the existence of taut structures, we consider first the simpler
case of a mapping torus. Then we will describe a procedure depending
on an initial choice of sutured manifold hierarchy of $M$.

\subsection{Manifolds fibering over $S^1$}\label{FIBER} 
Let $M$ be a compact oriented $3$-manifold bounded by tori and which
fibres over $S^1$ with fibers of negative Euler characteristic. So
there is a mapping torus realization $M \cong (Z\times [0,1])/\psi$,
where $Z$ is an oriented compact connected surface with non empty
boundary, $\chi(Z)<0$, and $\psi$ is an orientation preserving
diffeomorphism of $Z$. The surface $\hat Z$ has ideal
triangulations. Then, following \cite{L} we can construct taut
triangulations $(T,\omega)$ of $\hat M$ by implementing the following
procedure:
\begin{itemize} 
\item {\it Choose} an ideal triangulation $\Ss$ of $\hat Z$;
\item {\it Choose} a finite sequence of flips $ \Ss:=\Ss_0 \to \Ss_1
  \to \Ss_2 \to \dots \Ss_k:= \psi(\Ss_0)$ connecting $\Ss$ to
  $\psi(\Ss)$; if some edge of $\Ss$ is left unchanged by the sequence, add to it consecutively a flip and an inverse flip at the edge. 
\item Consider the four triangles involved in a flip $\Ss_j \to \Ss_{j+1}$ as the boundary of a pre-branched tetrahedron $(\Delta_j,\omega_j)$, so that the two edges exchanged by 
  the flip are its diagonal edges, and the co-orientations of the triangles of $\Ss_j$ (resp. $\Ss_{j+1}$) are ingoing (resp. outgoing) $\Delta_j$. Then we get an ideal triangulation $\tilde T$ of the space $\hat Z \times [0,1]$, with boundary triangulations $\Ss$ and $\psi(\Ss)$. Define $(T,\omega)$ as the projection of $\tilde T$ to $ \hat M=(\hat Z \times [0,1])/\psi$.
\end{itemize}

We call {\it layered triangulation} any taut triangulation of $\hat M$ obtained in this way. The first two steps of the above construction contain arbitrary choices, but we are going to see that the layered triangulations of a given mapping torus realization of $M$ define nevertheless the same taut structure. At this point we need to recall some fundamental results of Thurston \cite{Th} and Fried \cite{F}.

Let $W$ be a compact connected oriented irreducible $3$-manifold which fibers over $S^1$. Denote by $B_W$ the unit ball of the Thurston norm of $W$, $\vert\vert \cdot \vert\vert : H^1(W;\mz)\rightarrow [0;+\infty[$. The fibrations of $W$ are in $1$-to-$1$ correspondence with the integral points of a union of cones over some open top dimensional faces $F_W$ of $B_W$, called its {\it fibred faces} \cite{Th}. There are flows $(\psi_t)$ of diffeomorphisms of $W$ which are positively transverse to any fibration over a same fibred face (that is, $d(f\circ \psi_t(x))/dt >0$ for every $x\in W$ and every fibration $f:W\rightarrow S^1$ such that
$f^*[d\theta]$ is in the cone $C(F_W)$ over $F_W$, for some open fibred face $F_W$). Such a flow $(\psi_t)$ determines the fibered face by the condition that the normal plane bundle of the vector field $d(\psi_t)/dt$ has Euler class $$\chi_{(\psi)_t}(u) = -\vert\vert u \vert\vert$$ for all $u$ in the cone over the open face (here we view $\chi_{(\psi)_t}$ in $H^2(W,\partial W;\mz) \cong H_1(W;\mz)$ as a linear functional on $H^1(W;\mz)$). Moreover, if the interior of $W$ has a finite volume complete hyperbolic structure, then every fibered face $F_W$ of $B_W$ determines a unique isotopy class of {\it pseudo-Anosov} flows $(\psi_t)$ inducing a pseudo-Anosov return map on any fiber of any fibration over $F_W$ (see Theorem 7 of \cite{F}; the results of that paper are formulated for a closed compact $W$, but the arguments work verbatim when $W$ is a cusped manifolds). 

\begin{prop} \label{fibercase1} (1) For every realization of $M$ as a mapping torus $M_\psi:=(Z\times [0,1])/\psi$ such that $\chi(Z)<0$, all layered triangulations of $M_\psi$ determine the same element $\sG_\psi$ of $\tau(M)$. 

(2) For any two mapping torus realizations $M_\psi$ and $M_\phi$ of $M$ corresponding to fibrations lying on a same ray from the origin of $H_2(M,\partial M;\R)$, we have $\sG_\phi = \sG_\psi$.
  
(3) Any two mapping torus realizations $M_\psi$ and $M_\phi$ of $M$ such that $\sG_\phi = \sG_\psi$ lie over the same fibered face of the Thurston ball of $M$. 
\end{prop}
  \Dim (1) The Ptolemy groupoid of the set of (naked) ideal triangulations of $\hat Z$
  is generated by the flips modulo the ``square'' and ``pentagon''
  relations. Then it is enough to check how these relations modify
  the construction of layered triangulations of $M_\psi$. One
  realizes easily that the square relation gives rise to a non ambiguous lune
  move, and the pentagon relation gives rise to a non ambiguous
  $2\leftrightarrow 3$ transit. 
  
(2) It is well-known that any element $a\in H_2(M,\partial M;\Z)$ is
represented by an oriented and properly embedded surface $S$ in $M$, and that, if
$a=kb$ with $b$ a primitive element, then $S$ is the union of $k$
connected components, each representing $b$ (see eg. \cite{Th}, Lemma
1). In particular, let $a$ and $b$ be the classes of (the fibers of)
$M_\psi$ and $M_\phi$, respectively. Fibers are incompressible, and
moreover, since in any fibered $3$-manifold any incompressible surface
in the homology class of a fiber is isotopic to a fiber (\cite{Th},
Theorem 4), the components representing $b$ are isotopic. Hence, given a layered
triangulation $T$ of $M_\phi$, one obtains a layered triangulation of
$M_\psi$ by gluing $k$ copies of $\tilde T$ (ie. $T$ cut along $\hat Z$) along $k-1$ boundary components by the identity map, and then gluing
back the ends by $\phi$. Clearly, it is in the same non ambiguous class as $T$.

(3) Let $T_\psi$ be a layered triangulation of $M_\psi$. Denote by $\Sigma_\psi$ the underlying branched oriented surface (see Remark \ref{on taut}). The set of real $2$-cycles on $\Sigma_\psi$ surjects on $H_2(M,\partial M;\mz)$, and the subset $H^+(T_\psi)$ of non
negative $2$-cycles is a piecewise linear rational cone in $H_2(M,\partial M;\mz)$. The interior $C_\psi$ of $H^+(T_\psi)$ consists of the {\it full} $2$-cycles (with all positive coefficients). We claim that it is non empty. Indeed, any fiber $Z$ of $M_\psi$ is represented by a triangulated copy on $\Sigma_\psi$, that is, by a non negative $2$-cycle with coefficients $0$ or $1$. Then, by doing flips (ie. going through some tetrahedra of $T_\psi$), one realizes new triangulated copies which may have some triangles in common; all together they eventually determine a full $2$-cycle, which is homologically an integer multiple of $Z$. Conversely, any full $2$-cycle on $\Sigma_\psi$ is dual to an integral non singular closed $1$-form on $M$, and hence is a fiber \cite{Th} (this follows also from the fact that the complement of a full $2$-cycle on $\Sigma_\psi$ is a disjoint union of products). Then, any flow $(\psi_t)$ on $M$ positively transverse to $Z$ is positively transverse to $\Sigma_\psi$, and hence to any surface represented by an element of $H^+(T_\psi)$. Let $T_\phi$ be a layered triangulation of $M_\phi$. Since $\sG_\psi=\sG_\phi$ we have $H^+(T_\psi)\cong H^+(T_\phi)$ (Proposition \ref{invariant H+}), so $(\psi_t)$ is positively transverse to the fibers of $M_\phi$. Then, by Fried's result $M_\phi$ and $M_\psi$ must lie in the cone over a same fibered face.\cvd


\medskip

This achieves point (1) in Theorem \ref{existNA}.

\medskip

It can happen that mapping tori corresponding to different fibrations
of $M$, even with non homeomorphic fibres, have nevertheless a common
layered triangulation. Hence they determine the same taut structure. By Theorem \ref{Agol}, for many fibered cusped manifolds $M$ this happens over the whole of some fibered faces of the Thurston ball. We are going to discuss the proof, by following Agol's arguments \cite{A2}.

Let us use again the notations introduced before Theorem \ref{fibercase1}. Assume that the interior of $W$ has a finite volume complete hyperbolic structure. The monodromy of a fibration $f$ in $C(F_W)$ is given by the isotopy
class of the first return map $R_{(\psi_t),S}:S \rightarrow S$,
$x\mapsto \psi_{t(x)}(x)$, where $(\psi_t)$ is the pseudo-Anosov flow associated to $F_W$, $S$ is a fiber of $f$, and $t(x)>0$ is
the smallest time such that $\psi_{t(x)}(x)\in S$. Denote by $\phi_S$ the (unique) pseudo-Anosov homeomorphism of $S$ isotopic to $R_{(\psi_t),S}$. The suspension of the singular points
of $\phi_S$ is a link in $W= S\times [0,1]/\phi_S$ transverse to the
fibers. By varying $S$ or the fibration $f$ in $C(F_W)$, one obtains
isotopic links. Hence, the cusped manifold $M$ obtained from $W$ by
drilling out a small regular neighborhood of these links is uniquely
determined by the face $F_W$, and by restricting to $M$ the fibrations
$f$, one obtains fibrations of $M$ lying in a determined face $F_M$ of
the Thurston ball $B_M$ of $M$. The layered triangulations of $M$
associated to a fiber $Z$ of a fibration in the open cone $C(F_M)$
have pre-branching co-orientations which are compatible with the flow
induced by $(\psi_t)$ on $M$.

Recall the definition of a
veering triangulation (see Remark \ref{vering}). The proof of Theorem \ref{Agol} follows from:
\begin{prop}\label{fibercase2} There is a unique triangulation of $\hat M$ which is veering and layered for any fibration over
  the face $F_M$.
\end{prop}

\Dim Consider a fiber $S$ of a fibration $f$ in $C(F_W)$ as above, and
the pseudo-Anosov homeomorphism $\phi_S: S\rightarrow S$. Let
$(\theta,\mu)$ be any measured train track on $S$ carrying the stable
foliation $L^s_{\phi_S}$ of $\phi_S$. Consider the sequence of {\it
  maximal splittings} starting at $(\theta,\mu)$, that is, the
sequence of positive forced ambiguous flips at edges of maximal
$\mu$-weights (these are the flips going from the left or right picture to the
middle one on the bottom of Figure \ref{measured_flip}; the sequence
is unique up to permutation of flips of maximal $\mu$-weight). In
\cite{A}, Theorem 3.5, Agol shows that this sequence becomes periodic
at some stage, up to the action of $f$ and rescaling of the measure
$\mu$. Moreover, there is a common measured train track in any two
such sequences associated to two measured train tracks carrying
$L^s_{\phi_S}$ (\cite{A}, Corollary 3.4). Hence the periodic sequence
of maximal splittings does not depend on the choice of
$(\theta,\mu)$. It gives rise to a taut ideal triangulation of the
cusped manifold $M$, which is layered for $Z:=S\setminus
Sing(\phi_S)$, as described before Proposition
\ref{fibercase1}. Moreover, this triangulation is veering, and
conversely, any taut ideal triangulation which is veering and layered
for some fibration is associated to a periodic sequence of maximal
splittings (\cite{A}, Proposition 4.2). One has the same result if one
uses the unstable foliation $L^u_{\phi_S}$ of $\phi_S$, instead of
$L^s_{\phi_S}$; the two layered veering triangulations for
$L^u_{\phi_S}$ and $L^s_{\phi_S}$ coincide after a pre-branching total inversion. By starting with $\phi_S^{-1}$ one
obtains the same layered veering triangulations, and taking a
conjugate of $\phi_S$ yields a layered veering triangulation which
differs at most by a simplicial isomorphism. Hence, the layered veering
triangulation is uniquely determined by the fibration $f$. Let us
denote it $T_f$.

As in the proof of Theorem \ref{fibercase1} (3), consider the cones $C_f ={\rm Int}(H^+(T_f))$ in $H_2(M,\partial M;\mz)$ associated to the fibrations $f\in C(F_M)$. They form an open cover of $C(F_M)$. Consider one of these cones. Since it is rational, its boundary contains an integral
point. If this point lies in the interior of $C(F_M)$, then it
corresponds to some fibration $g$. As above, a multiple of it will be
fully carried by $\Sigma_f$ and $\Sigma_g$, and so $C_f \cap C_g\ne
\emptyset$. If $C_f \neq C_g$, an integral intersection point should
be fully carried by the two triangulations $T_f$ and $T_g$, which is impossible. Hence
the integral boundary points of $C_f$ lie in $\partial C(F_M)$, that
is, $C_f$ is the cone over $F_M$ for every $f\in C(F_M)$. \cvd

\subsection{Sutured manifold hierarchies} 
The layered triangulations of mapping tori can be recasted in the more
general framework of sutured $3$-manifolds \cite{L}: the initial
sutured manifold $(M_0,\gamma_0)$ is the mapping torus itself,
$M_0=Z\times [0,1]/\phi$, and the suture is the family of tori
$T(\gamma_0)$ formed by the boundary components of $M_0$ (no family of
annuli $A(\gamma_0)$); a fibre $Z_0$ is a ``styled surface'' for
$(M_0,\gamma_0)$, and we have a sutured manifold decomposition $
(M_0,\gamma_0) \rightarrow_{Z_0} (M_1,\gamma_1)$ along $Z_0$, where:
\begin{itemize}
\item $M_1 = M_0 \setminus U(Z_0)$ where $U(Z_0)$ is the interior of a
  regular neighbourhood of $(Z_0,\partial Z_0)$ in $(M_0,\partial
  M_0)$.
\item The $\gamma_1$-decomposition of $\partial M_1$ consists of the
  family $A(\gamma_1)$ of annuli which form the components of
  $\partial M_0 \setminus (U(Z_0)\cap \partial M_0)$, while $\Rr_\pm$
  consists of two parallel copies of $Z_0$ transversely oriented by
  the orientation of $[0,1]$.
\end{itemize}
The construction of a layered triangulation of $M_0$ can be decomposed 
in two steps:
\begin{itemize}
\item Define the {\it ideal region} $\delta$ of $(M_1,\gamma_1)$ as
  the union of the annuli $A(\gamma_1)$. Then use the construction of
  Section \ref{FIBER} to produce an ideal triangulation $\tilde T$ of
  the space obtained from the sutured manifold $(M_1,\gamma_1)$ by
  collapsing to one point each component of the ideal region
  $\delta$. The surfaces $\hat \Rr_\pm$ are eventually unions of ideal
  triangles.
\item By means of the same procedure, fill $U(Z_0)$ to get a layered
  triangulation of $M_0$.
\end{itemize}

The general construction of taut ideal triangulations in \cite{L} holds for any
irreducible $3$-manifold $M_0$ such that $\partial M_0$ is a non empty
collection of incompressible tori, and which admits furthermore a {\it
  sutured manifold hierarchy}
$$(M_0,\gamma_0,\delta_0)\rightarrow_{Z_0}  (M_1,\gamma_1,\delta_1) \dots 
\rightarrow_{Z_n} (M_{n+1},\gamma_{n+1},\delta_n).$$
Let us recall a few main features of sutured manifold hierarchies:
\begin{enumerate}
\item[(a)] The suture $\gamma_0$ consists of the family $T(\gamma_0)$
  of boundary components of $M_0$ (so there is no family of annuli
  $A(\gamma_0)$), and $\delta_0=\gamma_0$.
\item[(b)] Every $(M_j,\gamma_j)$ is a {\it taut sutured manifold}
  with non empty boundary (in the sense of Definition 2.2 of
  \cite{S}), such that no component of $\partial M_j$ is left
  untouched by the union of the family of tori $T(\gamma_j)$ and the
  family of annuli $A(\gamma_j)$.
\item[(c)] Every $\delta_j$ is an {\it ideal region} of
  $(M_j,\gamma_j)$ consisting of all of $T(\gamma_j)$, some components
  of $A(\gamma_j)$, and some ``squares'', i.e.  regions between two
  transverse arcs in some component of $A(\gamma_j)$, and it verifies
  the following conditions:
\begin{itemize}
\item no component of $A(\gamma_j)$ is left untouched by
$\delta_j$; 
\item The spaces obtained from the surfaces $\Rr_{\pm,j}$ by
  collapsing to one point each component of $\delta_j$ admit
  triangulations whose vertices are among the collapsed components.
\item $M_j$ contains no $\delta_j$-essential annulus (see page 12 of \cite{L}).
\end{itemize}
\item[(d)] Every $(M_j,\gamma_j) \rightarrow_{Z_j}
  (M_{j+1},\gamma_{j+1})$ is the sutured manifold decomposition along
  a {\it styled} proper surface $Z_j$ in $(M_j,\gamma_j)$, and $\delta_1=
  A(\gamma_1)\cup T(\gamma_1)$, $\delta_{j+1}=\delta_j \cap
  (A(\gamma_j)\cup T(\gamma_j))$.
\item[(e)] $H_2(M_{n+1},\partial M_{n+1};\Z)$ is trivial.
 \end{enumerate}  
 Moreover the surface $Z_0$ verifies the following properties:
 \begin{itemize}  
 \item[{(i)}] Every component of $Z_0$ has negative Euler
   characteristic and non empty boundary.
 \item[{(ii)}] $Z_0$ realizes the Thurston norm of its class in
   $H_2(M,\partial M;\Z)$.
 \item[{(iii)}] For every component $T$ of $\partial M$, $T\cap Z_0$
   is made of essential and coherently oriented parallel simple
   curves.
 \item[{(iv)}] There is no properly embedded essential annulus in
   $M_0$ disjoint from $Z_0$.
\end{itemize}     
\begin{teo}{\rm \cite{L}} If $M_0$ is irreducible, has non empty
  boundary formed by incompressible tori, and there is no properly
  embedded essential annuli in $M$, then every non trivial element of
  the (necessarily non trivial) image of $H_2(M_0,\partial M_0;\Z)$ in
  $H_1(\partial M_0;\Z)$ can be represented by a surface $Z_0$
  verifying the properties (i)-(iv). Moreover, every such a couple
  $(M_0,Z_0)$ can be included in a sutured manifold hierarchy.
\end{teo}
Note that the hypotheses on $M_0$ are verified if the interior of
$M_0$ carries a finite volume complete hyperbolic structure.  \medskip
  
Given a sutured manifold hierarchy as above, the construction of a
taut triangulation of $\hat M_0$ goes by induction backwards along the
hierarchy. The last step is as follows: by induction there is a taut
triangulation of $M_1$ with ideal region $\delta_1$. Then $R_-$ and
$R_+$, which are two parallel copies of $Z_0$, inherit ideal
triangulations. Fix a sequence of flips connecting the two
triangulations, and fill in a taut triangulated cylinder obtained via the
same construction as for the layered triangulations. We
eventually get a taut triangulation of $\hat M_0$.

Let us call it {\it a taut triangulation of $\hat M_0$ dominated by
  the given hierarchy}. The hierarchy, and in particular the sequence
of decomposing surfaces $Z_0,\dots, Z_n$, is an {\it intrinsic
  structure} obtained from $M_0$. The arbitrary choices that
eventually can produce different taut ideal triangulations of $M_0$
are of the same type already encountered with layered triangulations. So we
have:
\begin{prop}\label{prophier} For every sutured manifold hierarchy
   emanating from $M_0$, all the taut triangulations of $\hat M_0$
  dominated by the hierarchy determine the same element of $\tau(M)$.
\end{prop}

This achieves point (2) in Theorem  \ref{existNA}.

\subsection{On $\Z/2$-taut structures}\label{Z2taut}
There is simple procedure to construct manifolds whose boundary is not
necessarily made of tori, and carrying $\Z/2$-taut triangulations. Let
$(T,\omega)$ be a taut triangulation of $\hat M$. Let $e$ be an edge
of $T$ and $\bar e= e\cap M$. Then $(T,\omega)$ lifts to a $\Z/2$-taut
triangulation $(\tilde T,\tilde \omega)$ of $\hat {\tilde M}$, where
$\tilde M$ is a cyclic covering of $M$ branched along $e$ (if
any). For example, let us place ourselves in the situation of the
beginning of Section \ref{FIBER}. We say that a properly embedded arc
$\gamma$ in $Z$ is essential if it can be realised as the truncature
of an edge in an ideal triangulation of $\hat Z$. Consider a layered
triangulation $(T,\omega)$ of a mapping torus $M_\psi$ constructed by
using such a triangulation.  Then $\gamma$ is the truncature of an
edge of $T$ and we can perform the above construction along $\gamma$,
obtaining a $\Z/2$-taut triangulation $(\tilde T, \tilde \omega)$. One
can prove that this defines a $\Z/2$-taut structure
$\sG_{\psi,\gamma}$ on the so obtained manifold $\tilde M$.

\section{Relative non ambiguous and taut structures on 
pairs $(M,L)$}\label{ML} 
Assume that $M$ is closed and $L$ is a {\it non empty} link in $M$.
We are going to outline a theory of relative non ambiguous structures,
parallel to what we have already done in the ideal case.  In
particular we will introduce a notion of {\it relative taut
  structure}.
\smallskip

\begin{figure}[ht]
\begin{center}
\includegraphics[width=6cm]{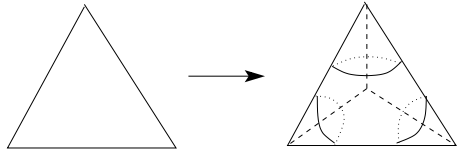}
\caption{\label{bubble} Naked bubble move}
\end{center}
\end{figure}

{\bf On the bubble move.} Dealing with possibly {\it non} ideal naked triangulations of $\hat M$
(for instance when $\partial M = \emptyset$), 
we must complete the naked transit equivalence with the 
$0\leftrightarrow 2$ {\it bubble move} which modifies the set of vertices 
(see Figure \ref{bubble}). 

A positive naked bubble move $T\rightarrow T'$ applies at a $2$-face
$F$ of $T$, and produces a $3$-ball $B$ triangulated by two tetrahedra
glued along three $2$-faces, so that the boundary of $B$ is
triangulated by two copies of $F$ glued along their boundary.  

\smallskip

We say that $(T,\omega)\rightarrow (T',\omega')$ is a {\it bubble
  pre-branching (pb-) transit} if for every $2$-face $F$ which is common to
$T$ and $T'$ the $\omega$ and $\omega'$ co-orientations of $F$
coincide, and if the restriction of $\omega'$ on the boundary of $B$
consists of two copies of the restriction of $\omega$ to $F$. For a
negative bubble move, the latter condition is replaced by: the
restriction of $\omega'$ on the boundary of $B$ consists of two copies
of a same co-orientation on $F$.  It is easy to see that every
pre-branching $\omega$ on $T$ gives rise to three positive bubble
pre-branching transits $(T,\omega)\rightarrow (T',\omega'_j)$,
$j=1,2,3$; in every case the two pre-branched tetrahedra of
$(T',\omega'_j)$ share their four square edges, two of them being
internal to $B$, two being boundary edges of the two copies of $F$ in
$T'$.  Hence the three prebranching transits are obtained one from
another by simultaneous reversal of the $\omega_*'$-co-orientations
of two oppositely co-oriented $2$-faces among the three inner
$2$-faces of the $3$-ball $B$. Although these are related to each
other by such evident symmetries, {\it in a strict sense no
  bubble $pb$-transit is ``non-ambiguous''}. We need some further
input in order to select one. We are going to do it in the framework
of (relative) distinguished triangulations $(T,H)$ of pairs
$(M,L)$. This kind of triangulations has been already used to
construct the QHI for $(M,L)$ \cite{Top}.  

\smallskip

By definition a {\it distinguished triangulation} $(T,H)$ of
$(M,L)$ is formed by a triangulation $T$ of the closed manifold $M$
and a {\it Hamiltonian} subcomplex $H$ of $T^{(1)}$ isotopic to
$L$. Then we have to consider naked distinguished triangulation moves $(T,H)\leftrightarrow (T',H')$. 

\smallskip

For every positive $2\rightarrow 3$ or lune move $T\rightarrow T'$,
every pair $(T,H)$ gives rise to a distinguished move $(T,H)\rightarrow (T',H')$ where $H=H'$. The inverse moves
$(T,H)\leftarrow (T',H')$ are by definition the negative distinguished $3\rightarrow 2$ or lune moves. 

\smallskip

In order to enhance a positive bubble move $T\rightarrow T'$ to some distinguished bubble move $(T,H)\rightarrow (T',H')$, we require that at least one edge $e$ of
the $2$-face of $T$ involved in the move belongs to $H$; $e$
determines one internal $2$-face of the two new tetrahedra of $T'$,
having three boundary edges $e$, $e'$, and $e''$.  Then we set $H'=
(H\setminus e)\cup e'\cup e''$. Recall (see \cite{Top}, and also \cite{GA}
in a more general setting) that the quotient set of distinguished
triangulations of $(M,L)$ up to distinguished moves
encodes the topological type of $(M,L)$.

\begin{defi}\label{NAbubble}{\rm Given a positive naked distinguished
    bubble move $(T,H)\rightarrow (T',H')$ and a pre-branching
    $(T,\omega)$, the {\it non ambiguous} pre-branching ($pb$-) enhancement $(T,H,\omega)\to
    (T',H',\omega')$ is the one among the three possible pre-branchings $w'$ such
    that there are no abstract diagonal edges of $(B,\omega')$ at both edges
    $e'$ and $e''$, while the two abstract edges of $(B,\omega')$ at $e$ are
    diagonal.  The negative non ambiguous distinguished bubble $pb$-transit is defined
    consequently. A (non ambiguous) distinguished ideal $pb$-transit is the relative enhancement of an ordinary (non ambiguous) ideal $pb$-transit.}
\end{defi}

\begin{defi}\label{NAML}{\rm The {\it non ambiguous distinguished $pb$-transit
      equivalence} on the set of distinguished triangulations
    of $(M,L)$ is generated by isotopy and non ambiguous distinguished $pb$-transits. We denote by $\Nn\Aa(M,L)$ the quotient set. Every
    coset is a {\it relative non ambiguous structure} on $(M,L)$.}
\end{defi}

\begin{defi}\label{tautML}{\rm A distinguished triangulation
$(T,H,\omega)$ of $(M,L)$ is {\it relative taut} if $(T,\omega)$ is a pre-branched
triangulation, around every edge $e\in T^{(1)} \setminus H$ there are
exactly two diagonal abstract edges, and around every edge $e\in H$ there are
no diagonal edges. We say that $(T,H,\omega)$ is {\it relative} $\Z/2$-{\it taut} if $(T,\omega)$
is locally taut and around every edge  $e\in T^{(1)} \setminus H$ there are
diagonal edges.}
\end{defi}

 Similarly to point (3) of Proposition \ref{N2_via_pb} we easily
have:
\begin{prop} If 
  $(T,H,\omega)$ is a relative $\Z/2$-taut triangulation of $(M,L)$ 
and there is no diagonal edge around any edge of $H$, then
it is a relative taut triangulation.
\end{prop}
Let us come to the relative version of ``charges''. Every triangulation
$T$ of $M$ can be considered as an ideal triangulation with respect to its set of vertices.
So the notions of $\Z$- and $\Z/2$-charges of Section \ref{seccharges}
could be adopted as well. However we have to modify the definitions
because of the presence of the subcomplex $H$.

\smallskip

A $\Z$-{\it charge on $(T,H)$} is defined as in Section
\ref{seccharges}, with the difference that for every edge $e\subset H$
we require that the total charge $C(e)=0$ instead of $C(e)=2$.
Every $\Z/2$-charge on $T$ is also a $\Z/2$-charge on $(T,H)$.

 Similarly to point (4) of Proposition \ref{N2_via_pb} we
have:
\begin{prop} Every distinguished triangulation $(T,H)$ of $(M,L)$
carries $\Z$-charges.
 \end{prop}

 The proof is not easy (see \cite{Top}, and the references given for (4) of Proposition \ref{N2_via_pb}).

\smallskip
 
We say that a locally taut $\mz/2$-charge $\gamma$ on $T$ is
$\Z/2$-{\it taut} for the pair $(T,H)$ if around every edge of
$T^{(1)}\setminus H$, there are at least two $1$-colored abstract
edges. Hence, $(T,H,\omega)$ is a relative $\Z/2$-taut triangulation of $(M,L)$
if and only if the $\Z/2$-charge $\gamma_\omega$ is taut for
$(T,H)$. If $(T,H,\omega)$ is a relative taut triangulation, then
$\gamma_\omega$ is a $\Z$-charge.

We can define the quotient sets $c(M,L,\Z/2)$
and $c(M,L,\Z)$ of $\mz/2$- or $\mz$-charged distinguished triangulations of $(M,L)$ up to charge transits, in the usual way. Similarly to Proposition \ref{Z2cmeaning} we have:

 \begin{prop} The sets $c(M,L,\Z/2)$ and $c(M,L,\Z)$ encode
   $H^1(M;\Z/2)$: in particular, for every class $\alpha \in H^1(M;\Z/2)$
   and every distinguished triangulation $(T,H)$ of $(M,L)$ there is a
   $\Z$-charge on $(T,H)$ that realizes $\alpha$.
\end{prop}  

Similarly to Lemmas \ref{taut-preserved} and \ref{on_transits} we have:
\begin{lem} 
Let $(T,H,\omega)\rightarrow (T',H',\omega')$ be any distinguished $pb$-transit (bubble included). Then the following
facts are equivalent:
\smallskip

(1) The transit is non ambiguous.

\smallskip

(2)The transit sends relative $\Z/2$-taut triangulations to
 relative $\Z/2$-taut triangulations.

\smallskip

(3) The transit sends relative taut triangulations to
relative taut triangulations.

\end{lem}

\begin{lem} Every (necessarily non ambiguous)
  $pb$-transit $(T,H,\omega)\rightarrow (T',H',\omega')$ (bubble
  included) of relative taut (resp. $\Z/2$-taut) triangulations
  $(T,H,\omega)\to (T',H',\omega')$ induces a transits of $\Z$- (resp. $\Z/2$-taut) charges on $(T,H)$.
\end{lem}

 Hence, similarly to Proposition \ref{taut-and-homology}, we
 can define the set $\tau(M,L)\subset \Nn\Aa(M,L)$ of {\it relative
   taut structures} on $(M,L)$, and there is a well defined map 
$$\hG: \Nn\Aa(M,L)\to H^1(M;\Z/2). $$

Let $(T,H,\omega)$ be a distinguished pre-branched triangulation of
$(M,L)$. Formally consider $T$ as an ideal triangulation of $M$
obtained by removing a small $3$-ball around every vertex of $T$.
Then $\partial M $ is a union of $2$-spheres $S_j$, with induced
branched triangulations $(K_j,b_j)$.  For every $\Z$-charge $c$ on
$(T,H)$ and every sphere $S_j$, there are exactly two vertices
$v^j_{\pm}$ of $K_j$ with total charge $C(v^j_{\pm})=0$, while the
other vertices have total charge equal to $2$.  If the triangulation
is relative taut, then the combing associated to $(K_j,b_j)$ on every sphere
$S_j$ can be represented by a field having one source singular point at one of
these special vertices, say $v^j_+$, one pit at $v^j_-$, and is non
singular elsewhere.  

\medskip

In the rest of the section we show how to construct 
relative taut triangulations. Let us consider the following
situation: 

\begin{itemize}

\item $Y$ is a $3$-manifold such that $\partial Y$ is a collection of
  tori.  For simplicity we assume that $\partial Y$ is connected (the
  general case would be treated similarly).
\item We denote by $m$ a simple essential curve on $\partial Y$,
  considered up to isotopy.  Then $M$ is the manifold obtained by Dehn
  filling of $Y$, by attaching a solid torus $\ \Uu= D\times S^1$
  along $\partial Y$ so that $m$ is identified with a meridian
  $\partial D \times \{s_0\}$ of $\Uu$; $L$ is the knot in $M$ formed
  by the core of the solid torus. So we deal with this couple $(M,L)$.
\item We assume that $\hat Y$ is endowed with a {\it taut}
  triangulation $(T_Y,\omega_Y)$.
\end{itemize}
Our aim is to construct a distinguished triangulation $(T,H,\omega)$
of $(M,L)$ ``strictly related'' to $(T_Y,\omega_Y)$ and (relatively)
taut.  

\medskip

Let us consider the standard spine $P_Y$ of $Y$ made by the
$2$-skeleton of the cell decomposition of $\hat Y$ dual to $T_Y$.  Fix
a ``normal retraction" $r: Y \to P_Y$ (the inverse image of a regular
point of $P_Y$ is an interval, for a point in the interior of an edge
of $P_Y$ it is a cone over $3$ points, and for a vertex of $P_Y$ it is
a cone over $4$ points). Let $D:= D\times \{s_0\}$ be a properly
embedded meridian disk of $\Uu$ as above with boundary $\partial D=m
\subset \partial Y$.

\begin{defi}\label{well-placed} {\rm We say that the curve $m$ is {\it
      well placed (with respect to $r$)} if it verifies the following
    properties:
\begin{enumerate}
\item The restriction of $r$ to $m$ is an embedding in $P_Y$.
\item  The curve $r(m)$ does not pass through any vertex of
$P_Y$ and is transversal to the edges of ${\rm Sing}(P_Y)$.
\item The curve $r(m)$ is subdivided in arcs, such that each arc
has its endpoints on ${\rm Sing}(P_Y)$ and its interior is contained in one
component of  $P_Y \setminus {\rm Sing}(P_Y)$.
\item Every component of $P_Y \setminus ({\rm Sing}(P_Y) \cup  r(m))$ is an 
open $2$-disk.
\item The mapping cylinder $\Cc$ of $r_{|m}$ is embedded in
$Y$ and intersects $P_Y$ only along $r(m)$.
\end{enumerate}
}
\end{defi}

\begin{lem}\label{r(m)} Up to isotopy we can assume that
  the curve $m$ is well placed.
\end{lem}
\Dim A simple argument of general position does not exclude the
presence of possible simple self-crossings of $r(m)$ in some regions
of $P_Y \setminus {\rm Sing}(P_Y)$.  However such crossing points can
be eliminated by sliding $r(m)$ on ${\rm Sing}(P_Y)$ and introducing more
points in ${\rm Sing}(P_Y) \cap r(m)$. So eventually (up to isotopy)
$r(m)$ is embedded.\cvd

\medskip

The following Lemma is clear.

\begin{lem}\label{spineM} Let  the essential curve $m\subset \partial Y$ 
  be well placed with respect to the normal retraction $r:Y\to P_Y$.
  Let $\Dd= \Cc \cup D$.  Then $P_M := P_Y\cup \Dd$ is the standard
  spine dual to a triangulation $T$ of $M$ with one vertex. The knot
  $L$ is realized by the Hamiltonian subcomplex $H$ of $T^{(1)}$
  formed by one edge $e$ (with identified endpoints) dual to $\Dd$.
\end{lem}

So we have constructed a distinguished triangulation $(T,H)$ of
$(M,L)$.  By duality, the other abstract edges of $T$ correspond to
the components of $P_Y \setminus ({\rm Sing}(P_Y) \cup r(m))$, the
tetrahedra of $T$ correspond to the vertices of $P_M$. Every vertex of
$P_Y$ persists in $P_M$, while there is a ``new'' vertex each time the
curve $r(m)$ crosses an edge of ${\rm Sing}(P_Y)$. Hence we have a
natural inclusion of the set of naked abstract tetrahedra of $T_Y$
into the set of naked abstract tetrahedra of $T$ such that $H$ is
separated from the abstract edges of $T_Y$.

Next we want to define a pre-branching $\omega$ on $(T,H)$ in such a
way that it agrees with $\omega_Y$ on the abstract tetrahedra of
$T_Y$.  Fix a weak branching $\tilde b$ on $T_Y$ which induces
$\omega_Y$. We can assume that the germs of $r(m)\cup {\rm Sing}(P_Y)$
at every point of $r(m)\cap {\rm Sing}(P_Y)$ are contained in a
disjoint union of branched ``butterfly'' neighbourhoods of some
vertices of $P_Y$.  Fix an orientation of the region $\Dd$ of $P_M$,
hence of its boundary $r(m)$.  It results that the portion of $P_M$
formed by the union of $\Dd$, those butterflies, and a regular
neighbourhood of $r(m)$ in $P_Y$ is a {\it branched} surface with
boundary, having one vertex at each point of $r(m)\cap {\rm
  Sing}(P_Y)$. It induces a pre-branching on every ``new'' tetrahedron
of $T$; by keeping the pre-branching $\omega_Y$ on every ``old''
tetrahedron of $T_Y$ we have eventually defined a pre-branching
$(T,H,\omega)$ with the desired properties.  The following Proposition
summarizes the properties of a distinguished pre-branched
triangulation of $(M,L)$ such that $T$ is dual to the spine $ P_M :=
P_Y\cup \Dd$ and $\omega$ is obtained by implementing the above
procedure.

\begin{prop}\label{pre-taut}  (1) 
  The set of abstract pre-branched tetrahedra of $(T_Y,\omega_Y)$ is
  included into the set of abstract pre-branched tetrahedra of
  $(T,H,\omega)$.

  (2) The abstract edges of $(T,H,\omega)$ around $H$ do not belong to
  the abstract edges of $(T_Y,\omega_Y)$.

(3) There is no $\omega$-diagonal edge at $H$.
\end{prop}
\Dim Points (1) and (2) follow directly from the construction. As for
(3), forgetting the orientation of $r(m)$ and the
$\omega_Y$-orientation of ${\rm Sing}(P_Y)$, at a branched butterfly
$\Bb$ around a point $x$ of $r(m)\cap {\rm Sing}(P_Y)$ there are two
possible configurations:
\begin{itemize}
\item The arc of $r(m)$ is {\it smoothly} embedded in $\Bb$ (with respect
to the branched $C^1$-structure of $\Bb$).
\item  The arc of $r(m)$ has a ``{\it cusp}'' at $x$.
\end{itemize}
By taking into account the orientations, there are $4$ possible
configurations.  It is easy to realize that in any case the abstract
edge dual to $\Dd$ is {\it not} a diagonal edge, and so we get the
desired behaviour at $H$.\cvd 

\medskip

Clearly point (3) above is among the tautness conditions. One might
wonder if any $(T,H,\omega)$ constructed in this way is a taut
triangulation. Unfortunately we have to strengthen our assumption. Let
us analyse the possible position of a well placed curve $m$ with
respect to the oriented train track $\theta$ on $\partial Y$
associated to $(T_Y,\omega_Y)$. As this last is taut, every region,
say $\Rr$, of $\partial Y \setminus \theta$ is a {\it bigon} with two
cusp points corresponding to the two diagonal corners at the vertex of
$\partial T_Y$ dual to $\Rr$.  The boundary of $\Rr$ is formed by two
oriented smooth arcs joining the cusp points. The retraction $r$ maps
the region $\Rr$ onto a region $r(\Rr)$ of $P_Y \setminus {\rm
  Sing}(P_Y)$; $r^{-1}(r(\Rr))\cap \partial Y$ consists of two regions
$\Rr$ and $\Rr'$, dual to the endpoints of a same truncated edge of
$T_Y$; every arc of $r(m)$ which intersects $r(\Rr)$ lifts to an arc
of $m$ in $\Rr$ or $\Rr'$ joining two (generic) points of $\partial
\Rr$ or $\partial \Rr'$. So we can state:
\begin{defi}\label{verywellplaced}{\rm An essential curve $m
    \in \partial Y$ is {\it very well placed} (with respect to the
    retraction $r$) if it is well placed and moreover, for every
    region $\Rr$ of $\partial Y \setminus \theta$ and every arc
    $\gamma$ of $m$ traversing $\Rr$, one of the following situations
    is realized:}
  \begin{enumerate}
  \item {\rm The arc $\gamma$ separates the two cusp points of
      $\partial \Rr$.}
  \item {\rm If the arc $\gamma$ does not separate the two cusp points
      of $\partial \Rr$, then together with a sub-arc $\sigma$ of
      $\partial \Rr$, $\gamma \cup \sigma$ is the boundary of a bigon
      contained in $\Rr$. Recall that both $\sigma$ and $\gamma$ are
      oriented. So we require that one point of $\gamma \cap \sigma$
      is a source while the other is a pit.}
 \end{enumerate}
 \end{defi}

Finally we have:

\begin{prop}\label{taut-taut} 
  Let the essential curve $m\subset \partial Y$ be very well placed
  with respect to the normal retraction $r: Y\to P_Y$. Let
  $(T,H,\omega)$ be a distinguished pre-branched triangulation of
  $(M,L)$ such that $T$ is dual to the spine $P_M:=P_Y\cup \Dd$ and
  $\omega$ is obtained as described before Proposition
  \ref{pre-taut}. Then $(T,H,\omega)$ is a relative taut triangulation
  of $(M,L)$.
\end{prop}
\Dim We already know that every edge of $H$ satisfies the tautness
condition. It remains to check it at every edge of $T$ dual to a
region of $P_Y \setminus ({\rm Sing}(P_Y) \cup r(m))$. We can do it
inductively by implementing the following procedure. The initial
step: select an arc $\gamma_0$ of $r(m)$ which intersects a region
$\Rr_0$ of $P_Y \setminus {\rm Sing}(P_Y)$, with endpoints on ${\rm
  Sing}(P_Y)$, and such that $\gamma_0$ is ``innermost'' among the
arcs of $r(m)$ intersecting the region $\Rr_0$; that is, there is a
subregion, say $\Rr'_0$, bounded by the union of $\gamma_0$ and an arc
$\sigma_0$ contained in ${\rm Sing}(P_Y)$, and there is no arc of
$r(m)$ in the interior of $\Rr'_0$.

Since $\Rr'_0$ cannot be cut by $r(m)$, it is dual to an edge $e'$ of
$T$. By using the fact that $(T_Y,\omega_Y)$ is taut, looking at the
two possible configurations of the inverse image $\gamma$ of
$\gamma_0$ in some region $\Rr$ of $\partial Y$, and the fact that $m$
is very well placed, it is not hard to verify that $e'$ verifies the
tautness condition, as well as the persistent edges of $T_Y$. We iterate
the procedure by adapting the above ``innermost'' criterium to the
partial subdivision $P_Y \setminus ({\rm Sing}(P_Y)\cup \gamma_0)$,
and so on. One by one we create the edges dual to the regions of $P_Y
\setminus ({\rm Sing}(P_Y)\cup r(m))$, verifying at each step the
tautness condition.\cvd

\medskip

The hypothesis of the last Proposition is not too demanding. For example
we have:
 
\begin{prop}\label{transverseOK} 
  Let the essential curve $m\subset \partial Y$ be positively
  transverse to the oriented train track $\theta$ of $(\partial
  T_Y,\partial \omega_Y)$. Then $m$ is very well placed.  Moreover, by
  varying $m$ by an isotopy through curves positively transverse to
  $\theta$, the associated taut triangulations of $(M,L)$ determine
  the same relative taut structure.
\end{prop}
\Dim As for the first claim, it is clear that only arcs of $m$
which separate two cusp points do occur.  As for the second, such
isotopies are generated by the local ones where $m$ crosses a
vertex of $\theta$.  One realizes easily that the corresponding taut
triangulations are related by a non ambiguous transit. \cvd
\medskip

{\bf Remarks and questions.}  (1) Starting with a taut triangulation
$(T,H,\omega)$ as above, via non ambiguous transits (in particular the
bubble ones) we can reach taut triangulations of $(M,L)$ verifying
further conditions, like having all edges with distinct endpoints.

(2) A main question is: {\it which essential curves $m$ on $\partial Y$
  can be made very well placed, possibly by modifying $(T_Y,\omega_Y)$
  via non ambiguous transits ?}

(3) Another question concerns to what extent the taut class of
$(T,H,\omega)$ can vary (even assuming that the class of
$(T_Y,\omega_Y)$ is fixed). Note that there are several arbitrary
choices along the construction of $(T,H,\omega)$, the most relevant one
being the specific way $m$ is made well placed and possibly very well
placed.

\section{The symmetry defect}\label{SYMD} 
We need to recall a few features of the QH triangulations of $M$ or
$(M,L)$ (see eg. \cite{AGT}).  We assume at first that $M$ is a cusped
manifold, and postpone the case of pairs $(M,L)$ to Section \ref{MLdefect}.

\medskip

A QH tetrahedron $(\Delta,b,w,f,c)$ consists of an oriented branched
tetrahedron $(\Delta,b)$ (hence with ordered vertices
$v_0,v_1,v_2,v_3$) together with an ordered decoration
$d=(d_0,d_1,d_2)$, where $d_j=(w_j,f_j,c_j)$, $j\in \{0,1,2\}$, is
associated to the triple of couples of opposite edges; in fact we
have: 
\begin{enumerate}
\item $d_0$ labels the edge $[v_0,v_1]$ (and the opposite edge), $d_1$
  labels $[v_1,v_2]$, and $d_2$ labels $[v_0,v_2]$.
\item The {\it shape parameters} $w_j\in \C \setminus \{0,1\}$ verify 
$w_{j+1}=(1-w_j)^{-1}$ (indices $j$ mod$(3)$);
\item The {\it flattenings} $f_j \in \Z$ determine {\it log branches}
  $ l_j:= \log(w_j)+f_ji\pi$ such that $l_0+l_1+l_2=0$.
\item The charges $c_j \in \Z$ satisfy $c_0+c_1+c_2=1$.
\end{enumerate}
By taking into account the $3$-simplex sign $*_b$, we say that $w_j^{*_b}$ and
$*_bl_j$ are {\it signed} shape parameters and log branches respectively. 
\smallskip

A QH triangulation $(T,\tilde b,w,f,c)$ is a weakly branched ideal
triangulation $(T,\tilde b)$ made of QH tetrahedra satisfying the
following global conditions at every edge $e$
of $T$:
\begin{itemize}
\item the {\it total signed shape parameter} at $e$, ie. the product of the 
signed shape parameters around $e$, is $W(e)=1$;
\item the {\it total signed log branch} at $e$, ie. the sum of
the signed log branches around $e$, is $L(e)=0$;
\item $c$ is a $\Z$-charge on $T$, as in Section \ref{seccharges}.
\end{itemize}
The conditions $W(e)=1$ mean that the system of signed shape parameters is a point of the gluing variety (ie. the algebraic set defined by the Thurston edge equations) supported by the (naked) ideal triangulation $T$. 
\smallskip

For every QH tetrahedron $(\Delta,b,w,f,c)$ define the triple $\bw:=(\bw_0,\bw_1,\bw_2)$ of $N$th-roots of the
shape parameters by
\begin{equation}\label{qshapes}
\bw_k := \exp\left(\frac{\log(w_k) +\pi i (N+1)(f_k-*_bc_k)}{N}\right), 
\ k =0,1,2.
\end{equation}
Similarly as above, $\bw_k^{*_b}$ is called a {\it signed ``quantum'' shape parameter}. Call {\it total signed quantum shape parameter} $\WW(e)$ of an edge $e$ of $T$ the product of the 
signed quantum shape parameters around $e$. Set $\zeta_N := \exp (2i\pi/N)$. We have the relations:
\begin{itemize}
\item $\bw_0\bw_1\bw_2=-\zeta_N^{*_b\frac{N-1}{2}}$;    
\item $\WW(e)= \zeta_N^{-1}$ around every edge $e$ of $T$.  
\end{itemize}

The local symmetrization factor of $(\Delta,b,w,f,c)$ is defined by
$$\alpha_N(\Delta,b,w,f,c):= \bigl(\bw_0 ^{-c_1}\bw_1^{c_0}\bigr)^{\frac{N-1}{2}}$$
and the global symmetrization factor of $\Tt:=(T,\tilde b,w,f,c)$ is 
$$\alpha_N(\Tt) := \prod_{\Delta\in \Tt^{(3)}}\alpha_N(\Delta,b,w,f,c).$$

\subsection{Pre-branching dependence} We consider at 
first the behaviour of $\alpha_N(\Tt)$ when only the weak branching varies.

\begin{lem}\label{orC} 
  Let $\Tt$ and $\Tt'$ be QH triangulations which differ only by their
  weak branchings $\tilde b$ and $\tilde b'$.  If $\tilde b$ and
  $\tilde b'$ induce the same pre-branching, then
  $\alpha_N(\Tt)=\alpha_N(\Tt')$.
\end{lem}
\Dim This statement is of local nature, as we can pass from $\Tt$ to
$\Tt'$ by changing the local branchings at some tetrahedra without
changing the induced pre-branching. Given a local pre-branching on an
oriented tetrahedron $\Delta$ there are four local branchings that
induce it. The set of such local branchings is obtained from a single
one, say $b$, by reordering the vertices $v_0,v_1,v_2,v_3$ of $\Delta$
using the four permutations of the set $J_4:=\{0,1,2,3\}$ that form
the cyclic subgroup $<\sigma>$ of the symmetric group $S(J_4)$,
generated by the cycle $\sigma:=(0,1,2,3)$. Let us assume for instance
that $*_b=1$ and compare $\alpha_N(\Delta,b,w,f,c)$ and
$\alpha_N(\Delta,b_\sigma,w_\sigma,f_\sigma,c_\sigma)$, where
$b_\sigma$ is the branching obtained by applying $\sigma$ on $b$,
while $w_\sigma,f_\sigma,c_\sigma$ denote the result of reordering
$w,f,c$ with respect to $b_\sigma$. Note that $*_{b_\sigma}=-1$. So we
have:
$$w_\sigma =(w_1^{-1},w_0^{-1},w_2^{-1}),\  f_\sigma=(-f_1,-f_0,-f_2), \
c_\sigma=(c_1,c_0,c_2). $$
Hence $\alpha_N(\Delta,b_\sigma,w_\sigma,f_\sigma,c_\sigma)$ is equal to
\begin{multline*}
\bigl(\exp(\frac{\log(w_1^{-1})- \pi i(N+1)(f_1-c_1)}{N}\bigr)^{-c_0}\exp
\bigl(\frac{\log(w_0^{-1})- \pi i(N+1)(f_0-c_0)}{N})^{c_1}
\bigr)^{\frac{N-1}{2}} \\ = 
\bigl(\exp(\frac{\log(w_1)+ \pi i(N+1)(f_1-c_1)}{N})^{c_0}
\exp (\frac{\log(w_0)+ \pi i(N+1)(f_0-c_0)}{N})^{-c_1}\bigr)^{\frac{N-1}{2}}
\end{multline*}
which in turn is equal to $\alpha_N(\Delta,b,w,f,c)$. One proceeds 
similarly if at the beginning we assume that $*_b=-1$. \cvd

\medskip

So, given a QH triangulation $\Tt=(T,\tilde b,w,f,c)$, we see that
{\it $\alpha_N(\Tt)$ depends on $\tilde b$ only through the underlying
  pre-branching $\omega$.}  Let us analyze now the effect of performing
the total inversion of the pre-branching.
For every branched tetrahedron $(\Delta,b)$
of $(T,\tilde b)$ with $b$-ordered vertices $v_0,v_1,v_2,v_3$, the
total inversion changes $b$ by reordering the vertices by the
permutation $\tau=(0,1)(2,3)$. Let us compare
$\alpha_N(\Delta,b,w,f,c)$ and $\alpha_N(\Delta, b_\tau ,w_\tau ,
f_\tau , c_\tau )$. Assume for example that $*_b = 1$; then also
$*_{b_\tau} = 1$. We see that
$$ w_\tau = w, \ f_\tau = f,\  c_\tau = c$$
so we have immediately:
\begin{lem}\label{cor-orC}  
  Let $\Tt$ and $\Tt'$ be QH triangulations which differ only 
by their weak branchings $\tilde b$ and $\tilde
  b'$. If the underlying pre-branchings differ by the total inversion,
  then $\alpha_N (\Tt )= \alpha_N (\Tt' )$.
\end{lem}

\subsection{Transit invariance} 
Now we consider the behaviour of $\alpha_N(\Tt)$ under transits of QH
triangulations $(T,\tilde b,w,f,c)\rightarrow (T',\tilde
b',w',f',c')$. Such a transit is supported by a $wb$-transit
$(T,\tilde b)\to (T', \tilde b')$ and involves natural transition
rules of the whole decoration $(w,f,c)$ (see \cite{GT}, Section
2.1.3).  For the charges,
these rules have been reminded in Section \ref{seccharges}; the
transition rules for the log branches associated to $(w,f)$ are
formally the same, and act as a ``logarithm'' of the transition rules
of the shape parameters.  As we can freely change the weak branching
without modifying the induced pre-branching, we can apply Lemma
\ref{wb} and assume that the $wb$-transit is locally branched.

\begin{prop}\label{NA-invariance} 
  Let $(T,\tilde b,w,f,c)\rightarrow (T',\tilde b',w',f',c')$ be a
  $2\leftrightarrow 3$ or lune QH transit.  Then we have
  $\alpha_N(T,\tilde b,w,f,c) = \alpha_N(T',\tilde b',w',f',c')$ if
  and only if $(T,\tilde b)\rightarrow (T',\tilde b')$ is non
  ambiguous. 
\end{prop} 

\Dim By Lemmas \ref{orC} and \ref{wb}, it is enough to check the claim
for an arbitrarily chosen locally branched transit
$(T,b)\leftrightarrow (T',b')$ which covers the given pre-branching
transit. Assume at first that the latter is $2\leftrightarrow 3$ non
ambiguous. Then we can restrict to the remarkable $b$-transits
pointed out before Lemma \ref{wb}. The invariance of $\alpha_N(\Tt)$
in this case was already an important fact used in Theorem 5.7 of
\cite{GT}, where we sketched the proof. For the sake of completeness
we give here all the details in the case of Figure \ref{NAb} (the
other case is similar). Let us write the local symmetrization factor
of a QH tetrahedron $(\Delta,b,w,f,c)$ with $*_b=1$ in the form
$$\alpha_N(\Delta,b,w,f,c):=\bigl(\bw_0 ^{-c_1}\bw_1^{c_0}\bigr)^{\frac{N-1}{2}} =  
\exp\bigl(\frac{N-1}{2}(-c_1\frac{l_0}{N} +
c_0\frac{l_1}{N})\bigr)\exp\bigl(i\pi\frac{N-1}{2} (-c_1f_0
+c_0f_1)\bigr).$$ Denote by $l^j_k$ the $k$-th log branch of the
tetrahedron opposite to the $j$-th vertex according to the vertex
ordering defined by the branching of Figure \ref{NAb}, and similarly $c^j_k$
and $f^j_k$ for the charges and the flattenings. The products of the
local symmetrization factors corresponding to the tetrahedra involved
in the QH transit are
$$(-1)^{\frac{N-1}{2}(f_0^1c_1^1+f_1^1c_0^1+f_0^3c_1^3+f_1^3c_0^3)}\ 
\exp\biggl(\frac{N-1}{2N}\bigl(-c_1^1l_0^1+ c_0^1l_1^1-c_1^3
l_0^3 +c_0^3l_1^3\bigr)\biggr)$$
before the transit, and 
$$(-1)^{\frac{N-1}{2}(f_0^0c_1^0+f_1^0c_0^0+f_0^2c_1^2+
  f_1^2c_0^2+f_0^4c_1^4+f_1^4c_0^4)}\
\exp\biggl(\frac{N-1}{2N}
\bigl(-c_1^0l_0^0+c_0^0l_1^0-c_1^2l_0^2+c_0^2l_1^2-c_1^4l_0^4+c_0^4l_1^4\bigr) 
\biggr)$$ 
after the transit. Let us consider the exponential terms and prove that 
\begin{equation}\label{eq0}
-c_1^1l_0^1+ c_0^1l_1^1-c_1^3
l_0^3 +c_0^3l_1^3 = -c_1^0l_0^0+c_0^0l_1^0-c_1^2l_0^2+c_0^2l_1^2-c_1^4l_0^4
+c_0^4l_1^4.
\end{equation}
The QH transit implies the relations
$$l^3_0 = l^2_0+l^4_0\ ,\ c^3_0 = c^2_0+c^4_0\ ,\ l^3_1 = l^0_0+l^4_1\ ,\ c^3_1 
= c^0_0+c^4_1\ ,\ l^4_0 = l^1_1-l^0_1\ ,\ c^4_0 = c^1_1-c^0_1$$
$$l^2_0 = l^1_0-l^0_0\ ,\ c^2_0 = c^1_0-c^0_0\ ,\ l^2_1 = -l^0_2-l^4_2 
= l^0_0+l^0_1+l^4_0+l^4_1\ ,\ c^2_1 = 2-c^0_2-c^4_2 =
c^0_0+c^0_1+c^4_0+c^4_1.$$ Taking the difference of both sides of
\eqref{eq0} and substituing these identities we get
\begin{multline*} -c_1^1l_0^1+c_0^1l_1^1+c_1^0l_0^0-c_0^0l_1^0 -
(c^0_0+c^4_1)(l^1_0-l^0_0+l^1_1-l^0_1) + (l^0_0+l^4_1)(c^1_0-c^0_0+c^1_1-c^0_1) 
+ \\ + (c^0_0+c^0_1+c^4_0+c^4_1)(l^1_0-l^0_0) - 
(c^1_0-c^0_0)(l^0_0+l^0_1+l^4_0+l^4_1) + c^4_1(l^1_1-l^0_1)-l^4_1(c^1_1-c^0_1).
\end{multline*}
Cancelling terms this gives
\begin{multline*} -c_1^1l_0^1+c_0^1l_1^1-c_0^0l_1^1+c_1^1l_0^0+
(c^0_1+c^4_0)(l^1_0-l^0_0)-(c^1_0-c^0_0)(l^0_1+l^4_0) \\ =  
c^1_1(l^2_0-l^1_0)+l^1_1(c^1_0-c^2_0)-c_0^0l_1^1+c_1^1l_0^0 = 0.
\end{multline*}
The same argument proves that the signs are equal, by working with
coefficients mod$(2)$ and replacing the log branches with the
flattenings.

Consider now a $2\leftrightarrow 3$ QH transit with underlying ambiguous 
pre-branching transit. Again it is enough to consider any compatible 
branching transit. An instance is obtained by applying the transposition 
$(12)$ on the vertex ordering induced by the branching of Figure \ref{NAb}. 
This transposition alters only the local symmetrization factors of the 
QH tetrahedra opposite to the $3$-, $0$- and $4$-vertex, that it 
multiplies by 
$$(\bw_1^3)^{\frac{1-N}{2}}\ ,\ (\bw_0^0)^{\frac{N-1}{2}}  \ ,
\ (\bw_1^4)^{\frac{1-N}{2}}$$
respectively. Since $\bw_1^3 = \bw_0^0\bw_1^4$, we have 
$(\bw_1^3)^{\frac{1-N}{2}}\ne (\bw_0^0)^{\frac{N-1}{2}}(\bw_1^4)^{\frac{1-N}{2}}$, 
and so $\alpha_N(\Tt)$ is not invariant in this case. 

The proof for the lune transits is based 
on similar computations.\cvd

\subsection{Boundary QH triangulations }\label{boundaryQH} 
Let $\Tt=(T,\tilde b,w,f,c)$ be a QH triangulation of $M$. Our next
task is to study the behaviour of $\alpha_N(\Tt)$ when $c$ varies. Let
us recall from \cite{AGT} that: 
\begin{itemize}
\item the system of shape parameters $w$ encodes a
$PSL(2,\C)$-character $\rho:=\rho(w)$ of $M$; 
\item the charge $c$ encodes a pair of {\it bulk} and {\it boundary weights}, 
$h_c\in H^1(M;\Z/2\Z)$ and $k_c\in H^1(\partial M; \Z)$;
\item the pair $(w,f)$, and equivalently the system of log-branches, encodes another pair 
of {\it bulk} 
and {\it boundary weights}, 
$h_f\in H^1(M;\Z/2\Z)$ and $k_f\in H^1(\partial M;\C)$.
\end{itemize}
The weights $h_c$, $k_c$ (resp. $h_f$, $k_f$) satisfy certain natural
compatibility conditions.  
\smallskip

Let us remind a few details about the actual computation of these
weights from $(w,f,c)$. Concerning the charge $c$, we did it in the proof of Proposition \ref{Z2cmeaning}; with the same notations, simply set $k_c=\gamma(c)$ and $h_c= \gamma_2'(c)$. Given
$(w,f)$, the construction of $(k_f,h_f)$ is similar. Namely, as in the proof
of Proposition \ref{Z2cmeaning}: represent any non zero class in
$H_1(\partial M;\mz)$ by {\it normal loops}, say $C$, with respect to
$\partial T$; if the triangle $F$ of $\partial T$ is a cusp section of
the tetrahedron $\Delta$ of $T$, for every vertex $v$ of $F$ we denote
by $E_v $ the edge of $\Delta$ containing $v$, and define the index
$ind(C,v)$.  Moreover, we denote by $*_v$ the branching sign of
$\Delta$. Finally one defines the cohomology class $k_f \in
H^1(\partial M;\mc)$ by setting
\begin{align}
  k_f([C]) & := \sum_{C \ra E_v} *_v \ ind(C,v) l(E_v) \label{weightwf}\\
  & = 
\sum_{C \ra E_v} *_v \ ind(C,v)\left(\log(w(E_v)) +\pi \sqrt{-1}f(E_v)\right) 
\notag
\end{align}
Hence, the formula is very similar to the one given for the charge in
Section \ref{Z2cmeaning}, simply here we take the signs
$*_b$ into account. The other class $h_f \in H^1(M;\mz/2\mz)$
is defined similarly, by using normal loops in $T$ and taking the sum
mod$(2)$ of the flattenings of the edges we
face along the loops. Denoting by $d_w\in
H^1(\partial M; \C/2i\pi \Z)$ the log of the linear part of
the restriction of $\rho(w)$ to $\pi_1(\partial M)$, for all $a\in
H_1(\partial M; \Z)$ we have
\begin{equation}\label{fconstraint2}
k_f(a)= d_w(a) \ {\rm mod}(i\pi)\ ,\
(k_f(a)-d_w(a))/i\pi = h_f(a) \ \ {\rm mod} (2) . 
\end{equation}
\smallskip

We will study $\alpha_N(\Tt)$ again via a holographic approach. Let
$(T,\tilde b,w,f,c)$ be a QH triangulation of $M$.  By Lemma \ref{orC}
we know that the symmetrization factors $\alpha_N(\Tt)$ can be
computed by using any weak branching $\tilde b$ compatible with the
underlying pre-branching $\omega$. So {\it we can normalize the choice
  of} $\tilde b$ by requiring that $*_b=1$ on every branched
tetrahedron $(\Delta,b)$ of $(T,\tilde b)$.

\begin{figure}[ht]
\begin{center}
 \includegraphics[width=10cm]{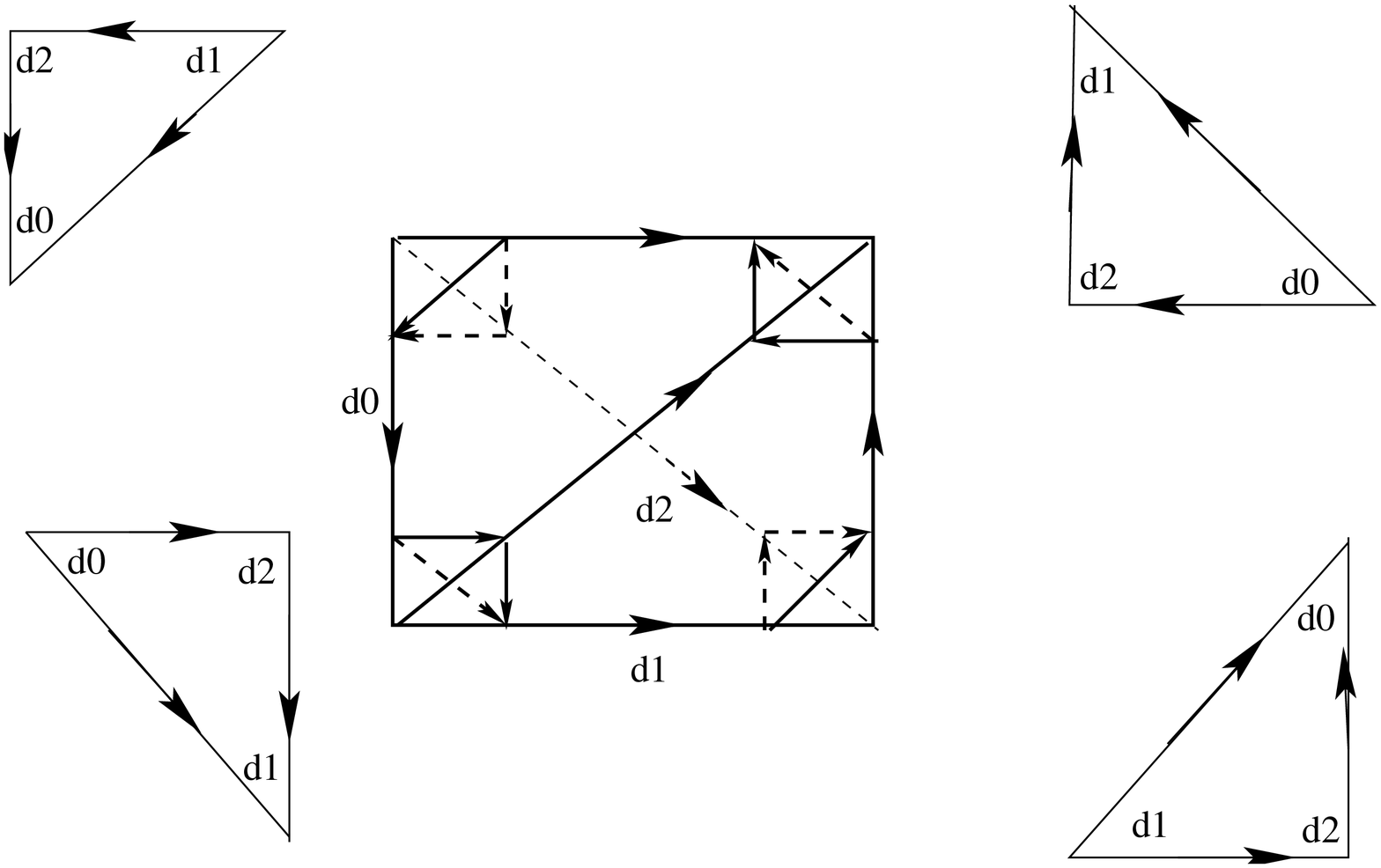}
\caption{\label{pb-boundary3} Boundary QH triangulations.} 
\end{center}
\end{figure}   

Referring to Figure \ref{pb-boundary3}, near each of the four vertices of
$(\Delta,b)$, the picture shows the corresponding branched triangle of
$(\partial T,\partial \omega)$, with the naturally induced ordered
decoration of the corners, denoted by $\partial d$. Note that we have
chosen one of the two positively branched tetrahedra that induce the
same pre-branching; it is immediate that the boundary corner
decoration {\it does not depend on this choice}.  So we have defined a
correspondence
$$ \Tt=(T,\tilde b,d)\rightarrow 
(\partial T,\partial \omega, \partial d) =: \partial \Tt . $$
For every
triangle $\tau:=(t,\partial \omega,\partial d )$ of $\partial \Tt$ set
$$\bw_k := \exp\left(\frac{\log(w_k) +\pi i (N+1)(f_k-c_k)}{N}\right), 
\ k =0,1,2 $$
$$\alpha^0_N(\tau)=\alpha^0_N(t,\partial \omega,\partial d ):= 
\bigl(\bw_0 ^{-c_1}\bw_1^{c_0}\bigr)^{\frac{N-1}{2}}$$
and define
\begin{equation}\label{bsym}
\alpha^0_N(\partial \Tt)= \prod_{\tau\in \partial \Tt}\alpha^0_N(\tau). 
\end{equation}
The following Lemma is immediate.
\begin{lem}\label{boundary_vs_bulk_alfa}
$\alpha^0_N(\partial \Tt)= \alpha_N(\Tt)^4$.
\end{lem}  
So, up to multiplication by $4$th-roots of $1$, the symmetrization
factor is completely determined by $\partial \Tt$.

\subsection{Flattening and charge invariance} 
We use the notations of Section \ref{boundaryQH}.

\begin{figure}[ht]
\begin{center}
 \includegraphics[width=9cm]{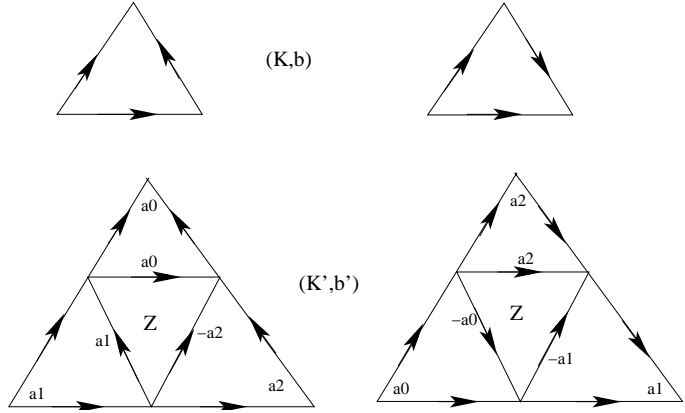}
\caption{\label{subdivide} The subcomplex $Z$.} 
\end{center}
\end{figure}

\begin{prop}\label{3D-weight-inv} Let $\Tt=(T,\tilde b,w,f,c)$ and
  $\Tt^*=(T,\tilde b,w,f,c^*)$ be two QH triangulations encoding a
  same tuple $(M,\rho,h,k)$, and differing only by the charges
  $c$ and $c^*$. Then $\alpha^0_N(\partial \Tt) = \alpha^0_N(\partial
  \Tt^*)$. Moreover, up to multiplication by a $N$-th root of $1$ the
  value of $\alpha^0_N(\partial \Tt)$ is the same for any flattening
  $f$ (hence also for any weight $(h_f,k_f)$ compatible with the holonomy $\rho$).
\end{prop}
\Dim The last claim is clear from the formula of $\alpha^0_N(\partial
\Tt)$. Let us consider the situation where $c^*\neq c$. As explained
in Section \ref{boundaryQH} we can assume that $(T,\tilde b)$ is
normalized by requiring that $*_b=1$ on every branched tetrahedron
$(\Delta,b)$. Let us subdivide $(K,b)=(\partial T, \partial \omega)$
by a branched triangulation $(K',b')$, as suggested in Figure
\ref{subdivide}. Denote by $(Z,b'_{Z})$ the subcomplex of $(K',b')$
formed by the union of the ``central'' triangles $\tau'$ of $K'$, lying
inside all triangles $\tau$ of $K$.  By confusing a complex with its
support, we have an inclusion $i: Z \to \partial M$, which induces a map
$i^*:H^1(\partial M;\C) \to H^1(Z;\C)$.

Set $a:=c-c^*$, that is, label each abstract edge of $T$ (hence each
corner of each triangle of $\partial T$) with the difference of its
values by $c$ and $c^*$. Define a $1$-cochain $\gamma_a$ on $Z$ as
indicated in Figure \ref{subdivide}, where the ordering of the labels
$a_j$ follows from the one shown in Figure \ref{pb-boundary3}. Actually
$\gamma_a$ is a $1$-cocycle.  Define another $1$-cochain $\gamma_l$ on
$Z$ by replacing $a_j$ with the log branch $l_j$ for every $j$.  By
restricting $\rho$ to $\partial M$ we get a
conjugacy class of representations $\rho_w: \pi_1(\partial
M)\cong\Z\times \Z \rightarrow {\rm Aff}(\C)$, where ${\rm Aff}(\C)$
is the group of complex affine transformations of the plane. The
linear part of $\rho_w$ defines a class in $H^1(\partial M;\C^*)$. Let
$\delta_w\in H^1(\partial M;\C/2i\pi \Z)$ be the $\log$ of this class,
with imaginary part in $]-\pi,\pi]$. The boundary weight $k_f\in
H^1(\partial M;\C)$ is a lift of the class $\delta_w$ in the sense
that, for every $a\in H_1(\partial M;\Z)$, we have (see \cite{AGT},
Section 4C)
$$k_f(a)=\delta_w(a)\quad {\rm mod}(i\pi).$$
Then we realize that
$$[\gamma_l]= i^*(k_f)\ , \ [\gamma_a]= i^*(k_c-k_{c^*}) . $$ 
Define the $2$-chain
$$c_Z=\sum_{\tau'\in Z} *_\tau \tau'. $$
Using the factorization formula of $\alpha_N(\Delta,b,w,f,c)$ at the
beginning of the proof of Proposition \ref{NA-invariance}, we see that
$$\alpha_N^0(\partial \Tt)\alpha_N^0(\partial \Tt^*)^{-1} =
\exp\biggl(2\frac{(N-1)}{N} (\gamma_l\cup \gamma_a - \gamma_a\cup
\gamma_l)(c_Z)\biggr).$$ Assume that $k_c=k_{c^*}$; then there exists
a $0$-cochain $\mu$ on $Z$ such that $\gamma_a = \delta \mu$. By the
formula of the cup product of cochains and the fact that $\gamma_l$ is
a cocycle we deduce
\begin{equation}\label{ratioSD}
  \alpha_N^0(\partial\Tt)\alpha_N^0(\partial\Tt^*)^{-1}  = 
\exp\biggl(2\frac{(N-1)}{N} (-\gamma_l\cup \mu
  - \mu \cup \gamma_l)(\partial c_Z)\biggr).
\end{equation}
Now, $(\gamma_l\cup \mu + \mu \cup \gamma_l)(\partial c_Z)$ is the sum
of the scalars $\mu(v)\gamma_l(e)$, where $v$ ranges over the vertices
of $Z$ and $e$ is an edge of $Z$ having $v$ among its endpoints (so
there are four such scalars for each $v$).  Here is a complete
description of $\mu$. Since $h_c=h_{c^*}$ and $k_c = k_{c^*}$, one can
obtain the charge $c$ from $c^*$ by local modifications taking place
in the abstract stars of some edges of $T$ (see the proof of
Proposition \ref{Z2cmeaning}). As viewed from $\partial M$, such a
modification can be described equivalently on the abstract stars in
$\partial T$ of the endpoints $x$ of the edge, in terms of $1$-{\it
  cochains} on $Z$. Namely, consider the cochain $\gamma_{c^*}$
encoding $k_{c^*}$, defined similarly as $\gamma_a$ above, by
replacing $a_j$ by $c_j^*$ for every $j$. Then a
modification about $x$ changes $\gamma_{c^*}$ by the coboundary of a
$\mz$-valued $0$-cochain $\mu_x$, supported by the vertices of $Z$
``around'' $x$, that is, the endpoints of the edges of $Z$ ``facing"
$x$ in the triangles $\tau\in {\rm Star}_{\partial T}(x)$. By
definition, $\mu_x$ takes the same value on all such vertices. We have
$\textstyle \mu=\sum_x \mu_x$, where $x$ ranges over the endpoints of
the edges about which $c^*$ is modified.

For instance, if the bottom left triangle of Figure \ref{subdivide}
represents $\tau$, and $x$ is its bottom left vertex, then $\mu_x$ is
supported in $\tau$ by the endpoints of the segment labelled by $a_1$;
if $\mu_x$ equals $n\in \mz$ at these vertices, then its coboundary
changes $\gamma_{c^*}$ by substracting $n$ on the edges
labelled by $a_0$ and $-a_2$.

One readily checks that the contribution of $(\gamma_l\cup \mu + \mu \cup \gamma_l)(\partial c_Z)$
coming from the edges of $Z$ in the star of $x$ is a multiple of $2\textstyle \sum_e
\gamma_l(e)$, where $e$ spans the set of edges of $Z$ facing $x$. By the log-branch condition about the edges of $T$, this
is $0$. Since these contributions are disjoint (there is no common
contribution), this proves $(\gamma_l\cup \mu + \mu \cup
\gamma_l)(\partial c_Z)=0$.\cvd
\medskip

We stress that the proof of Proposition \ref{3D-weight-inv} is {\it
  not} of cohomological nature (the $2$-chain $c_Z$ is not a
cycle). We have worked at the chain-cochain level, and this reflects
the fact that the symmetry defects (and the reduced QHI) are genuinely
geometric invariants.

\begin{remarks}\label{imp} {\rm The definition of $ \alpha^0_N(\tau)$
    depends only on the boundary trace $\partial \Tt$ of the QH
    triangulation $\Tt$, so that its global $3D$ entanglement is
    immaterial. This suggests that the functions $\alpha^0_N$ might
    belong to a certain ``free'' $2D$ QH theory.  However, though they
    verify a statement analogous to (1) of Theorem \ref{Intro-inv},
    there is no clean statement analogous to (2) because there are
    $b$-sliding transits between 2D QH triangulations for which
    $\alpha_N^{0}$ is not invariant. It is the entanglement of the
    system of $2D$ transits induced by $3D$ transits which ensures the
    invariance of $\alpha^0_N$ in that case (by Proposition
    \ref{NA-invariance}).  A quite complicated ``free'' $2D$ theory
    for $\alpha_N^0$ can nevertheless be developed on some quotient
    set of $2D$ QH triangulations.}
\end{remarks}

\subsection{Proofs of Theorem \ref {Intro-inv} and \ref {Intro-invbis}} 
 Theorem \ref {Intro-inv} (1) follows directly from
Lemma \ref{orC}, Lemma \ref{boundary_vs_bulk_alfa} and Proposition
\ref{3D-weight-inv}, and (2) follows from Proposition
\ref{NA-invariance}.  The proof of the first claim of (3) is then
immediate, by using the factorization formula \eqref{facto} in the
Introduction. As for the second claim, note that the reduced
state sums are just functions of the quantum shape parameters \eqref{qshapes}, and that these depend only on the system of {\it fused} flat/charges $f_k -*_bc_k$
mod$(N)$ (the only independent contributions of the charges in the
unreduced state sums being concentrated in the symmetrization
factors). On another hand, any two pairs of flattening weights
$(h_f,k_f)$ and $(h_f',k_f)$ compatible with $\rho$ and differing only
by the bulk weights $h_f$ and $h_f'$ are realized on every tetrahedron
by flattenings $f_k$ and $f_k':=f_k+a_k$, $a_k\in \mz$, such that the
weights associated to the system $a$ of integers $a_k$ are
$\gamma(a)=0\in H^1(\partial M;\mz)$ and $\gamma_2'(a) = h_f'-h_f\in
H^1(M;\mz/2\mz)$ (this follows from the arguments at the end of the proof of Proposition 4.8
(2) and 4.9 (2) in \cite{AGT}, that we have recalled briefly in Proposition \ref{Z2cmeaning} and \ref{3D-weight-inv}).  Since $N$ is odd, the system $Na$ has the same
weights, and hence also the flattenings $f_k+Na_k$ realize the weight $(h_f',k_f)$. However, by the previous observation the reduced
state sums take the same values on these flattenings and on the
flattenings $f_k$. Hence they do not depend on $h_f$. The same
argument works for the charge weights.  Finally, the formulas \eqref{weightc} and \eqref{weightwf} imply immediately that
the boundary weight associated to the fused flat/charges $f_k -*_bc_k$
is $\kappa:=k_f-\pi i k_c$. This concludes the proof of Theorem \ref
{Intro-inv}.

Consider now Theorem \ref {Intro-invbis}. If one applies to
$\alpha_N^0(\partial\Tt)$ the arguments of Proposition \ref{3D-weight-inv}, but fixing $c$ and changing the choice of log-branch system $l^*$ to $l$, we get the formula
\begin{equation}\label{ratioSD-2}
  \alpha_N^0(\partial\Tt)\alpha_N^0(\partial\Tt^*)^{-1}  = 
\exp\biggl(2\frac{(N-1)}{N} \left( (\nu\cup \gamma_c  + \gamma_c \cup \nu)(\partial c_Z)- (\nu \cup \delta\gamma_c + \delta\gamma_c\cup \nu)(c_Z)\right) \biggr)
\end{equation}
where $\nu$ is a $0$-cochain on $Z$
such that $\delta \nu=\gamma_{l-l^*}$. Then, by following the rest of the
proof with these substitutions we see that the contribution of $(\nu\cup \gamma_c +
\gamma_c \cup \nu)(\partial c_Z)$ on the edges $e$ of $Z$ around a
vertex $x$ is a multiple of $\textstyle \sum_e \gamma_c(e)=2\neq
0$. If only $f$ varies when changing $l^*$ to $l$, then $\nu$ takes values in $\pi i \mz$, and hence
$(\nu\cup \gamma_c + \gamma_c \cup \nu)(\partial c_Z)\in 2\pi i \mz$, like the other term in the exponential. This makes $\alpha_N^0(\partial\Tt)$ well
defined up to an $N$-th root of $1$. On another hand, if also $w$
varies, then $\nu$ takes values in $\mc$, and we have no control on
the variation of $\alpha_N^0(\partial\Tt)$. We correct this issue by
considering $\alpha_{N,c_0}^0(\partial\Tt):= \alpha_{N}^0(\partial\Tt)/\alpha_{N}^0(\partial\Tt_{c_0})$, instead of $\alpha_N^0(\partial \Tt)$. From \eqref{ratioSD-2} we deduce \begin{equation}\label{ratioSDN} 
  \alpha_{N,c_0}^0(\partial\Tt)\alpha_{N,c_0}^0(\partial\Tt^*)^{-1}  = 
  \exp\biggl(2\frac{(N-1)}{N} (\nu \cup \gamma_a + 
  \gamma_a\cup \nu)(\partial c_Z)\biggr) \end{equation}
where $a:=c-c_0$; note that $\gamma_a$ is a $1$-cocycle. Since $\textstyle \sum_e \gamma_a(e) =0$, by following the proof as 
before we get that $\alpha_{N,c_0}^0(\partial\Tt)$ is eventually a function 
of  $k_f$ and $k_c-k_{c}^0$. \cvd 
\medskip

\subsection{Proofs of Corollary \ref{corinv2} and \ref{corinv2bis}} 
\label{secondproofs} Consider the proof of Corollary \ref{corinv2}. Let $(T,\tilde
b,w,f,c)$ and $(T',\tilde b',w',f',c')$ be two QH triangulations
encoding $(M,\rho,h,k)$ and such that $(T,\omega_{\tilde b})$ and
$(T',\omega_{\tilde b'})$ represent the non-ambiguous structure
$\sG$. Let $(T,\tilde b)\rightarrow \ldots\rightarrow (T',\tilde b')$
be a sequence of non-ambiguous transits. By Theorem 6.8 of \cite{GT}
(in the case of cusped manifolds) and the results of Section 7 of
\cite{AGT} (that allow to deal with weak branchings $\tilde b$ which are not
global branchings), one derives a sequence of QH transits $(T,\tilde
b,w,f,c)\rightarrow \ldots\rightarrow (T',\tilde b',w',f'',c'')$,
where possibly $f''$ and $c''$ differ from $f'$ and $c'$. On another
hand, the system of shape parameters $w'$ is determined uniquely, as
follows (see the proof of Theorem 4.16 of \cite{AGT} for details). Suppose first that $M$ has a single cusp. Then one assumes that the gluing variety of $T$ (and similarly $T'$) has some {\it rich} irreducible component $Z$, which means that:
$Z$ contains a point $w_{hyp}$ representing the complete hyperbolic
structure of $M$; the regular map assigning its holonomy to a point $w\in Z$ is a birational isomorphism $\chi: Z \rightarrow X$, where $X$ is the geometric component of the variety of $PSL(2,\mc)$-characters of $M$ (for instance, a rich component exists if $T$ is a subdivision of the Epstein-Penner cellulation of $M$). Then, by taking $w \in \chi^{-1}(\rho)$, $\rho\in X$, the uniqueness of $w'$ follows from the fact
that the QH transits preserve the holonomy. If $M$ has several cusps, rich components map birationally to the {\it eigenvalue subvariety} $E$ of $X$ (see eg. \cite{KT} for the definition), and the conclusion is the same for $w \in \chi^{-1}(\rho)$, $\rho\in E$. Finally, by (2) and (1) of
Theorem \ref{Intro-inv} we have
$$\alpha_N(T,\tilde b,w,f,c) = \alpha_N(T',\tilde b',w',f'',c'')$$
and 
$$(\alpha_N(T',\tilde
b',w',f'',c''))^{4N} = (\alpha_N(T',\tilde b',w',f',c'))^{4N}.$$ This
proves the first claim of Corollary \ref{corinv2} (1), by noting that
Lemma \ref{boundary_vs_bulk_alfa} and Theorem \ref{Intro-inv} (1)
imply the invariance of $\alpha_N$ with respect to the choice of
$(h_c,h_f,k_f)$.  The second claim follows from it and the formula \eqref{facto} in the
introduction. 

Point (2) is merely a reformalization of (1). When $M$ has a single cusp, the covering space $\tilde
X_N$ is defined precisely in \cite{AGT}, Theorem 1.1. Briefly, by varying $\rho$ in $X$
and $\kappa$ among the fused weights compatible with $\rho$,
one describes a covering space of $X$ made of the sheets of the
(multivalued) map $X \rightarrow \mc^* \times \mc^*$ assigning to
$\rho$ a pair of $N$-th roots of eigenvalues of $\rho$ on basis elements of
$\pi_1(\partial M)$. These $N$-roots are rational expressions in the
quantum shape parameters $\bw$, $w\in \chi^{-1}(\rho)$. When $M$ has $p>1$ cusps, we get similarly a covering space of the eigenvalue variety $E$ of $X$, made of the sheets of the product map $X \rightarrow (\mc^* \times \mc^*)^h$ over the boundary components of $M$.

Corollary \ref{corinv2bis} follows from Proposition \ref{NA-invariance}, applied to the normalized symmetrization factors, and Theorem \ref{Intro-invbis}. Indeed, the latter shows that $\alpha_{N,c_0}^0(\partial\Tt)$ depends only on $k_f$ and $k_c-k_{c}^0$, these are preserved by QH transits, and any two QH triangulations encoding a same tuple $(M,k_f,k_c,\sG)$, with $M$ as in the statement, can be connected via sequences of non ambiguous QH transits. This last fact is proved as explained above, but without the need of the rigidity argument because only the weights contribute to $\alpha_{N,c_0}^0(\partial\Tt)$. Finally, the last claim is proved like (2) above, noting that $\alpha_{N,c_0}^0$ factors through the map $\{$QH triangulations of $M$ with fixed weight $k_c-k_c^0\} \rightarrow H^1(\partial M,\mc^*)$ given by $\Tt \mapsto \exp(4k_f/N)$.
\cvd

\begin{remark}{\rm There is an important difference between the
    invariance proof of the QHI and the proof of Corollary
    \ref{corinv2} (for the reduced QHI). In the latter one uses the 
invariance of the symmetry defect, for which a crucial role is played by
    {\it global} considerations like the charge invariance in
    Proposition \ref{3D-weight-inv}. On another hand for the QHI one
    uses the full transit invariance of the QH state sums. They may
    involve ambiguous transits, whereas by Proposition
    \ref{NA-invariance} these do not preserve the values of the
    symmetry defect, hence also of the reduced QHI.}
\end{remark}
   
\subsection{The case of pairs $(M,L)$}\label{MLdefect} Similarly to
Proposition \ref{NA-invariance}, for bubble transits we have (the
proof follows by very similar computations):

\begin{lem}\label{NA-bubble-invariance} 
  Let $(T,H,\tilde b,w,f,c)\rightarrow (T',H',\tilde b',w',f',c')$ be a
  QH bubble transit. If $(T,\tilde b)\rightarrow (T',\tilde b')$
  is non ambiguous, then $\alpha_N(T,\tilde b,w,f,c) =
  \alpha_N(T',\tilde b',w',f',c')$.
\end{lem} 

We show at first that the {\it non normalized} symmetry defects of
pairs $(M,L)$ are in general ill-defined. Consider a knot $L$ in
$S^3$, represented as the closure of a planar ``vertical'' braid
diagram $\Bb$ oriented from bottom to top. In \cite{LINK}, Lemma 3.2,
we described a procedure which associates to $\Bb$ a QH distinguished
triangulation $\Tt_\Bb$, suited to compute the QHI of $(S^3,L)$ (note
that $(\rho,h,k)$ is trivial in that case). The charge of $\Tt_\Bb$ is
the ''Yang-Baxter charge'' fixed in Section 3.3 of
\cite{LINK}. Moreover, every branched tetrahedron of $\Tt_\Bb$ carries
the same shape parameters and flattenings, given by
$(w_0,f_0,f_1)=(2,0,-1)$. If the knot diagram of $L$ is alternating,
then in Section 2.5 of \cite{LINK} we show that these shape parameters
and flattenings are the particular case at $a=-i$ of a family
$(w_0,f_0,f_1)_a$, depending on a generic complex parameter $a$, and
defining QH triangulations $\Tt_\Bb(a)$ that work as
well. 
Eventually we have $w_0 = 4a/(a+1)^2$ (so that $w_1=(a+1)^2/(a-1)^2$),
while $(f_0,f_1)$ is given by a so called {\it canonical
  log-branch}. For simplicity let us just compute the absolute value
$|\alpha_N(\Tt_\Bb(a))|$. This is a function of $|w_0|$ and $|w_1|$,
and we can forget the flattening. By a direct computation we get
$$|\alpha_N(\Tt_\Bb(a)|= \left( \frac{1}{\sqrt[N]{|w_0|^2}}
\left(\sqrt[N]{|\frac{w_1}{w_0}|}\right)^{2C}
\left(\sqrt[N]{|w_1|}\right)^{2s} \right)^{\frac{N-1}{2}}$$
where $C$ is the number of crossing of $\Bb$ and $s$ the number of
strands.  It is clear that $|\alpha_N(\Tt_\Bb(a)|$ varies with $a$.  The same examples also prove the claim (2) in Remark \ref{sQH}.
\smallskip

Next we prove Proposition \ref{invcaseML}. Consider a tuple
$(M,L,\rho,h;\sG)$; recall that $\rho$ is an {\it arbitrary} conjugacy
class of representations of $\pi_1(M)$ into $PSL(2,\mc)$, and
$h=(h_f,h_c)$ is an {\it arbitrary} pair of bulk weights on
$(M,L)$. Then, there exist QH distinguished triangulations $\Tt$
encoding $(M,L,\rho,h;\sG)$ \cite{Top,GT}; the only modification of
this notion with respect to the case of a cusped manifold $M$ is that
now $c$ is a $\mz$-charge of the pair $(M,L)$, as in Section
\ref{ML}. The proof of Theorem \ref{Intro-invbis} applies now to the
spherical components of $\partial M$. Then $k_f = k_c=k_c^0=0$. In
particular $\alpha_{N,c_0}(\Tt)^4$ does not depend on $c$, and so, by
taking $c=c_0$, we get $\alpha_{N,c_0}(\Tt)^4=1$. If $\sG$ is a taut
structure and $c=c_0$ is the charge tautologically carried by $\sG$,
this yields $\alpha_{N}(\Tt)^4=1$. The conclusion then follows from
\eqref{facto} (the independence on $h$ is proved like in Theorem
\ref{Intro-inv} (3)). \cvd

\section{Examples}\label{EXAMPLES}

We will use different ways to encode the weak branchings or the
pre-branchings, and in particular the so called $\Nn$-{\it graphs}
defined in \cite{AGT}. We refer also to \cite{AGT} for details about
gluing varieties. Concerning the symmetry defects, let us recall that,
by Corollary \ref{corinv3}, in the case of a taut structure $\sG$ they
are computed up to $4$-th roots of $1$ by the symmetrization factor of
any QH triangulation representing $\sG$. For non taut structures the symmetry defects are defined only up to multiplication by $4N$-th roots of $1$.

\subsection{A ``trivial'' example} For every $n\geq 1$, denote by
$F_n$ an oriented surface of genus $1$ with $n$ boundary components,
and by $\hat F_n$ the surface obtained by collapsing to one point each
boundary component of $F_n$. Hence $\hat F_n$ is a torus $F$ with a
set $P_n$ of $n$ marked points. The product manifold $M:= F_1 \times
S^1$ is the simplest example of a manifold with boundary a single
torus $S$, and carrying infinitely many different
fibrations. 

\smallskip

{\it Many fibrations.}  We consider a family
$\{M_n=M_{\psi_n}\}_{n\geq 1}$ of realizations of $M$ as mapping tori
of orientation preserving diffeomorphisms $\psi_n:F_n\to F_n$, which
extend to homeomorphisms $\hat \psi_n$ of $\hat F_n$ so that the set
$P_n$ is $\hat \psi_n$-invariant. They are defined as follows. We put $\psi_1 = {\rm id}$; 
so $M_1$ is the tautological product fibration of $M$. For every
$n\geq 1$ realize $\hat F_n$ as the quotient $ \R^2/(n\Z\times \Z)$; a
fundamental domain of the action is $R_n = [0,n]\times [0,1]$. The set $P_n$ is
given by the $n$ classes of the points in $R_n \cap \Z^2$. The
translation $\tau(x,y)=(x+1,y)$ of $\R^2$ induces the automorphism
$\hat \psi_n$ of $\hat F_n$.

\smallskip

{\it The layered triangulations.} For every $n\geq 1$, decompose the
domain $R_n$ in $n$ squares $Q_j = [j,j+1]\times [0,1]$, $j=0,\dots ,
n-1$. By means of the diagonal edge $[(j,1),(j+1,0)]$, every $Q_j$ is
triangulated by two triangles.  This induces a triangulation of the
domain $R_n$ and eventually an ideal triangulation $H_n$ of $\hat
F_n$.
Starting with $H_n$, we can
construct taut triangulations $\Tt_n$ of $\hat M$ which are
layered with respect to the mapping torus structure $M_n$.  As $H_n$
is $\hat \psi_n$-invariant, a ``full'' sequence of flips which
connects $H_n$ to $\hat \psi_n (H_n)$ can be formed by $3n$
couples of flips and inverse flips at every edge of $H_n$.  In order to specify
$\Tt_n$ it is then enough to fix an ordering $e_1,\dots e_{3n}$ of the
edges of $H_n$.  We fix such an ordering by lifting it to $R_n$,
accordingly to the following rules:
 \begin{itemize}
 \item In every square $Q_j$ we consider the ordered edges
   $[(j,1),(j+1,0)], \ [(j,0), (j+1,0)], \ [(j+1,0), (j+1,1)]$.
 \item These $3n$ edges of $R_n$ correspond bijectively to the $3n$
   edges of $H_n$. If $E$, $E'$ are among them, $E \subset Q_j$ and
   $E' \subset Q_i$ for $j<i$, then $E<E'$.
 \end{itemize} 
 
 For every edge $e_k$, $k=1,\dots , 3n$ of $H_n$, the associated
 couple of flip and inverse flip carries a polyhedron $B_k$ made of
 two pre-branched tetrahedra of $\Tt_n$, like the one resulting from a
 positive non ambiguous lune move.
 We denote by $\Delta_{\pm,k}$ the two tetrahedra of $B_k$, where the
 label ``$+$'' (resp. ``$-$'') indicates that the co-orientation of
 the two free $2$-faces goes in (resp. out) $B_k$. We fix a weak
 branching $\tilde b$ on $\Tt_n$ which is compatible with the
 pre-branching and such that the local branching $b_k$ induced on
 every tetrahedron has sign $*_b=1$. We use it to ``name'' the usual
 decorations on the tetrahedra; for instance
$$w(\pm,k)=( w(\pm,k)_0,w(\pm,k)_1,w(\pm,k)_2)$$ 
will denote the shape parameters on $(\Delta_{\pm,k},b_k)$, and so
on. Every truncated polyhedron $B_k$ carries $8$ triangles of the boundary
triangulation $\partial \Tt_n$.

We denote by $\delta(\pm,k)_r$, $r=0,1,2,3$, the triangle obtained by
truncation of the $r$th vertex of the branched tetrahedron
$(\Delta_{\pm,k},b_{\pm,k})$. We can assume that $\delta(+,k)_1$ and
$\delta(-,k)_2)$ (resp.  $\delta(+,k)_3$ and $\delta(-,k)_0$) have two
common edges and their union is like the result of a sliding $2D$
bubble move; $\delta(+,k)_0$ and $\delta(-,k)_1$ (resp.
$\delta(+,k)_2)$ and $\delta(-,k)_3$) have one common edge.

\begin{figure}[ht]
\begin{center}
 \includegraphics[width=3cm]{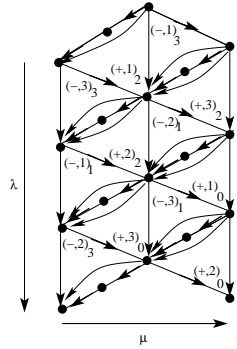}
\caption{\label{2_M1} The triangulation $\partial \Tt_1$.}
\end{center}
\end{figure}


\begin{figure}[ht]
\begin{center}
 \includegraphics[width=6cm]{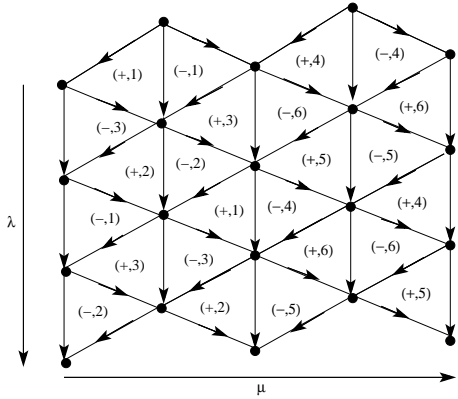}
\caption{\label{M2} The triangulation $\partial \Tt'_2$.}
\end{center}
\end{figure}

\smallskip

{\it The boundary triangulations.} In Figure \ref{2_M1}  we show a
planar fundamental domain $\Rr_1$ for $\partial \Tt_1$.
As usual, the opposite vertical boundary sides are identified, 
as well as the opposite horizontal boundary sides, in order to produce a
triangulated torus $S$. Let us denote by  $\partial \Tt'_1$ the triangulation 
obtained from $\partial \Tt_1$ by performing all evident negative $2D$ sliding 
bubble moves. 


In Figure \ref{M2} we show a fundamental domain $\Rr'_2$ for $\partial
\Tt'_2$, defined similarly. The rule to go back to a domain $\Rr_2$
relative to $\partial \Tt_2$ is clear.  In both figures, if a triangle
is labeled by $(\pm,k)$, it means that it has been obtained by
truncating a vertex of the tetrahedron $\Delta_{\pm,k}$. For the pairs
of triangles which form the bubbles this is understood.


On $\Rr_1$  we have also indicated some vertex indices (by omitting to write ``$\delta$''). 
On $\Rr'_2$ we omit the vertex indices because they are not essential in
the next discussion. 
$\Rr_1$ or $\Rr'_1$
plays the role of a ``basic tile'':  by forgetting the labels,
$\Rr'_2$ is formed by two copies of $\Rr'_1$ adjacent along a
vertical side, and in general $\Rr_n$ or $\Rr'_n$ is formed by $n$ copies
of $\Rr_1$ or $\Rr'_1$ adjacent along $n-1$ vertical sides. The rule
for the triangle labels can be naturally extended from $\Rr_2$ to an
arbitrary $\Rr_n$. In the figures we have also indicated a
``vertical'' generator $\lambda$ and a ``horizontal'' generator $\mu$
of the fundamental group of the torus $S$. They do not depend on
$n$. In particular $\lambda$ is isotopic to $\partial F_1
\subset \partial M$.  For every $n$ the {\it boundary combing} is the
one induced by the foliation by simple curves parallel to $\lambda$.

\smallskip

{\it One taut structure.} It is clear that all these mapping tori
$M_n$ belong to the same ray in $H_2(M,\partial M;\Z)$ and that the
layered triangulation constructed so far for $M_{nm}$ is also a layered
triangulation for both $M_n$ and $M_m$. Hence all these non equivalent 
fibrations of $M$ (they have non homeomorphic fibres) actually represent
a single taut structure on $M$. 

\smallskip

{\it The gluing varieties.} We encode characters $\rho$ of the mapping tori $M_n$ by points in the gluing varieties associated to the triangulations $\Tt_n$. The gluing variety $G(\Tt_1)$ is eventually defined by
the following system of independent equations:
$$ w(+,k)_{j+1} =1/(1-w(+,k)_j),  \ j=0,1,2 \ {\rm mod}(3)$$
$$  w(-,k)_2=w(+,k)_2^{-1}$$
where $k\in \{1,2,3\}$. So it is a non singular algebraic set of complex dimension
$3$, isomorphic to $(\C \setminus \{0,1\})^{3n}$.  The gluing variety
$G(\Tt_2)$ contains the set defined by the equations

$$ w(+,k)_{j+1} =1/(1-w(+,k)_j) , \ j=0,1,2 \ {\rm mod}(3)$$
$$  w(-,k)_2=w(+,k)_2^{-1}$$ 
$$w(\pm,k)_j=w(\pm,k+3)_j, \ j=0,1,2 \ {\rm mod}(3)$$
where $k\in \{1,2,3\}$. Again this set has complex dimension $3$. In
general, every gluing variety $G(\Tt_n)$ contains a similar
$3$-dimensional subset.  For every $n$, the data
$w(+,k)_2=\exp(i\pi/3)$, $k=1,\dots, 3n$, determine a point in such a
subset. Let us call it the {\it special shape system}. For every $n$
it encodes a fixed character $\rho_0$ of $\pi_1(M)$.
 \smallskip
 
 {\it On flattenings.} For every $n$, let us denote by $c_\tau$ the
 $\Z$-charge tautologically carried by the taut triangulation $\Tt_n$. It is easy to see that $f_\tau $ defined by
$$ f_\tau(\pm,k)= \mp c_\tau(\pm,k)$$  
is a flattening of the special shape system. The boundary weight
determined by the associated log-branches is trivially equal to zero
for every $n$.
 \smallskip
  
 {\it On charges.} Every charge on $\Tt_n$ is of the form $c=c_\tau +
 \gamma$, $\gamma$ being a solution of the {\it homogeneous} linear
 system associated to the defining equations of $\mz$-charges.  We fix a
 distinguished set of charges $c=c(a,b)$ by setting, for every $k=1,\dots, 3n$, 
$$ \gamma(+,k)_2= a,\ \gamma(+,k)_0= b, \gamma(+,k)_1= -(a+b)$$
$$ \gamma(-,k)_j=-\gamma(+,k)_j , \ j=0,1,2 $$
where $a$ and $b$ are free parameters. As for every $j\in \{0,1,2\}$
there is an equal number of terms $\gamma(+,.)_j$ and $\gamma(-,.)_j$,
the charge edge relations are automatically satisfied. The boundary
weight of $c(a,b)$ is given by
$$k_c(\lambda)=-6a, \ k_c(\mu)= -(2a+4b)n. $$

\smallskip

{\it On the symmetry defect.} Let us endow $\Tt_n$ with the special
shape system, the special flattening $f_\tau$, and a charge
$c(a,b)$. Then $\Tt_n$ becomes an ideal QH triangulation
$\Tt_n(a,b)$ of $\hat M$, which encodes some tuple $(M,\rho_0,h,k)$
whose $c$-weights vary with $(a,b)$.


Let us compute the symmetrization factors. For simplicity we omit the exponent $(N-1)/2$ overall, and we specialize the charges by setting $a=0$. Then the contribution of
every $B_k$ to $\alpha_N(\Tt_n(0,b))$ is equal to
\begin{multline}
\exp(\frac{i\pi}{3N} - \frac{i\pi(N+1)b}{N})^b
\exp(\frac{i\pi}{3N} + \frac{i\pi(N+1)b}{N})^b)\times
\\ \times \exp(\frac{-i\pi}{3N} + \frac{i\pi(N+1)b}{N})^{-b}
\exp(\frac{-i\pi}{3N} - \frac{i\pi(N+1)b}{N})^{-b}= \exp(\frac{4bi\pi}{3N}).
\end{multline}
Hence
$$\alpha_N(\Tt_n(0,b))= (\exp(\frac{4bi\pi}{3N}))^{3n}\ . $$ 
Given $\Tt_n$ and $\Tt_m$, $n\neq m$, then $\Tt_n(0,m)$ and $\Tt_m(0,n)$
encode the same boundary $c$-weight and eventually
$$\alpha_N(\Tt_n(0,m))=(\exp(\frac{4mi\pi}{3N}))^{3n}= (\exp(\frac{4ni\pi}{3N}))^{3m}= 
\alpha_N(\Tt_m(0,n))$$
as it must be! It is not hard to see that the same happens for every $\Tt_n(a,b)$ and $\Tt_m(a,b')$, for every choice of $a$, $b$ and $b'$ such that $(2a+b)n=(2a+b')m$.

\subsection{Three taut structures on the figure-eight knot complement}\label{eightsec}
Let $M_8$ be the complement of an open tubular neighbourhood of the
figure eight knot in $S^3$. It is a cusped manifold with a single
realization as a mapping torus $M_\psi$, with once punctured genus $1$
fibre. The minimal ideal triangulation $T$ of $\hat M_8$ has two
tetrahedra.  We are going to show that it carries three taut
structures represented by pre-branchings $(T,\omega_j)$, $j\in
\{0,1,2\}$, where $[(T,\omega_0)]=\sG_\psi$.  In doing it, we will
test the various invariants of non ambiguous structures that we have
introduced.

\begin{figure}[ht]
\begin{center}
 \includegraphics[width=9cm]{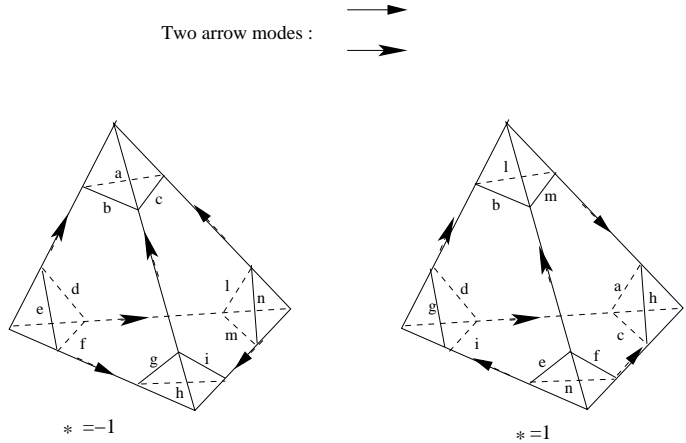}
\caption{\label{8-1} $(T,b)$ inducing $(T,\omega_1)$.}
\end{center}
\end{figure}

\begin{figure}[ht]
\begin{center}
 \includegraphics[width=7cm]{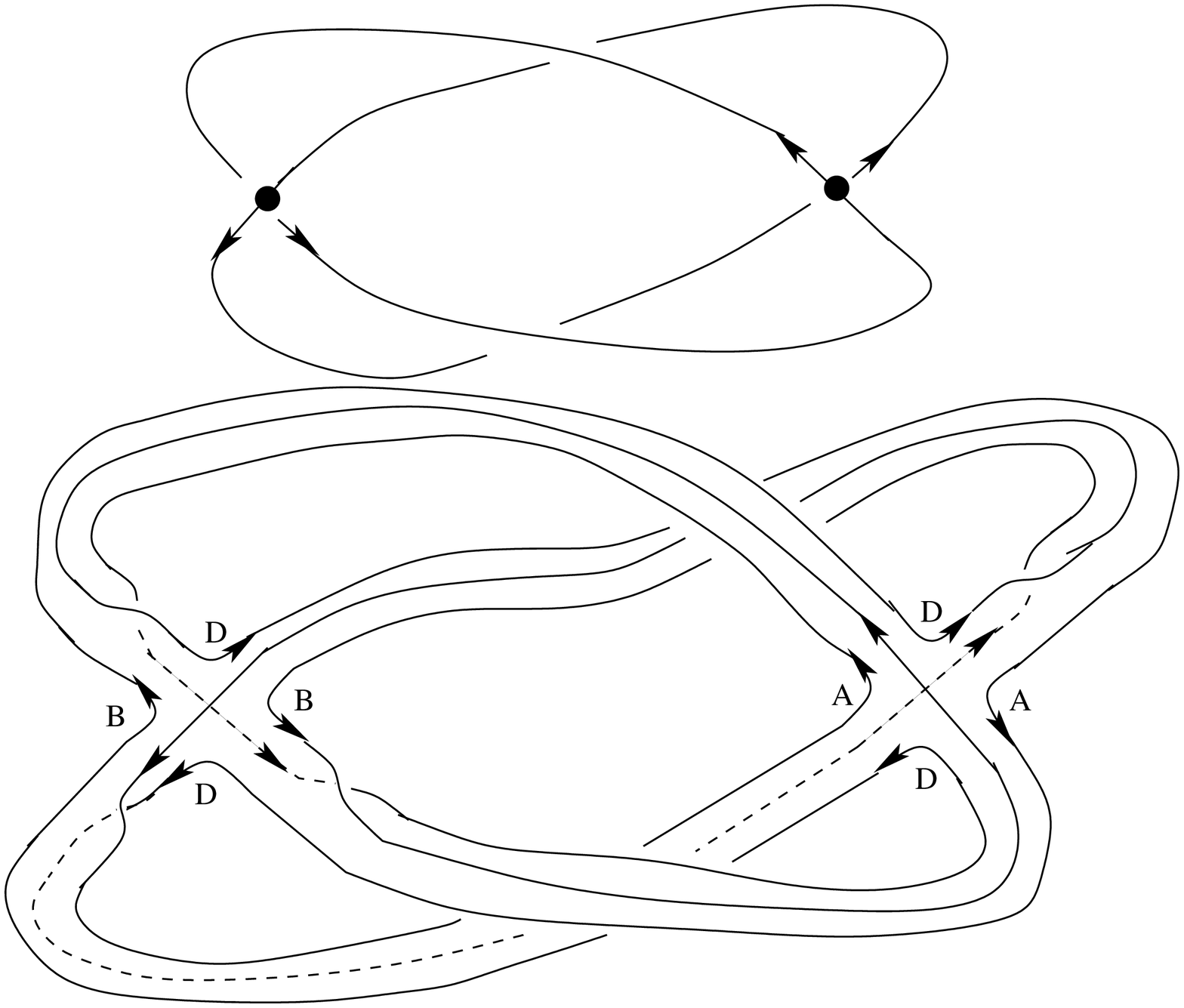}
\caption{\label{8-2} The $\Nn$-graph $\Gamma_1$.}
\end{center}
\end{figure}

Let us start with $(T,\omega_1)$. In Figure \ref{8-1} the system of
edge arrows specifies the $2$-face identifications that
produce $\hat M_8$, as well as a (genuine) branching $(T,b)$. Set
$(T,\omega_1):=(T,\omega_b)$.  The labels $a,b,\dots, m,n$ in the
figure show also the identifications of the sides of the triangles that form $\partial T$. From now on we will
refer to the tetrahedron on the right (resp. left) side as $\Delta_+$
(rest. $\Delta_-$). The various decorations of the edges of $\Delta_+$
(resp. $\Delta_-$) shall be denoted by capital (resp. small)
letters. The top picture in Figure \ref{8-2} shows an $\Nn$-graph
$\Gamma_1$ representing $(T,b)$; as $b$ is a branching, the edge
$\Z/3$-colors are equal to $0$, hence omitted. The bottom picture
shows the decoding of the graph in terms of a regular neighbourhood
$N(P)$ of ${\rm Sing}(P)$ in $P$, the internal spine of $M_8$ dual to
$(T,b)$. The triangulation $T$ has two edges as well as $N(P)$ has two
boundary components, each one being the boundary of an embedded disk
contained in $P\setminus {\rm Sing}(P)$, and dual to one edge of
$T$. At every vertex of $N(P)$ we see a distribution of the
pre-branching colors $A,B,D$; the rule depends on the sign $*_b=\pm
1$, and the color on the opposite branch of $N(P)$ at the crossing is understood. The
colors along every component of $\partial N(P)$ reflect the
distribution of the abstract edges of $(T,\omega_1)$. These facts hold in general for every $\Nn$-graph
representing a weakly branched triangulation $(T,\tilde b)$ (see
\cite{AGT}). Then one easily checks that $(T,\omega_1)$ is taut, and
moreover it verifies the property that the two diagonal (abstract)
edges at every edge of $T$ belong to a single tetrahedron.

\begin{figure}[ht]
\begin{center}
 \includegraphics[width=7cm]{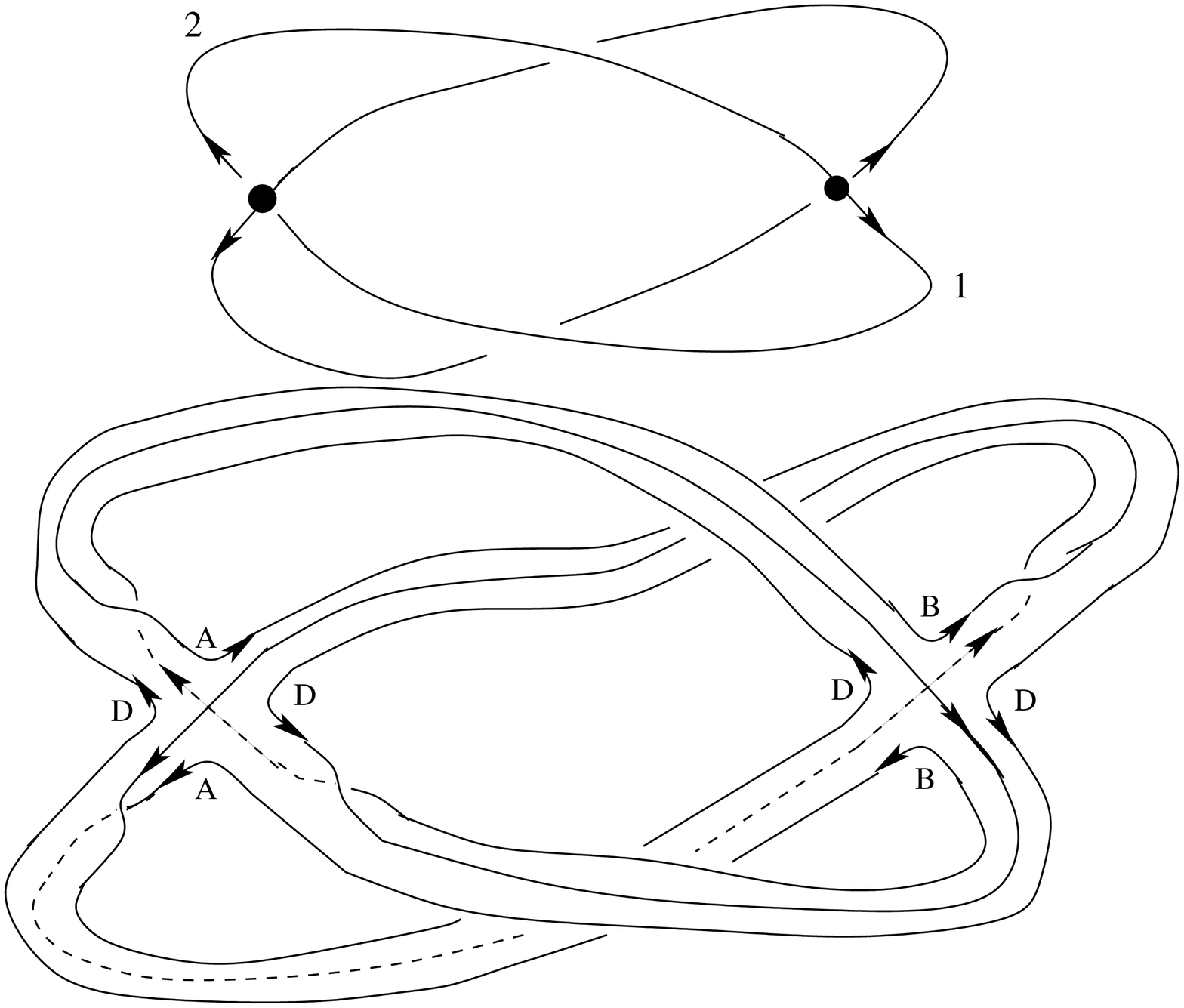}
\caption{\label{8-3} The $\Nn$-graph $\Gamma_0$.}
\end{center}
\end{figure}

Next we turn to $(T,\omega_0)$. Every pair of $\Nn$-graphs representing a given
triangulation $T$ are connected by a finite sequence of moves of two
types: either local moves at a dotted crossing that do not change the
underlying pre-branching; or inversion of the orientation of a circuit
of the graph together with a correction of the $\Z/3$-colors on the
edges that form the circuit (see \cite{AGT} and \cite{BP}). By just reversing
the orientation of one circuit in $\Gamma_1$, we get the $\Nn$-graph
$\Gamma_0$ of a weakly branched triangulation $(T,\tilde b)$; see
Figure \ref{8-3}. Set $(T,\omega_0):= (T,\omega_{\tilde b})$. As
$\tilde b$ is only a weak branching, there are non trivial
$\Z/3$-colors and orientation conflicts along some components of
$\partial N(P)$. We see that also $(T,\omega_0)$ is taut. Differently
from $(T,\omega_1)$, the two diagonal edges at each edge of
$(T,\omega_0)$ belong to different tetrahedra.  By looking at the
triangulations themselves, without passing to the dual picture, we see
that $\tilde b$ is obtained from $b$ by applying the transposition
$(0,1) $ on the (locally) ordered vertices $v_0,\dots , v_3$ of the
tetrahedron $\Delta_-$ of $(T,b)$, and the transposition $(2,3)$ on
$\Delta_+$.

Finally let us apply respectively the permutation $(0,2,1)$ on the
vertices of $\Delta_-$ and the permutation $(0,3,2)$ on the vertices
of $\Delta_+$ of $(T,b)$.  In this way we get a further
weak branching $(T,\tilde b')$. Set $(T,\omega_2):=
(T,\omega_{\tilde b'})$. It is also taut, and has the same properties as 
$(T,\omega_1)$ stated above.

\begin{figure}[ht]
\begin{center}
 \includegraphics[width=9cm]{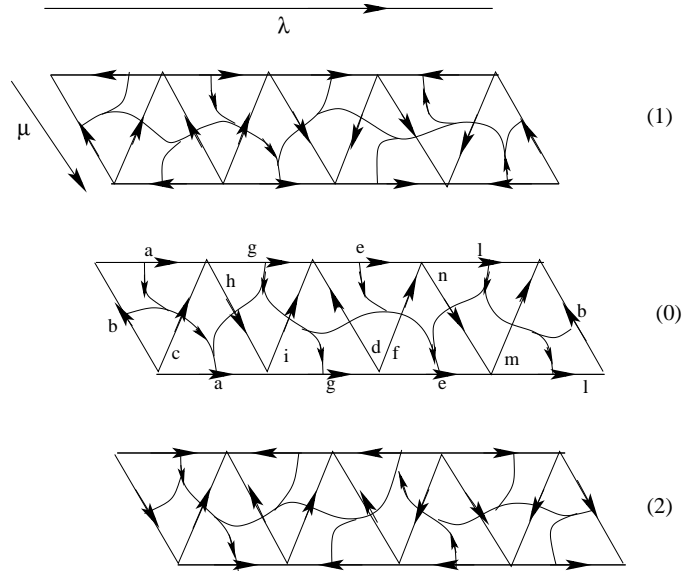}
\caption{\label{8-4} Boundary triangulations.}
\end{center}
\end{figure}  

In Figure \ref{8-4} we show (a fundamental parallelogram of) the
branched boundary triangulations $(\partial T, \partial \omega_j)$,
for $j\in\{1,0,2\}$, from top to bottom in that order. On $(\partial
T, \partial \omega_0)$ the edges are labelled by $a,b,c,\dots, m,n$,
referring to Figure \ref{8-1} (this is understood in the other two
cases).  The figure shows also the train track $\theta_j$ carried by
$(\partial T,\partial \omega_j)$, and a couple
$\lambda, \mu$ of generators of $H_1(\partial M_8;\Z)$ 
(the train track orientation is
shown only on a few arcs for future use).

\smallskip

{\it Computation of $\partial H^+(T,\omega_j)$.} Consider the truncated cell
decomposition of $M$ derived from $T$ as in Section \ref{2-cycles}, and the integer weights of the truncated $2$-faces, ie. hexagons. The integer weight of such a hexagon $H$ induces a labelling of the three edges of $\partial T$ contained in $\partial H$. So let us fix the following correspondence
``triples of triangle edges $\mapsto$ indeterminate integer weights":
$$ \{a,d,l\} \to x, \ \{c,i,n\} \to y, \ \{b,g,e\} \to z, 
\ \{h,f,m\}\to t  \ . $$
The relations between tuples $x,y,z,t \geq 0$ that define the non
  negative $1$-cycles on $(\partial T,\partial \omega_j)$ are
  respectively:
\begin{itemize}
\item $j=1$ : $z=-x, \ t= -y$ ;
\item $j=0$:  $z=x, \ t= y$;
\item $j=2$: $z=-x, \ t= -y$.
\end{itemize}
Then it is easy to verify that $\partial H^+(T,\omega_0)$ spans the
rank $1$ submodule of $H_1(\partial M_8;\Z)$ generated by $\lambda$,
while $\partial H^+(T,\omega_j)=\{0\}$ for $j\in \{1,2\}$. So, certainly
$[(T,\omega_0)]\neq [(T,\omega_j)]$ for $j\in\{1,2\}$, but for the moment the
last two are confused.

\smallskip

{\it The taut structure $[(T,\omega_0)]=\sG_\psi$.} It follows
immediately that the branched surface of $(T,\omega_0)$ carries a full
positive $2$-cycle.  Hence, accordingly with \cite{L}, it determines a
fibration which in this case is unique, and the full cycle is a
multiple of the fibre. The branched surface carries also a non
negative cycle with boundary $\lambda$ that realizes the fibre, hence
$\lambda$ is the canonical longitude.  Strictly speaking,
$(T,\omega_0)$ is not a layered triangulation of $M_\psi$ in the sense
of Section \ref{FIBER}, but it becomes one via a non ambiguous lune
move. So $[(T,\omega_0)]=\sG_\psi$.


\begin{remark}\label{othersurface}{\rm As $\partial
    H^+(T,\omega_j)=\{0\}$ for $j\in \{1,2\}$, these two taut structures on
    $M_8$ cannot be obtained as in Theorem \ref{existNA} (1) or
    (2). It is well-known that $M_8$ contains another properly
    embedded orientable incompressible and $\partial$-incompressible
    surface, which has the boundary slope $\gamma=\pm 4\lambda + \mu$,
    where the sign depends on the meridian orientation (see Chapter 4
    of Thurston's notes). This surface satisfies the hypothesis of
    \cite{L}, Theorem 2.  So it can be incorporated in a hierarchy
    $\Zz$ as in Theorem \ref{existNA} (2). Then $\sG_\psi \neq
    \sG_\Zz$ because $\gamma$ does not belong to $\partial
    H^+(T,\omega_0)$.}
\end{remark}
{\it Cohomological invariants.} Let $c^j$ be the canonical
$\Z$-charge carried by $(T,\omega_j)$. By using Figure \ref{8-4} it is
easy to compute the corresponding boundary weight $k_{c^j}$. We have:
\begin{itemize}
\item For every $j\in \{0,1,2\}$, $k_{c^j}(\mu)=0$.
\item $k_{c^0}(\lambda)=0$, $k_{c^1}(\lambda)= 2$ and $k_{c^2}(\lambda)=-2$.
\end{itemize}
Hence, the cohomological invariants $\partial\hG(*)$ distinguish the
three taut structures $[(T,\omega_j)]$, and it is also
excluded that one is obtained from another via the total inversion involution.
\smallskip

{\it Comparison of boundary combings.} The combing determined by the
boundary oriented train track of the branched surface of
$(T,\omega_j)$ is transverse and equivalent to the combing
$\theta_j$. So we can compare these combings by considering the oriented
curves carried by the $\theta_j$'s. By using the curves of this type
shown in Figure \ref{8-4}, we realize that the primary
obstruction $\sigma(\theta_0 -\theta_j)=\mu$, $j\in \{1,2\}$. Hence, like
for $H^+(*)$, the boundary combings distinguish $[(T,\omega_0)]$ from
$[(T,\omega_j)]$, $j=1,2$, but these last two are not distinguished.

 \smallskip
 
 {\it Symmetry defects.} Let $c^0$ be as above the $\Z$-charge of
 $(T,\omega_0)$.  We know that the corresponding boundary weight
 $k_{c^0}=0$, and that for any QH triangulation $\Tt_0=(T,\tilde b,
 w,f,c^0)$, $\alpha_N(\Tt_0)=1$.  With respect to the ordering of the
 charge entries determined by $(T,\tilde b)$, and with the ``capital vs
 small letters'' convention stated above, we have that $c^0$ is given
 by
 $$ (c_0,c_1,c_2)=(C_0,C_1,C_2)=(0,0,1). $$
 Let us place ourselves at the complete structure, so that the ordered
 shape parameters on $(T,\tilde b)$ start with:
 $$ w_0= \exp(i\pi/3), \ W_0=\exp(-i\pi/3). $$
 Let us fix the flattening with ordered entries with respect to
 $(T,\tilde b)$ given by:
 $$ (f_0,f_1,f_2)=(0,0,-1), \   (F_0,F_1,F_2)=(0,0,1). $$
 Let us rename these data with respect to $(T,b)$ and $(T,\tilde b')$
 respectively. We have:
 \begin{itemize}
 \item On $(T,b)$:
 $$ w_0 = \exp(-i\pi/3) , \ W_0=\exp(i\pi/3)$$
 $$ (c_0,c_1,c_2)=(0,1,0) , \ (C_0,C_1,C_2)=(0,1,0), (f_0,f_1,f_2)=(0,1,0) , \ (F_0,F_1,F_2)=(0,-1,0). $$
 
\item On $(T,\tilde b')$:
$$ w_0 = \exp(-i\pi/3)\ , \ W_0=\exp(i\pi/3)$$
$$ (c_0,c_1,c_2)=(1,0,0) , \ (C_0,C_1,C_2)=(0,1,0), (f_0,f_1,f_2)=(1,0,0) , 
\ (F_0,F_1,F_2)=(0,-1,0). $$
\end{itemize}
Let us denote by $\Tt_0$, $\Tt_1$ and $\Tt_2$ the resulting QH triangulations
supported by $(T,\tilde b)$, $(T,b)$ and $(T,\tilde b')$. 
They obviously represent a 
same tuple $(M_8,\rho_{hyp},h,k)$. By computation we get
$$ \alpha_N(\Tt_1)=1, \ \alpha_N(\Tt_2)= 
\exp(-\frac{2i\pi}{3N})^{\frac{N-1}{2}} \ . $$
These values separate $[(T,\omega_1)]$ from  
$[(T,\omega_2)]$ and $[(T,\omega_0)]$ from  $[(T,\omega_2)]$.

Let us consider now the canonical charge $c^1$. We specify the
corresponding QH triangulations $\Tt'_0$ and $\Tt'_1$ supported by
$(T,\tilde b)$ and $(T,b)$ as follows.
\begin{itemize}
\item $\Tt'_1$ is given by:
$$ w_0 = \exp(-i\pi/3)\ , \ W_0=\exp(i\pi/3)$$
 $$ (c_0,c_1,c_2)=(0,0,1) , \ (C_0,C_1,C_2)=(0,0,1), (f_0,f_1,f_2)=(0,1,0) , 
\ (F_0,F_1,F_2)=(0,-1,0). $$
 
 \item $\Tt'_0$ is given by:
 $$ w_0= \exp(i\pi/3)\ , \ W_0=\exp(-i\pi/3)  $$
 $$ (c_0,c_1,c_2)=(0,1,0) , \ (C_0,C_1,C_2)=(1,0,0), 
(f_0,f_1,f_2)=(0,0,-1) , \   (F_0,F_1,F_2)=(0,0,1). $$
\end{itemize}
Then we get
$$ \alpha_N(\Tt'_1)=1, \ \alpha_N(\Tt'_0)= 
\exp(-\frac{2i\pi}{3N})^{\frac{N-1}{2}}. $$ Again, this separates
$[(T,\omega_1)]$ from $[(T,\omega_0)]$. Moreover, we can conclude
that if $i \neq j$ then $[(T,\omega_j)]$ is not related to
$[(T,\omega_i)]$ via the total inversion involution.

\begin{figure}[ht]
\begin{center}
 \includegraphics[width=7cm]{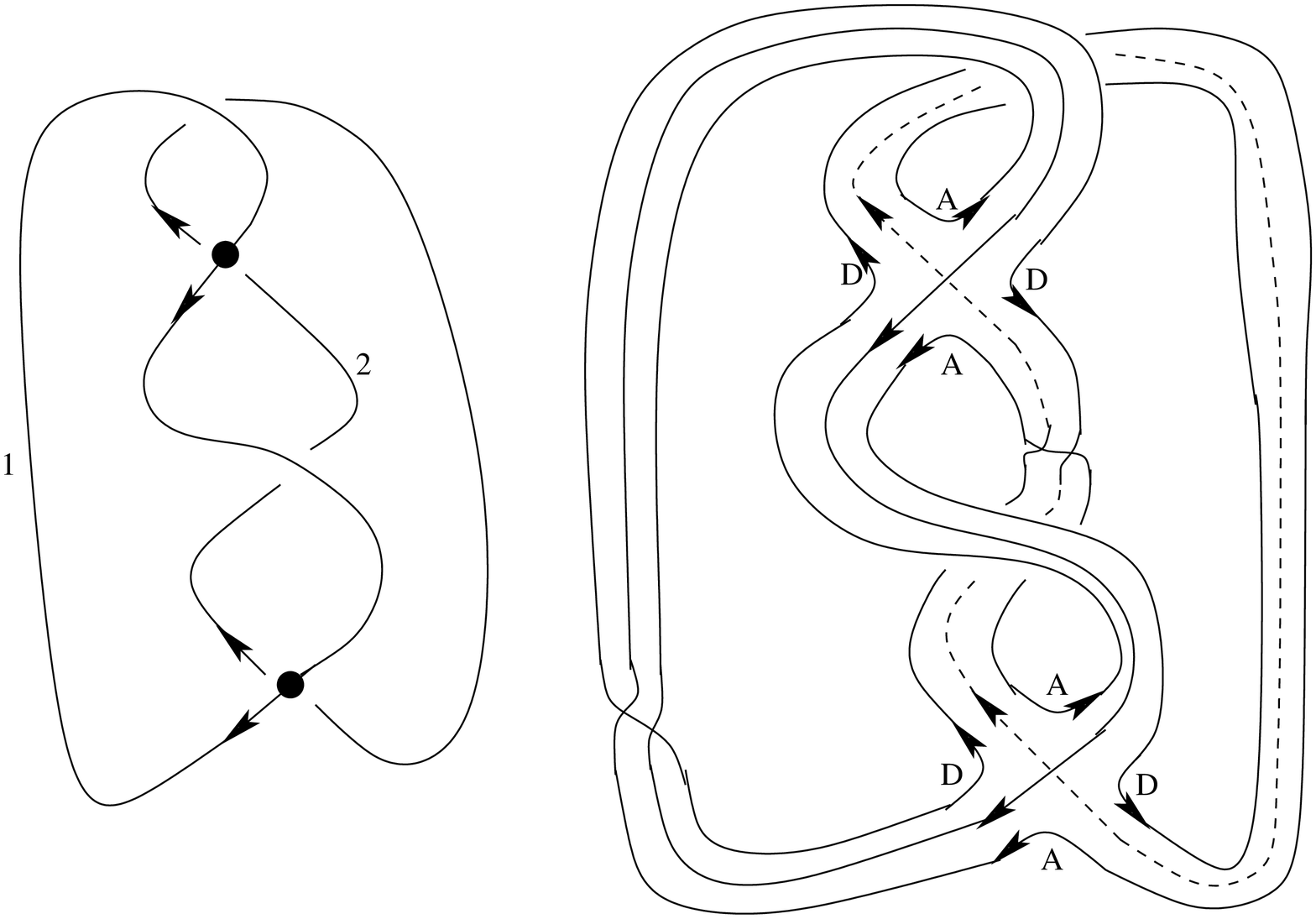}
\caption{\label{8S-1} The $\Nn$-graph $\Gamma_0$ 
for $M^s_8$ and its decoding.}
\end{center}
\end{figure}

\begin{figure}[ht]
\begin{center}
 \includegraphics[width=7cm]{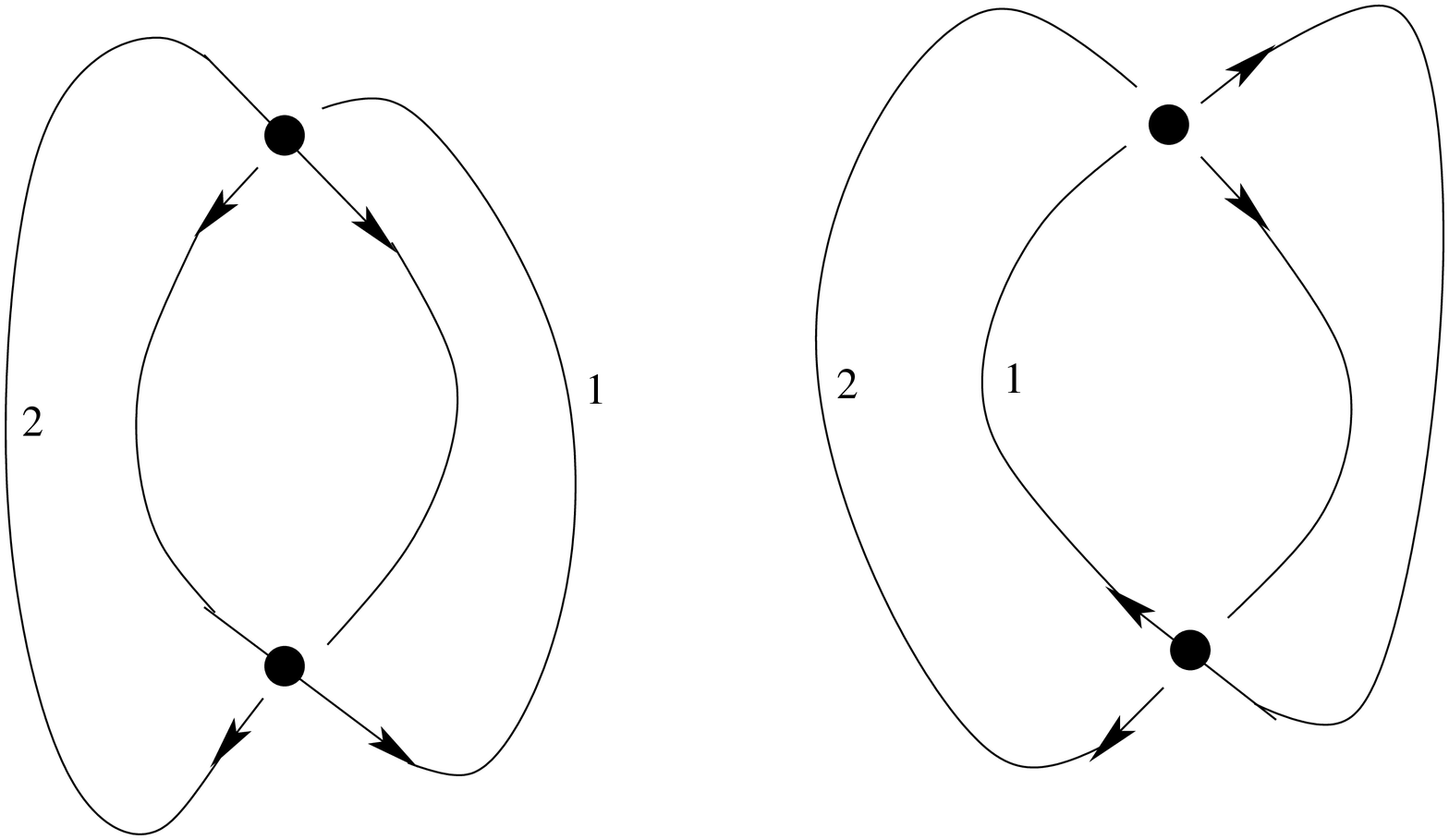}
\caption{\label{8S-2} Two more $\Nn$-graphs for $M^s_8$.}
\end{center}
\end{figure}  

\subsection{Three non ambiguous structures on the figure eight knot's
  sister}\label{sistersec}
Denote by $M^s_8$ this cusped manifold. In Figure \ref{8S-1} we show a 
$\Nn$-graph $\Gamma_0$ and its decoding, which represent a weakly branched
triangulation $(T,\tilde b_0)$ of $M^s_8$, supported by the minimal 
triangulation
$T$ with two tetrahedra. Note by the way that $T$ does not carry any
genuine branching.  Set $(T,\omega_0):= (T,\omega_{\tilde b_0})$. This
is the unique taut triangulation supported by $T$. One can see that
$H^+(T,\omega_0)=\{0\}$, so it does not represent the taut structure 
defined by the fibration of $M_8^s$.

In Figure \ref{8S-2} we show (left to right) two other $\Nn$-graphs
$\Gamma_j$, $j\in \{1,2\}$, representing weakly branched triangulations
$(T,\tilde b_j)$ of $M^s_8$. Set $(T,\omega_j)=(T,\omega_{\tilde
  b_j})$. They are not taut. Denote by $\Delta_+$ (resp. $\Delta_-$)
the tetrahedron of $T$ dual to the top (resp. bottom) dotted crossing
of the graphs.  We adopt the``capital vs small letters'' convention as
above.  Note that $(T,\tilde b_1)$ is obtained from $(T,\tilde b_0)$
by performing the permutation $(0,2,3)$ on the ordered vertices of
both $\Delta_+$ and $\Delta_-$. Also, $(T,\tilde b_2)$ is obtained by
performing the transposition $(2,3)$ on both tetrahedra $\Delta_+$ and
$\Delta_-$ of $(T,\tilde b_1)$.
  
 We are going to use the symmetry defects in order to distinguish the non 
 ambiguous structures represented by the pre-branched triangulations $(T,\omega_j)$. Let $c^0$ be the
  tautological $\Z$-charge on the taut triangulation $(T,\omega_0)$.
  As usual, for every QH triangulation $\Tt_0=(T,\tilde b_0,w,f,c^0)$ we have
  $\alpha_N(\Tt_0)=1$. Let us place ourselves at the complete structure
  on $M^s_8$. Consider the following specific QH triangulation $\Tt_1$
  supported by $(T,\tilde b_1)$ :
\begin{itemize}
\item $W_0=w_0= \exp(i\pi/3)$
\item $(C_0,C_1,C_2)=(c_0,c_1,c_2)=(1,0,0)$
\item  $(F_0,F_1,F_2)=(f_0,f_1,f_2)=(-1,0,0)$.
\end{itemize}
Clearly the charge is just $c^0$ renamed with respect to $(T,\tilde b_1)$.
Now let us rewrite these data with respect to $(T,\tilde b_2)$, getting
a QH triangulation $\Tt_2$:
 \begin{itemize}
\item $W_0=w_0= \exp(-i\pi/3)$
\item $(C_0,C_1,C_2)=(c_0,c_1,c_2)=(1,0,0)$
\item  $(F_0,F_1,F_2)=(f_0,f_1,f_2)=(1,0,0)$.
\end{itemize}

Finally 
$$ \alpha_N(\Tt_1)= \exp(-\frac{10i\pi}{3N})^{\frac{N-1}{2}}\ ,\ 
\ \alpha_N(\Tt_2)= \exp(\frac{10i\pi}{3N})^{\frac{N-1}{2}}.$$ Again
these values are different up to multiplication by $4N$-th roots of
$1$, so $\Tt_1$ and $\Tt_2$ represent distinct non ambiguous
structures.

\subsection{Six taut structures on the Whitehead link complement} 
\begin{figure}[ht]
\begin{center}
 \includegraphics[width=12cm]{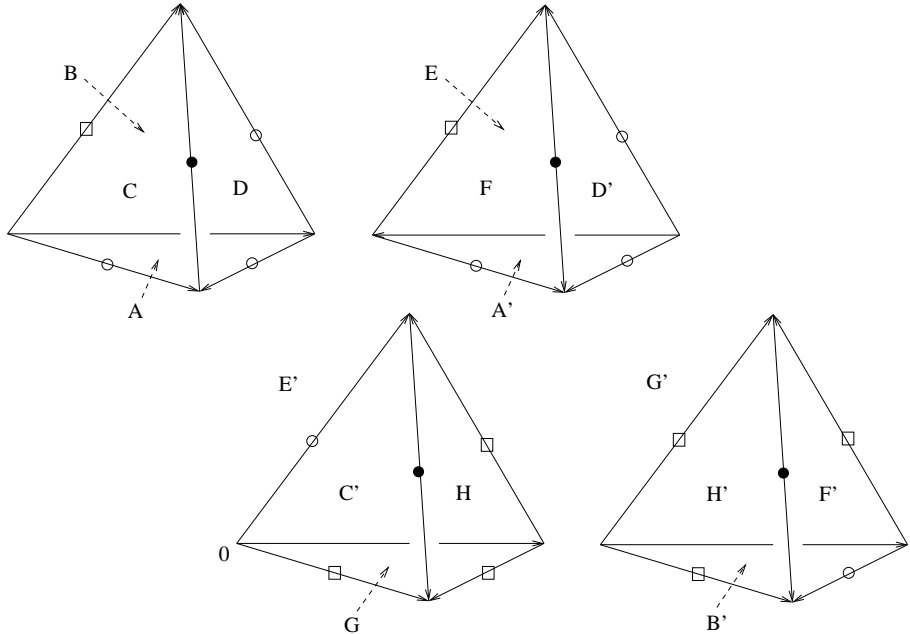}
\caption{\label{whiteheadtri} The triangulation $T$ of the Whitehead 
link complement.}
\end{center}
\end{figure}

Denote by $M$ the complement of an open tubular neighborhood of the
Whitehead link $L$ in $S^3$. We use the ideal triangulation $T$ of
$\hat M$ provided by SnapPea \cite{CDW} or Regina \cite{BBP}. It has
four tetrahedra $\Delta_0,\ldots,\Delta_3$ shown from left to right
and top to bottom in Figure \ref{whiteheadtri}; the face pairings are
indicated by the symbols and the orientations on the edges, and by the
letters on the faces.

\begin{figure}[ht]
\begin{center}
 \includegraphics[width=15cm]{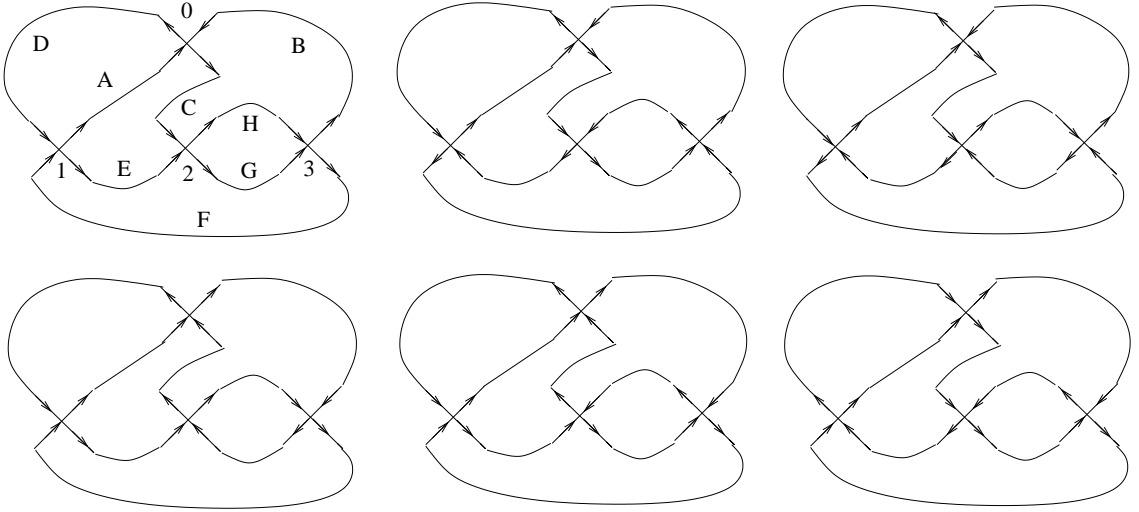}
\caption{\label{tautpreb} Graphs showing the six taut pre-branchings $\omega_j$ on $T$, up to total inversion.}
\end{center}
\end{figure}

By using eg. Regina, we see that $T$ carries $10$ taut angle
structures. Only six of them are compatible with some pre-branchings,
and for each one, there is only one compatible pre-branching up to the
the total inversion involution. These taut pre-branchings are depicted
in Figure \ref{tautpreb}; let us denote them $(T,\omega_j)$,
$j=0,\ldots,5$. For each graph, the vertices represent the tetrahedra
$\Delta_0,\ldots,\Delta_3$ and the edges represent the $2$-faces of
$T$, as labelled in the top left graph. So these graphs are just
obtained from usual $\Nn$-graphs by forgetting the informations at the
dotted vertices.

It is easy to determine the cyclically ordered sequences of $2$-faces
given by the pre-branching co-orientations about the edges of $T$. For
instance, the pre-branching $(T,\omega_0)$ gives (the numbers between
arrows correspond to the tetrahedra to which the faces belong to):
\begin{itemize}
\item sequence for the edge with no symbol: $1 \stackrel{E}{\longrightarrow} 2 \stackrel{G}{\longrightarrow} 3 \stackrel{B}{\longrightarrow} 0 \stackrel{A}{\longleftarrow} 1$; 
\item sequence for the edge with symbol $\bullet$: $0 \stackrel{C}{\longrightarrow} 2 \stackrel{H}{\longrightarrow} 3 \stackrel{F}{\longrightarrow} 1 \stackrel{D}{\longleftarrow} 0$;
\item sequence for the edge with symbol $\circ$: $3 \stackrel{F}{\longrightarrow} 1 \stackrel{A}{\longrightarrow} 0 \stackrel{D}{\longrightarrow} 1 \stackrel{E}{\longrightarrow} 2 \stackrel{C}{\longleftarrow} 0 \stackrel{A}{\longleftarrow} 1 \stackrel{D}{\longleftarrow} 0 \stackrel{B}{\longleftarrow} 3$;
\item sequence for the edge with symbol $\Box$: $2 \stackrel{G}{\longrightarrow} 3 \stackrel{F}{\longrightarrow} 1 \stackrel{E}{\longrightarrow} 2 \stackrel{H}{\longrightarrow} 3 \stackrel{G}{\longleftarrow} 2 \stackrel{C}{\longleftarrow} 0 \stackrel{B}{\longleftarrow} 3 \stackrel{H}{\longleftarrow} 2$.
\end{itemize}
In particular, the first two sequences show that the diagonal edges of
$\Delta_0$ are the common edges of the faces $A$, $B$ and $C$, $D$,
the last two show that the diagonal edges of $\Delta_3$ are the common
edges of the faces $B$, $F$ and $H$, $G$, and so on.

\smallskip

\begin{figure}[ht]
\begin{center}
 \includegraphics[width=15cm]{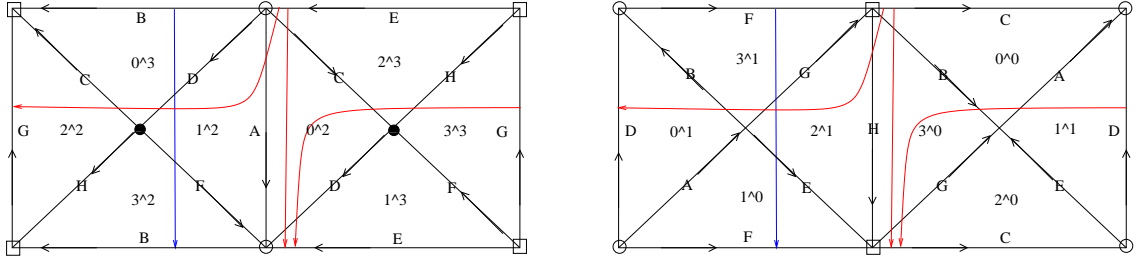}
\caption{\label{boundtri} The triangulations of the boundary components, with the branching $\partial \omega_0$.}
\end{center}
\end{figure}

{\it Computation of $H^+(T,\omega_j)$.} Denote by a small letter ``$h$" the integer weight of a hexagon $H$ of the cell decomposition of $M$ obtained from $T$, as in Section \ref{2-cycles}. The relations between the integer weights of the non-negative relative $2$-cycles on $(T,\omega_j)$ are respectively:
\begin{itemize}
\item $H^+(T,\omega_0)$: $a=e+g+b$, $d=c+h+f$, $b+c=f+e$;
\item $H^+(T,\omega_1)$: $b=c=d=e=f=h=0$, $a=g$;
\item $H^+(T,\omega_2)$: $b=c=a=e=f=g=0$, $h=d$;
\item $H^+(T,\omega_3)$: $b=c=a=e=f=g=0$, $h=d$;
\item $H^+(T,\omega_4)$: $b=c=d=e=f=h=0$, $a=g$;
\item $H^+(T,\omega_5)$: $g=a+b+e$, $h=c+d+f$, $b+f=c+e$.
\end{itemize}
For $j\in\{1,\ldots,4\}$, $H^+(T,\omega_j)$ is a cone of $H_2(M,\partial M;\mz)$ of rank $1$, generated by a surface of Euler characteristic $-1$ (by $a=g=1$ or $h=d=1$ it has $2$ triangles and $3$ edges); we will see that it has one boundary component, so it is a once punctured torus. For $j=0$ or $j=5$, $H^+(T,\omega_j)$ is more complicated. To identify the class in $H_2(M,\partial M;\mz)$ of a surface carried by $H^+(T,\omega_j)$, it is enough to look at the circuit defined by its boundary on the triangulations of the two boundary components of $M$ induced by $T$; these triangulations are shown in Figure \ref{boundtri}, with the meridians $\mu_0$ and $\mu_1$ (in red, $\mu_0$ being in the left picture) and the canonical longitude $\lambda_0$ and $\lambda_1$ (in blue, $\lambda_0$ being in the left picture). A label $i^j$ inside a triangle indicates that it is the boundary section of the $i$-th tetrahedron near its $j$-th vertex; the ordering of the vertices of tetrahedra is the same as for SnapPea: in Figure \ref{whiteheadtri}, for each tetrahedron the vertex $3$ is the top one, and the vertices $0$, $1$, and $2$ are in clockwise order starting from the bottom left vertex. The branching $\partial \omega_0$ on $\partial T$ is also shown. By the rule of Figure \ref{pb-boundary3} it is immediate to label the corners of each triangle with the decorations $d_j$. 
\smallskip

Doing similarly with each of the five other branchings $\partial \omega_j$, we find that the surfaces generating $H^+(T,\omega_j)$, $j\in \{1, 4\}$, are bounded by (a curve isotopic to) $-\lambda_1$, and those generating  $H^+(T,\omega_j)$, $j\in \{2, 3\}$, are bounded by $\lambda_0$. Let us denote them $b$ and $c$ respectively. A surface fully carried by $(T,\omega_0)$ is given by $b=c=e=g=h=1$, which gives $a=d=3$ and $f=1$. It is bounded by the curves $\mu_0+3\lambda_0$ and $-2\mu_1-2\lambda_1$, and so its class $a\in H_2(M,\partial M;\mz)\cong H_1(L;\mz)$ has the coordinates $(3,-2)$ in the basis $(\lambda_0, \lambda_1)$. Similarly, the class $d$ of the surface fully carried by $(T,\omega_0)$ given by $a=b=c=d=e=1$ has the coordinates $(-2,-3)$. 

\begin{figure}[ht]
\begin{center}
 \includegraphics[width=7cm]{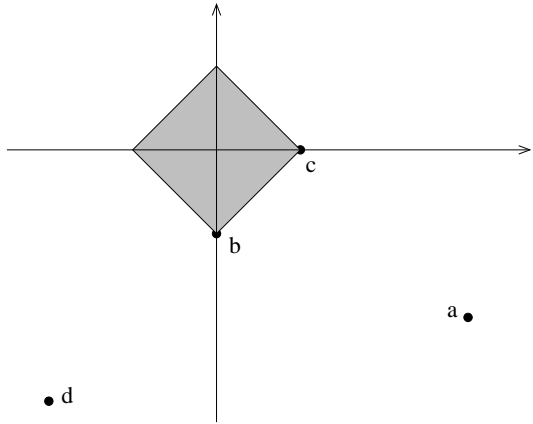}
\caption{\label{Thball} $B_M$ and a few classes of surfaces carried by one of the taut triangulations.}
\end{center}
\end{figure}

The classes $a$, $b$, $c$ and $d$ and the Thurston ball $B_M$ (determined in \cite{Th}) are shown in Figure \ref{Thball}. The classes $b$ and $c$ are vertices of the Thurston ball $B_M$ of $M$; $a$ and $d$ lie in the cones over two distinct and non opposite faces of $B_M$. The classes obtained from them by totally reversing the pre-branchings lie in the opposite vertices and fibered faces, respectively. 

\smallskip

{\it Cohomological invariants.} One obtains immediately the weight of the canonical charge $c^j$ of the taut triangulation $(T,\omega_j)$ by looking at Figure \ref{boundtri}, for $j=0$, and the analogous figures corresponding to the branchings $\partial \omega_j$, $j\in \{1,\ldots,5\}$ (the corners corresponding to diagonal edges, hence where $c^j=1$, have one incoming edge and one outgoing edge on the boundary of the triangle it belongs to). We get:
\begin{itemize}
\item $k_{c^0}$ is identically $0$;
\item $k_{c^1}(\mu_0) = k_{c^1}(\mu_1)=k_{c^1}(\lambda_1)=0$ and $k_{c^1}(\lambda_0)=2$;
\item $k_{c^2}(\mu_0) = k_{c^2}(\mu_1)=k_{c^2}(\lambda_0)=0$ and $k_{c^2}(\lambda_1)=2$;
\item $k_{c^3}(\mu_0) = k_{c^3}(\mu_1)=0$, $k_{c^3}(\lambda_0)=-1$ and $k_{c^3}(\lambda_1)=-2$;
\item $k_{c^4}(\mu_0) = k_{c^4}(\mu_1)=k_{c^4}(\lambda_1)=0$ and $k_{c^4}(\lambda_0)=-2$;
\item $k_{c^5}$ is identically $0$;
\end{itemize}
Hence the the invariants $\partial \mathfrak{h}(*)$ do not distinguish $(T,\omega_0)$ from $(T,\omega_5)$.
\smallskip

{\it The symmetry defects.} They are given on $(T,\omega_j)$ as functions of the charge $c$ and the system $\bw$ of quantum shape parameters. Assuming that we take compatible weak branchings $\tilde b_j$ normalized so that $*_b=1$ for each tetrahedron, the edge relations satisfied by $c$ and $\bw$ are  (first comes the equation associated to the edges $\bullet$ or with no symbol, and then comes the one for the edges $\circ$ or $\otimes$; see the relations before \eqref{qshapes} in Section \ref{SYMD})
$$\bw^2_0 \bw^0_2 \bw^1_2 \bw^3_0 = e^{-\frac{2i\pi}{N}}\quad , \quad (\bw^0_0)^2 \bw^1_1 \bw^2_2 \bw^3_2 (\bw^1_0)^2 \bw^0_1 = e^{-\frac{2i\pi}{N}}$$
and 
$$c^2_0 + c^0_2 + c^1_2 + c^3_0 = 2\quad , \quad 2c^0_0 + c^1_1 + c^2_2 + c^3_2 + 2 c^1_0 + c^0_1 = 2.$$
Besides these relations, there are also the tetrahedral ones, $c^k_0 + c^k_1 + c^k_2 =1$, $k\in \{0,\ldots,3\}$, and $\bw^k_0\bw^k_1\bw^k_2 = - e^{\pi/N}$. We see that the gluing variety of shape parameters has complex dimension $2$, and that the rank of the lattice of charges with given bulk and boundary weights $h_c$ and $k_c$ is equal to $2= 4\cdot (2 {\rm \ charges \ per\ tetrahedron})-(2\  {\rm edge \ relations})-(4\ {\rm weight \ relations})$. 

Let us denote by $\Tt_j = (T, \tilde b_j,w,f,c)$ a QH triangulation of $M$ supported by $T$ with induced pre-branching $\omega_j$. By using the tetrahedral relations and comparing the decorations $d_j$ as viewed from the different boundary branched triangulations $(\partial T,\partial \omega_j)$, we get easily (we put $\zeta:=\exp(2i\pi/N)$):
$$\alpha_N(\Tt_0)\alpha_N(\Tt_5)^{-1} = \left( \frac{\bw^2_1\bw^3_1}{\bw^0_0\bw^1_0} (-\zeta^\frac{N-1}{2})^{c^0_0+c^1_0-c^2_1-c^3_1} \right)^{\frac{N-1}{2}}$$
$$\alpha_N(\Tt_0)\alpha_N(\Tt_1)^{-1} = \left( \frac{\bw^1_1\bw^2_1}{\bw^3_0} (-\zeta^\frac{N-1}{2})^{c^3_0-c^1_1-c^2_1} \right)^{\frac{N-1}{2}}$$
$$\alpha_N(\Tt_0)\alpha_N(\Tt_2)^{-1} = \left( \frac{\bw^1_1\bw^3_1}{\bw^2_0} (-\zeta^\frac{N-1}{2})^{c^2_0-c^1_1-c^3_1} \right)^{\frac{N-1}{2}}$$
and so on. By specializing at the complete structure, for instance, where $w^2_0 = w^0_2 = w^1_2 = w^3_0 = i$, and choosing some charge and flattening, one can verify easily that the values are $\neq 1$, so that the symmetry defects distinguish the taut structures defined by the $\omega_j$'s.

\section{Appendix: The reduced Turaev-Viro
  invariants}\label{TV}
The construction of reduced Turaev-Viro (TV) invariants is much simpler than that of 
reduced QHI, because their arguments are just compact connected oriented $3$-manifolds $M$, possibly 
with non empty boundary. Indeed we slightly modify our set up
as follows:
\begin{itemize}
\item We consider only oriented compact $3$-manifolds $M$ with {\it
    non} empty boundary, and we stipulate that either $\partial M$ has
  no spherical component, or $\partial M$ is connected and consists of
  {\it one} spherical boundary component. The closed manifolds are
  covered by this latter case.

\item We use ideal triangulations of $\hat M$, hence we exclude the
  bubble move.

\item We deal with the TV invariants of manifolds with boundary
  already considered in \cite{BP1} (where they were related to the
  Witten-Reshetikhin-Turaev invariants of the ``double" $D(M)$), and
  recently reconsidered in \cite{CY}.
\end{itemize}

Let us recall the main features of these invariants.
Fix an integer {\it level} $r\geq 3$, and $q_0$ such that $q:=q_0^2$
is a primitive $r$th root of $1$. We have the set of half
integer colors
$$I_r=\{0,1/2,1,3/2,\cdots,(r-3)/2,(r-2)/2\} . $$ 
For every positive real number $x$, we denote by $x^{1/2}$ the positive
square root; then for $x<0$ we have $x^{1/2}=\sqrt{-1}|x|^{1/2}$. For every
integer $m\geq 1$, set
$$[m]=\frac{q_0^m-q_0^{-m}}{q_0-q_0^{-1}}\in \R \ ,
\ [m]!=[m][m-1]\cdots [1] .$$
For every $j\in I_r$, set
$$w_j=\sqrt{-1}^{2j}[2j+1]^{1/2}\ ,\ w=\sqrt{2r}/|q_0-q_0^{-1}| . $$

A {\it $r$-TV-tetrahedron} $(\Delta,b,\sigma)$ consists of: 
\begin{itemize}
\item A branched tetrahedron $(\Delta,b)$. 

\item An $I_r$-coloring $\sigma$ of the edges of $\Delta$.  We require
  furthermore that $\sigma$ is {\it admissible}, that is: 
\smallskip

\noindent {\it For every $2$-face $F$, the sum of the colors of the edges of $F$
  is $\leq r-2$, and the colors satisfy all triangular inequalities.}

\end{itemize}

As usual, the branching $b$ corresponds to an ordering $v_0,v_1,v_2,v_3$ of the
vertices of $\Delta$. If $F_k$ is the $2$-face opposite to the vertex
$v_k$, then the edges of $F_k$ are ordered as usual:
$e_{0,k}=[v_s,v_t]$, $e_{1,k}=[v_t,v_h]$, $e_{2,k}=[v_s,v_h]$, where
$s<t<h\in \{0,1,2,3\}\setminus \{k\}$. We will denote by
$\sigma_{i,j}$ the color that $\sigma$ gives to the edge $e_{i,j}$. 
\medskip

\begin{figure}[ht]
\begin{center}
 \includegraphics[width=8cm]{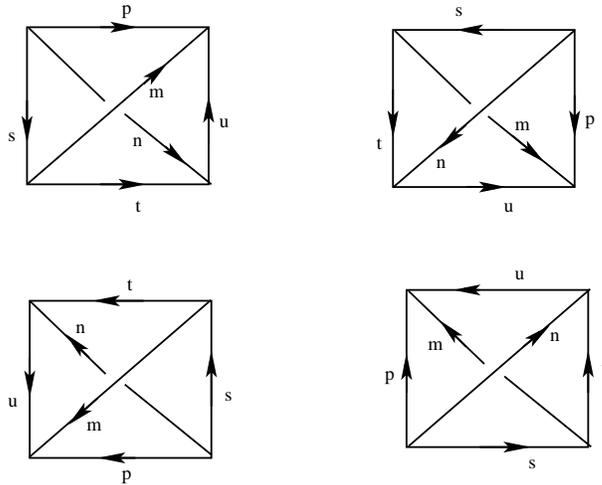}
\caption{\label{Icolors} $I$-labelling of $(\Delta,b)$.} 
\end{center}
\end{figure}

The {\it basic} $6j$-symbol $S_{q_0}(\Delta,b,\sigma)$ is a scalar denoted by
$$S_{q_0}(\Delta,b,\sigma) := \begin{Bmatrix} \sigma_{2,3} & \sigma_{1,2} &\sigma_{1,3} \\ 
\sigma_{0,1} &\sigma_{0,3}  &\sigma_{0,2} 
\end{Bmatrix}. $$ It is not important to give here the explicit
formula, which is derived from the representation theory of the
``small'' quantum group $\overline U_q(sl_2(\C))$ (see \cite{TV}). Figure
\ref{Icolors} shows four branched tetrahedra that share
a same pre-branching and carry the $6j$-symbol
$$\begin{Bmatrix} u & t & m \\ s & p  & n 
\end{Bmatrix}.$$ 

The {\it symmetrized} $6j$-symbol $Z_{q_0}(\Delta,b ,\sigma)$ is defined by: 
\begin{equation}\label{6j1}
Z_{q_0}(\Delta,b ,\sigma):= \begin{vmatrix} \sigma_{2,3} & \sigma_{1,2} &\sigma_{1,3} \\ 
\sigma_{0,1} &\sigma_{0,3}  &\sigma_{0,2} \end{vmatrix}:=
(w_{\sigma_{0,2}}w_{\sigma_{1,3}})^{-1}
\begin{Bmatrix} \sigma_{2,3} & \sigma_{1,2} &\sigma_{1,3} \\ 
\sigma_{0,1} &\sigma_{0,3}  &\sigma_{0,2}
\end{Bmatrix}.
\end{equation} 
Here, by ``symmetrized" we mean that the value of $Z_{q_0}(\Delta,b ,\sigma)$ is 
{\it the same for any branching $b$ on $\Delta$}.

\medskip

A $r$-{\it TV-triangulation} $\Tt=(T,\tilde b,\sigma)$ of $\hat M$
consists of a weakly branched triangulation $(T,\tilde b)$ equipped
with an $I_r$-coloring $\sigma$ of the edges of $T^{(1)}$, such that
every branched tetrahedron $(\Delta,b)$ of $(T,\tilde b)$ inherits a structure of
$r$-TV-tetrahedron. Each $\sigma$ is also called
a {\it $r$-state} of $(T,\tilde b)$. Clearly, for every $(T,\tilde
b)$ there is only a finite number of $r$-TV-triangulations
$(T,\tilde b,\sigma)$. Define the TV state sum and the {\it reduced}
TV state sum by
\begin{equation}\label{6j2}\textstyle 
|T,\tilde b|_{q_0}:= 
\sum_\sigma \left( \prod_{e\in T^{(1)}}w_{\sigma(e)}^2 \right)\  Z_{q_0}(T,\tilde b,\sigma)
\end{equation}
$$\textstyle |T,\tilde b|^{{\rm red}}_{q_0}:= 
\sum_\sigma S_{q_0}(T,\tilde b,\sigma)$$
where
$$ Z_{q_0}(T,\tilde b,\sigma)= \prod_{\Delta \in
  T^{(3)}}Z_{q_0}(\Delta,b,\sigma) \ , \ S_{q_0}(T,\tilde b,\sigma)=   
\prod_{\Delta \in T^{(3)}}S_{q_0}(\Delta,b,\sigma)$$ 
and $\sigma$ ranges over the set of $r$-states of $(T,\tilde b)$.
\begin{remark}{\rm There is an important difference between the TV and
the QHI local  symmetrization factors. In the TV case, they are the scalars
$ \alpha_{q_0}(\Delta,b,\sigma):= (w_{\sigma_{0,2}}w_{\sigma_{1,3}})^{-1}$ in \eqref{6j1}, 
which {\it depend on the $r$-state $\sigma$}. This and the 
factor $\textstyle \prod_{e\in T^{(1)}}w_{\sigma(e)}^2$ in \eqref{6j2} imply that it is no longer true
that the state sums $|T,\tilde b|_{q_0}$ differ from the reduced ones $|T,\tilde b|^{{\rm red}}_{q_0}$ just
by a global symmetrization factor.}
\end{remark}

The naked ideal triangulation moves can
be enhanced automatically to {\it TV-transits} $(T,\tilde
b,\sigma)\leftrightarrow (T',\tilde b',\sigma')$, defined by imposing that $\sigma'$ and 
$\sigma$ coincide on the common edges of $T$ and $T'$. The main result of
\cite{TV} is that $$|M|_{q_0} := |T,\tilde b|_{q_0}$$ is a well defined real
valued topological invariant of $M$. The proof is based on the fact
that the TV state sums are fully TV-transit invariant. 

\begin{remark}\label{q0}{\rm
The choice of $q_0$, and not only $r$, is far to be immaterial. For example, if
$r$ is odd, both $q_0=\exp(i\pi/r)$ and $q_0=\exp(2i\pi/r)$
are possible choices. Let $M$ have non trivial boundary, eg. take $M=M_8$ as 
in Section \ref{EXAMPLES}. In \cite{CY} one can find evidences that, for the first choice
of $q_0$, $|M|_{q_0}$ has polynomial growth when $r$ is odd and $r\to +\infty$, 
accordingly with the {\it Witten asymptotic expansion conjecture} for
the Witten-Reshetikhin-Turaev invariant $WRT_{q_0}(D(M))$ of the double of $M$ (see \cite{BP1}). On the other
hand, $|M|_{q_0}$ grows exponentially for the second choice of $q_0$.}
\end{remark}  
We can formulate the analog of Theorem \ref{Intro-inv} as follows.

\begin{teo}\label{TVinvariance} 
  (1) For every ideal triangulation $(T,\tilde b)$ of $\hat M$, the
  value of the reduced state sum $|T,\tilde b|_{q_0}^{{\rm red}}$ depends only the
  underlying pre-branched triangulation $(T,\omega_{\tilde b})$.  Moreover, it is
  invariant under the total inversion of the pre-branchings.

(2) If $(T,\tilde b)$ and $(T',\tilde b')$ are ideal triangulations
of $\hat M$ such that $(T,\omega_{\tilde b})$ and $(T',\omega_{\tilde b'})$
represent the same non ambiguous structure on $M$, say $\sG$, then
 $|T,\tilde b|_{q_0}^{{\rm red}} = |T',\tilde b'|_{q_0}^{{\rm red}}$.  
\end{teo}
Hence there are well defined {\it reduced} TV invariants
$|M;\sG|^{{\rm red}}_{q_0}$, which are invariant under the natural total
inversion involution on the set of non ambiguous structures.
\smallskip

\Dim A change of the weak branching that preserves or reverses totally
the pre-branching does not modify the TV local symmetrization factors,
because $(w_{\sigma_{0,2}}w_{\sigma_{1,3}})^{-1}$ does not depend on the colors of
the square edges. So Theorem \ref{TVinvariance} (1) follows
immediately.

As for the TV transits, let us just consider the $2\rightarrow 3$ one
(things are similar with the lune move). It is enough to prove that
\begin{equation}\label{TVSchaeffer}
S_{q_0}(T,b,\sigma) = \sum_{\sigma'} S_{q_0}(T',b',\sigma')
\end{equation}
whenever $(T,b)\to (T',b')$ is one of the two remarkable $b$-transits described 
before Lemma \ref{wb}, and the move $(T,b,\sigma)\rightarrow
(T',b',\sigma')$ varies among its TV enhancements.  We realize that \eqref{TVSchaeffer} 
is the Biedenharn-Elliot identity that arises from the representations theory
of $\overline U_q(sl_2(\C))$. This achieves also Theorem
\ref{TVinvariance} (2). \cvd \medskip

It is interesting to recover the proof of Theorem \ref{TVinvariance}
(2) in a way which points out also a {\it holographic content}, in the sense of what we have developed in Section \ref{2D}.  Every
triple $(T,\omega,\sigma)$ as above restricts naturally to a boundary
``TV branched triangulation'' $(\partial T,\partial \omega, \partial
\sigma)$, where $\partial \sigma$ is the $I_r$-coloring of the
vertices of $\partial T$ induced by $\sigma$. For every such a colored $2D$ branched
triangulation $(\partial T,\partial \omega, \partial
\sigma)$ and every vertex $v$ of $\partial T$, let $C_v$ be the $2$-cell dual to $v$. Define
$$\alpha_{TV}(\partial T,\partial \omega, \partial
\sigma)= \prod_{v\in \partial T} w_{\partial \sigma(v)}^{{\rm Eu}_b(C_v)}. $$
A simple rewriting of the formulas shows that:
\begin{lem} For every TV-triangulation
 $(T,\tilde b,\sigma)$ of $\hat M$ we have
 $$ Z_{q_0}(T,\tilde b,\sigma)=  
 \alpha_{TV}(\partial T, \partial \omega_{\tilde b}, 
 \partial \sigma) S_{q_0}(T,\tilde b,\sigma).$$
\end{lem}
Hence  $\alpha_{TV}(\partial T, \partial \omega_{\tilde b}, 
\partial \sigma) $ can be considered as the symmetrization
factor of the product of basic $6j$-symbols $S_{q_0}(T,\tilde b,\sigma)$, and is entirely determined by
the boundary TV triangulation. Now, in the remarkable $b$-transits described before Lemma \ref{wb}, for every $\sigma'$ we have
$$ \alpha_{TV}(\partial T',\partial \omega',\partial \sigma')=
\alpha_{TV}(\partial T,\partial \omega,\partial \sigma).$$ So,
up to an overall scalar factor, the Biedenharn-Elliot identity \eqref{TVSchaeffer} {\it
  coincides} with the identity obtained from it by replacing $S_{q_0}$ with $Z_{q_0}$. Ultimately, this
depends on the nice behaviour of the scalar $\alpha_{TV}(\partial T, \partial \omega, 
\partial \sigma)$
under $2D$ sliding transits.

\medskip

Finally, recall that the taut structures are in a sense the most natural non
ambiguous structures. To this respect we note that the reduced TV invariants are
disappointing. In fact $\alpha_{TV}(\partial T, \partial \omega, 
\partial \sigma) = 1$ for a taut triangulation $(T,\omega)$. Hence:
\begin{cor}\label{blindTV} For every taut structure $\sG$ on $M$ we have $|M;\sG|^{{\rm
      red}}_{q_0}=|M|_{q_0}$.
\end{cor} 
So the reduced TV invariants are completely blind with respect to the
taut structures (though they can distinguish other non ambiguous
structures on $M$). This can be regarded as a TV counterpart of the
fact that the QHI symmetrization factors $\alpha_N(T,\tilde b,w,f,c)=1$, whenever 
$(T,\omega_{\tilde b})$ is a taut triangulation and $c$ is
the tautological charge carried by $(T,\omega_{\tilde b})$. The advantage of the
QH framework is that we can vary also the $c$-weight $k_c$ in order to
distinguish different taut structures, as we have done in Section
\ref{EXAMPLES}. 

On the other hand, recall that the normalized QH
symmetry defects are blind with respect to the relative non ambiguous
structures on $(M,L)$. The situation for reduced TV invariants is
slightly better in the relative case.  Let us work for example with relative taut
triangulations $(T,H,\omega)$ of $(M,L)$ (if any) and the
corresponding relative taut structures $\sG$. Then the reduced TV
invariant $|M,L;\sG|^{\rm red}_{q_0}$ does not coincide in general with
$|M|_{q_0}$, and it should be sensitive to both $L$ and $\sG$. As an
example, let us outline a construction of invariants of fibred knots
$K$ in $S^3$. Let $Y$ be the complement of an open tubular
neighbourhood of $K$ in $S^3$. Let $M$ be obtained via Dehn filling of
$Y$ along the canonical longitude, say $m$, of $K$, and $L\subset M$
be the core of the attached solid torus. Since $Y$ has a single
fibration it has a canonical taut structure $\sG_Y$. By Section
\ref{ML} it extends to a relative taut structure $\sG_{K}$ on
$(M,L)$. Hence $|M,L;\sG_{K}|^{{\rm red}}_{q_0}$ is an invariant of
the fibred knot $K$.

\end{document}